\documentclass[a4paper,11pt]{article}
\usepackage{amssymb}
\usepackage{amsmath}
\usepackage{latexsym,a4wide}
\usepackage{hyperref}
\usepackage{tikz-cd,mathtools}
\usetikzlibrary{calc,arrows}
\tikzset{mydescription/.style={anchor=center,fill=white}}
\tikzcdset{%
    triple line/.code={\tikzset{%
        double equal sign distance, 
        double=\pgfkeysvalueof{/tikz/commutative diagrams/background color}}},
    quadruple line/.code={\tikzset{%
        double equal sign distance, 
        double=\pgfkeysvalueof{/tikz/commutative diagrams/background color}}},
    Rrightarrow/.code={\tikzcdset{triple line}\pgfsetarrows{tikzcd implies cap-tikzcd implies}},
    RRightarrow/.code={\tikzcdset{quadruple line}\pgfsetarrows{tikzcd implies cap-tikzcd implies}}
}
\newcommand*{\tarrow}[2][]{\arrow[Rrightarrow, #1]{#2}\arrow[dash, shorten >= 0.5pt, #1]{#2}}

\input xypic
\input xy
\xyoption{all}

\newtheorem{example}{Example}[section]
\newtheorem{defn}[example]{Definition}

\newtheorem{thm}[example]{Theorem}

\newtheorem{lem}[example]{Lemma}

\newtheorem{remark}[example]{Remark}
\newenvironment{pf}{\noindent \textbf{Proof:} }{$\Box$ \mbox{}}

\newcommand{\ke}{\mathrm{Ker}}
\begin{document}

\date{}
\author{Murat SARIKAYA and Erdal ULUALAN   }
\title{On Regular Representation of 2-Crossed Modules and Cat$^2$-Groups}
\maketitle

\begin{abstract}
In this paper, we describe a regular representation given by  Cayley's theorem for 2-crossed modules of groups and their associated Gray 3-(group) groupoids with a single 0-cell and equivalently cat$^2$-groups.


\textbf{AMS.Classification:} 18B40, 18G45, 20C99, 55U15, 55U35, 20L05.
\end{abstract}

\section*{Introduction}
The main aim in representation theory for groups is to find an appropriate group of linear transformations or permutations with the same structure as a given, abstract, group, see \cite{Ledermann}. A linear representation of a group $G$ is a homomorphism of groups $\phi$ from $G$ to $GL(V)$, where $GL(V)$ is the general linear group of linear isomorphisms $V\rightarrow V$ for a vector space $V$. Then, the representation $\phi$ assigns to any element $g$ of $G$ a linear isomorphism $\phi_g:V\rightarrow V$. A permutation representation of a group $G$ is related to actions of $G$ on a set $X$. An action of $G$ on $X$ gives a representation by mapping to each $g\in G$ the permutation $\phi(g)\in S_X$ that takes $x\in X$ to $g_x$. Conversely, if $\phi:G \rightarrow S_X$ is a representation, then $G$ acts on $X$ by $g_x:=\phi_g(x)$ where $\phi_g$ is a permutation of $X$. Linear matrix representations can be given as homomorphisms of $G$ into $GL_n(K)$ for a field $K$, giving invertible matrix for each element of $G$. Any group $G$ can be thought as a category with a single object ${*}$ and invertible morphisms, so a functor from $G$ to the category of sets corresponds to a representation of $G$.

The notions of cat$^1$-groups \cite{Loday}, \cite{WL} and crossed modules, \cite{W} are equivalent structures of 2-groups (cf. \cite{BL}) as 2-dimensional generalisations of groups. These structures are algebraic models for homotopy (connected) 2-types. This generalisation suggests that the linear and permutation representation of a group $G$ can be generalised to 2-dimensional versions of groups. Thus, representations of crossed modules or equivalently cat$^1$-groups should be 2-functors into a suitable category. To find this suitable category with the same structure as given, abstract, cat$^1$-groups, Barker, \cite{Barker} has investigated the 2-category $\mathbf{Ch}^1_K$ of length-1 chain complexes of vector spaces over a fixed field $K$. This category, $\mathbf{Ch}^1_K$ has a 2-groupoid structure. Therefore, a linear representation of a cat$^1$-group $\mathfrak{C}$ is a 2-functor $\mathfrak{C} \rightarrow \mathbf{Ch}^1_K$. The functorial image of $\mathfrak{C}$ lies within a sub 2-groupoid with a single object of $\mathbf{Ch}^1_K$. This sub 2-groupoid $\mathbf{Aut(\delta)}$, called automorphism cat$^1$-group, was obtained by considering only the invertible chain maps over a single chain complex of length-1; $\delta:C_1 \rightarrow C_0$ of vector spaces and has a cat$^1$-group structure. Thus, the 0-cell in $\mathbf{Aut(\delta)}$ is just the chain complex $\delta$ of length-1, and 1-cells in $\mathbf{Aut(\delta)}$ are chain maps $\delta \rightarrow \delta$ and 2-cells are 1-homotopies between chain maps. Since a linear representation of a group $G$ realises that group as a subgroup of $GL(V)$, so a linear representation of a cat$^1$-group $\mathfrak{C}$ realises that it as a cat$^1$-subgroup $\mathbf{Aut(\delta)}$ of $\mathbf{Ch}^1_K$.

A common approach to representations of groups is via modules over a group or an algebra \cite{Burrow}, \cite{Curtis}, \cite{Dornhoff}. Linear representations of a group $G$ are in one-to-one correspondence with modules over its group algebra, $K(G)$, see \cite{Barker}, where $K$ is the group algebra functor from the category of groups to that of algebras. Since a cat$^1$-group is a generalisation of a group; there must be a notion of cat$^1$-group algebra. For any cat$^1$-group $\mathfrak{C}$; by applying the group algebra functor $K$ to $\mathfrak{C}$, Barker, first obtained the structure $K(\mathfrak{C})$ as a pre-cat$^1$-algebra. In order to construct a cat$^1$-algebra from $K(\mathfrak{C})$, it is necessary to impose some relations so that the kernel condition is satisfied for the structural homomorphisms in $K(\mathfrak{C})$. Barker found suitable expressions in $K(\mathfrak{C})$ and by factoring $K(\mathfrak{C})$ by the ideal $J$ generated by these expressions, he showed that every cat$^1$-group, $\mathfrak{C}$, has an associated cat$^1$-group algebra $\overline{K(\mathfrak{C})}$. This structure was obtained by first applying the group algebra functor to $\mathfrak{C}$ and factoring the algebra of 1-cells by a given ideal $J$ in order to introduce relations necessary to satisfy the kernel condition in $\overline{K(\mathfrak{C})}$.

To construct a regular representation  given by Cayley's theorem for cat$^1$-groups; Barker has considered the chain complex; $\overline{\delta}$ of length-1 which is obtained from $\overline{K(\mathfrak{C})}$ and defined the structure $\mathbf{Aut(\overline{\delta})}$ as a subcat$^1$-group algebra of $\overline{K(\mathfrak{C})}$. Therefore, a regular representation of a cat$^1$-group $\mathfrak{C}$ is a 2-functor
$$\mathbf{\rho}:\mathfrak{C}\longrightarrow\mathbf{Aut(\overline{\delta})}.$$
Consequently, Barker's result can be summarized pictorially as;
$$
\begin{tikzcd}
\mathfrak{C} \arrow[r,"K"] \arrow[ddrr,dashed,"\rho"," \mathrm{(regular \ representation)}"{sloped,below=0.0ex,xshift=-0.0em}]  &K(\mathfrak{C}) \arrow[r,twoheadrightarrow,"q"] &\overline{K(\mathfrak{C})} \arrow[d]
\\
 & &\ \ \ \  \overline{\delta} \in \mathbf{Ch}^1_K \arrow[d]
\\
 & &\mathbf{Aut(\overline{ \delta})}
\end{tikzcd}
$$
where $\mathfrak{C}$ is any cat$^1$-group; $K$ is the group algebra functor, $q$ is a quotient functor from pre-cat$^1$-algebras to cat$^1$-algebras and $\rho$ is the right regular representation of $\mathfrak{C}$. Note that the functor $\rho$ defined in \cite{Barker} is contravariant in the horizontal direction but covariant in the vertical direction. Thus, this functor can be regarded as a functor $\mathfrak{C}^{op} \rightarrow \mathbf{Aut}(\overline\delta)\leqslant \mathbf{Ch}^1_K$, considering the convention given in \cite{Kelly} that any 2-category $\mathfrak{C}$ has a dual $\mathfrak{C}^{op}$ which reverses only the 1-cells and have affects both 1-cell and horizontal 2-cell composition. Another dual $\mathfrak{C}^{co}$ which reverses only the 2-cells and hence affects vertical 2-cell composition. Therefore, with these constructions, it is obtained that a cat$^1$-group version of Cayley's theorem exists in terms of linear regular representations.

In \cite{Elgueta}, Elgueta has obtained an alternative representation of 2-groups in the 2-category of finite dimensional 2-vector spaces (over a field $K$) as defined Kapranov and Voevodsky \cite{Kapranov}. The notion of 2-vector space introduced by these authors is different from that of Baez and Crans \cite{Baez1}. It was proven in \cite{Baez1} that the category of 2-vector spaces is equivalent to the category of 2-term chain complexes of vector spaces, which is of course none other than $\mathbf{Ch}^1_K$. Thus, Barker's result is related to the definition of Baez and Crans, and Elgueta's result is related to that of Kapranov and Voevodsky.

Gray in his lecture notes developed the notion of the tensor product for 2-categories (cf. \cite{Gray}). Restricting this construction to 2-groupoids gives a basic example of a monoidal category of 2-groupoids with monoidal structure given by Gray's tensor product. Joyal and Tierney, \cite{JT}, proved that Gray groupoids model homotopy 3-types. As an alternative algebraic model for homotopy 3-types, Conduch{\'e}, \cite{Con}, defined 2-crossed modules  and showed how to obtain a 2-crossed  module from a simplicial group. See also for this construction, \cite{mutpor2, mutpor3}, in terms of Carrasco-Cegarra pairing operators (cf. \cite{cc}). 2-Crossed modules are also equivalent to the notions of crossed squares introduced by Loday and Guin-Walery in \cite{WL}, and braided regular crossed modules introduced by Brown and Gilbert in \cite{BG}. For these connections, see also \cite{AU1} and  \cite{Con1}. A connection between 2-crossed modules of group(oid)s and Gray 3-groupoids was proven by Kamps and Porter in \cite{KP}. Using a different method; Martins and Picken, \cite{Martins}, gave the relationship between 2-crossed modules of groups and Gray-3-groupoids with a single 0-cell. This construction has also been found in Al-asady's Ph.D thesis (cf. \cite{Jinan}). See also \cite{Wang}. Kamps and Porter, \cite{KP}, have proven that the category of chain complexes of length-2 of vector spaces over a field $K$,\ $\mathbf{Ch}^2_K$, has a Gray category structure. In this constructions, 0-cells of $\mathbf{Ch}^2_K$ are chain complexes of length-2, 1-cells are chain maps, 2-cells are 1-homotopies and 3-cells are 2-homotopies between 1-homotopies. For any chain complex of length-2 within $\mathbf{Ch}^2_K$;
$
\begin{tikzcd}
             \delta:=C_{2}\ar[r,"{\delta_2}"]&C_{1}\ar[r,"{\delta_{1}}"]& C_{0}
\end{tikzcd}
$
of linear transformations, the structure $\mathbf{Aut(\delta)}$ as a cat$^2$-group, was introduced in \cite{Jinan}. This structure in fact is a Gray 3-groupoid with a single object set $\{\delta\}$. Using the equivalance between crossed squares and cat$^2$-groups (cf. \cite{WL}); Al-asady, in \cite{Jinan}, has constructed a linear representation of a cat$^2$-group $\mathfrak{C}^2$, as a lax 3-functor $\mathfrak{C}^2 \rightarrow \mathbf{Aut(\delta)}\leqslant \mathbf{Ch}^2_K$. In this construction; 0-cell of $\mathbf{Aut(\delta)}$ is $\delta$, 1-cells are chain maps $F:\delta \rightarrow \delta$, 2-cells are 1-homotopies $(H,F):F \Rightarrow G$ and 3-cells are 2-homotopies $(\alpha,H,F):(H,F)\Rrightarrow(K,F)$. Therefore, Al-asady, by generalising the 2-category of length-1 chain complexes $\mathbf{Ch}^1_K$ to a Gray category of length-2 chain complexes $\mathbf{Ch}^2_K$, has given the linear representation of a cat$^2$-group $\mathfrak{C}^2$.

In this paper, we give the regular representation of a cat$^2$-group or a Gray 3-groupoid with a single object $*$, $\mathfrak{C}^2$, obtained from a 2-crossed module $\mathfrak{X}$. By applying the group algebra functor to $\mathfrak{C}^2$, we obtain the structure $K(\mathfrak{C}^2)$. We see that $K(\mathfrak{C}^2)$ is certainly a pre-cat$^2$-algebra, but not a cat$^2$-algebra. To make it a cat$^2$-algebra; we need to give some relations for $K(\mathfrak{C}^2)$ to satisfy the kernel condition. By defining ideals $J_2$ and $J_1$ in $K(\mathfrak{C}^2)$ and factoring it by these ideals generated by these relations, we have obtained a Gray 3-group algebra groupoid with a single object or a cat$^2$-group algebra $\overline{K(\mathfrak{C}^2)}$. We construct a regular representation of $\mathfrak{C}^2$ as a 3-functor from $\mathfrak{C}^2$ to $\mathbf{Aut(\overline{ \delta})}$, where
$
\begin{tikzcd}
             \overline{\delta}:=\mathrm{K}_3\ar[r,"\overline{\tau}_3"]&\mathrm{K}_2\ar[r,"\overline{\tau}_2"]& \mathrm{K}_1
\end{tikzcd}
$
is the chain complex of length-2 obtained from $\overline{K(\mathfrak{C}^2)}$. Thus, the regular representation of $\mathfrak{C}^2$ is a 3-functor $\lambda:\mathfrak{C}^2 \rightarrow \mathbf{Aut(\overline{\delta})}\leqslant \mathbf{Ch}^2_K$. In the construction of this functor, for the 0-cell $*$ in $\mathfrak{C}^2$, we obtain that $\lambda_*$ is $\overline{\delta}$ and for a 1-cell in $\mathfrak{C}^2$; we obtain that $\lambda_n$ is a chain map $\overline{\delta} \rightarrow \overline{\delta}$, similarly for a 2-cell $(m,n)$ in $\mathfrak{C}^2$; we have that $\lambda_{m,n}$ is a 1-homotopy from $\lambda_n$ to $\lambda_{\partial_1mn}$ and for a 3-cell $(l,m,n)$ in $\mathfrak{C}^2$; we have that $\lambda_{l,m,n}$ is a 2-homotopy from $\lambda_{m,n}$ to $\lambda_{\partial_2lm,n}$. Thus, our results can be summarised pictorially as;
$$
\begin{tikzcd}
{\mathfrak{C}^2} \arrow[r,"K"] \arrow[ddrr,dashed,"\lambda","\mathrm{(regular \ representation)}"{sloped,below=0.0ex,xshift=-0.0em}]  & {K(\mathfrak{C}^2)} \arrow[r,twoheadrightarrow,"q"] &\overline{K(\mathfrak{C}^2)} \arrow[d]
\\
 & &\ \ \ \ \overline{\delta} \in \mathbf{Ch}^2_K \arrow[d]
\\
 & &\mathbf{Aut(\overline{ \delta})}
\end{tikzcd}
$$
where $\mathfrak{C}^2$ is a cat$^2$-group obtained from the 2-crossed module $\mathfrak{X}$ and $\mathbf{Aut(\overline{ \delta})}$ is the cat$^2$-group algebra obtained from the chain complex $\overline{\delta}$ within $\overline{K(\mathfrak{C}^2)}$.
\tableofcontents

\section{Preliminaries}

Recall that a \textit{small category} $\mathfrak{C}$ consists of an object set $C_0$, a set of morphisms $C_1$, source and target maps from $C_1$ to $C_0$, a map $i:C_{0}\rightarrow C_1$ which gives the identity morphisms at an object and a partially defined function $C{_1}\times C_{1}\rightarrow C_1$ which gives the composition of two morphisms. We will show a small category $(C_1 ,C_0 )$ and diagramatically as
$$
\xymatrix{C_1 \ar@<1ex>[r]^{s,t}\ar@<0ex>[r]&C_0. \ar@<1ex>[l]^{i}}
$$
For the set of morphisms $C_1$, and $x,y\in C_0$ the set of morphisms from $x$ to $y$ is written $C_{1}(x,y)$ and termed a hom-set. Then for $a\in C_{1}(x,y)$, we have $s(a)=x$ and $t(a)=y$. We will usually write $i_x$ for $i(x)$ and $b\circ a$ for the composite of the morphisms $a:x\rightarrow y$ and $b:y\rightarrow z$. The elements of $C_0$ are also called 0-cells and the elements of $C_1$ are called 1-cells between 0-cells.

A \textit{groupoid} is a small category in which every morphism (or every 1-cell) is an isomorphism (or invertible), that is, for any 1-cell  $(a:x\rightarrow y)\in C_{1}(x,y)$, there is a 1-cell  $(a^{-1}:y\rightarrow x)\in C_{1}(y,x)$, such that $a^{-1}\circ a=i_x$ and $a\circ a^{-1}=i_y$. A groupoid with a single 0-cell can be regarded as a group. For a survey of application of groupoids and introduction to their literature, see \cite{Brown, BG, bs}.

Cat$^1$-groups or categorical groups are group objects in the category of small categories (cf. \cite{Loday}). They are sometimes referred to simply as cat-groups. A cat$^1$-group $\mathfrak{C}:=(G,N,i,s,t)$ consists of groups $G$ and $N$ and embedding $i:N\rightarrow G$, and epimorphisms $s,t:G\rightarrow N$ satisfying the conditions: (i) $si=ti=id_N$ and (ii) $[\ke s,\ke t]={1_G}.$ Condition (ii) is called the \textit{kernel condition}. A structure with the same data as a cat$^1$-group and satisfying the first condition but not the kernel condition is called a \textit{pre-cat$^1$-group}.

Crossed modules were introduced by Whitehead in \cite{W}. A crossed module $\mathfrak{X}:=(M,N,\partial)$ consists of groups $M,N$ together with a homomorphism $\partial:M\rightarrow N$ and a left action $N\times M \rightarrow M$ of $N$ on $M$ given by $(n,m)\mapsto {^n{m}}$, satisfying the conditions: (i) $\partial(^n{m})=n\partial(m)n^{-1}$ and (ii) $^{\partial (m)}{m'}=mm'm^{-1}$ for all $n\in N,\ m,m'\in M $. Condition (ii) is called the \textit{Peiffer identity}. A structure with the same data as a crossed module and satisfying the first condition but not the Peiffer identity is called a \textit{pre-crossed module}.

If $M$ and $N$ are groups with a left action of $N$ on $M$, the semi-direct product of $M$ by $N$ is the group $M\rtimes N=\{(m,n):m \in M, n\in N\}$ with the multiplication $(m,n)(m',n')=(m^{n}{m'},nn')$. The inverse of $(m,n)$ is $(^{n^-1}m^{-1},n^{-1})$.

For a cat$^1$-group $\mathfrak{C}:=(G,N,i,s,t)$, it is well-known that $G\cong\ke s\rtimes N$, (cf. \cite{Brown-loday, bs}). From a crossed module $\mathfrak{X}:=(M,N,\partial)$, a cat group $\mathfrak{C(X)}$:
$$
\xymatrix{M\rtimes N \ar@<1ex>[r]^-{s,t}\ar@<0ex>[r]&N \ar@<1ex>[l]^-{i}}
$$
can be constructed. Here $s,t:M\rtimes N\rightarrow N$ ; $i:N\rightarrow M \rtimes N$ are defined as $s(m,n)=n$, $t(m,n)=\partial(m)n$ and $i(n)=(1,n)$. Then $si=ti=id$ and $[\ke s,\ke t]=1_{M\rtimes N}$. Thus, crossed modules and cat$^1$-groups are equivalent algebraic structures.

In a cat$^1$-group $\mathfrak{C}:=(G,N,i,s,t)$, since the big group $G$ can be decomposed as $\ke s \rtimes N$, we can describe a typical element is $(m,n)$, $m\in \ke s$ and $n\in N$. So we can view this as a 2-cell $(m,n):n\Rightarrow\partial_1mn$. In this case, the 0-cell of $\mathfrak{C}$ is $*$, and 1-cells are $n:*\rightarrow *$ with $n\in N$. 2-cells are $(m,n):n\Rightarrow\partial_1mn$. These are given pictorially as
$$
 \begin{tikzcd}[row sep=tiny,column sep=small]
                & \ar[dd, Rightarrow, "{\scriptscriptstyle(m,n)}"] \\
              {*}\ar[rr,bend left=50,"\scriptscriptstyle n"] \ar[rr,bend right=50,"\scriptscriptstyle \partial_1mn"']  &  & \ {*} \\
                 & \
    \end{tikzcd}
$$
Since $G$ is itself a group; it can be regarded as a category with a single object $*$ and its group operation is considered categorically as composition. Therefore,
$(m,n)\#_0(m',n')=(m,n)(m',n')=(m^{n}m',nn')$
and we can represent it pictorially as
$$
\begin{array}{c}
    \begin{tikzcd}[row sep=small,column sep=0.8cm]
                    &  \ar[dd,Rightarrow,"\scriptscriptstyle{(m,n)}"{description}] &    & \ar[dd,Rightarrow,"\scriptscriptstyle{(m',n')}"{description}]& \\
                 {*}\ar[rr,bend left=50,"\scriptscriptstyle n"] \ar[rr,bend right=50,"\scriptscriptstyle \partial_1mn"']& &{*}\ar[rr,bend left=50,"\scriptscriptstyle n'"] \ar[rr,bend right=50,"\scriptscriptstyle \partial_1m'n'"'] &  & {*} \\ \ & \ & \  & \
    \end{tikzcd}
\end{array}
{:=}
\begin{array}{c}
    \begin{tikzcd}[row sep=small,column sep=0.8cm]
                 & \ar[dd, Rightarrow, "\scriptscriptstyle{(m,n)\#_0(m',n')}"{description}] & \\
                 {*}\ar[rr,bend left=50,"\scriptscriptstyle nn'"] \ar[rr,bend right=50,"\scriptscriptstyle \partial_1mn\partial_1m'n'"'] &  &{*}\\
                 & \
    \end{tikzcd}
\end{array}
$$
This operation is called horizontal composition of 2-cells.
On the other hand, there is another composition in $\mathfrak{C}$. For any 2-cell $(m',\partial(m)n):\partial(m)n\Rightarrow\partial m'\partial mn$ where $m'\in \ke s$, we can compose $(m,n)$ and $(m',\partial mn)$ by defining $(m',\partial(m)n)\#_2(m,n)=(m'm,n)$. The kernel condition  establishes an interchange law between these compositions $\#_2$ and $\#_0$. Hence, a cat$^1$-group can be considered as a 2-category with a single object $*$. These compositions are invertible, so we can say that a cat$^1$-group has a structure of a 2-group, or equivalently has a structure of a 2-groupoid with a single 0-cell $*$.
\subsection{The Category $\mathbf{Ch}^1_K$ as a 2-Groupoid}

Barker, in \cite{Barker}, has constructed a 2-groupoid structure over chain complexes of length-1. In this section, we give a brief description of the category of chain complexes of vector spaces over a fixed field $K$ and its 1-truncation case $\mathbf{Ch}^1_K$. Now, suppose that $K$ is a field and $C_i$ is a $K$-vector space for each $i\in \mathbb{Z}$. Let
$$
\begin{tikzcd}
              \mathfrak{C}:=\cdots\ar[r]&C_{n}\ar[r,"{d_n}"]&C_{n-1}\ar[r,"{d_{n-1}}"]&\cdots C_{1}\ar[r,"{d_1}"]&C_{0}\ar[r,"{d_{0}}"]&C_{-1}\ar[r,"{d_{-1}}"]&\cdots
\end{tikzcd}
$$
be a chain complex of $K$-vector spaces and linear transformations between them. A chain map from $f:\mathfrak{C}\rightarrow \mathfrak{D}$ consists of components $f_i:C_i\rightarrow D_i$ such that $f_{i-1}d_i=d_if_i$ for all $i\in \mathbb{Z}$, i.e. the following diagram commutes:
$$
\begin{array}{c}
 \xymatrix{\mathfrak{C}\ar[d]_{f}\\
 \mathfrak{D}}
 \end{array}
 {:=}
 \begin{array}{c}
\xymatrix{ \cdots \ar[r]&C_{n+1} \ar[r]^-{d_{n+1}}\ar[d]^-{f_{n+1}}& C_{n} \ar[r]^-{d_{n}}\ar[d]^-{f_n}&C_{n-1} \ar[r]^-{d_{n-1}}\ar[d]^-{f_{n-1}}  & \cdots\\
 \cdots\ar[r]&D_{n+1} \ar[r]_-{d_{n+1}} & D_{n} \ar[r]_-{d_{n}} \ar[r] &D_{n-1} \ar[r]_-{d_{n-1}} & \cdots}
\end{array}
$$
where each $f_i$ is a linear map. The composition $g\#_0 f:\mathfrak{C}\rightarrow \mathfrak{E}$ of the chain maps $f:\mathfrak{C}\rightarrow \mathfrak{D}$ and $g:\mathfrak{D}\rightarrow \mathfrak{E}$ is defined by $(g\#_0 f)_i=g_if_i$ for all $i$, where the composition on the right hand side is the usual one for linear maps. Note that a chain isomorphism $f:\mathfrak{C}\rightarrow \mathfrak{D}$ is an invertible chain map in which each component is a linear isomorphism.

Let $f,g$ be chain maps from $\mathfrak{C}$ to $\mathfrak{D}$. A chain homotopy $H:f\simeq g$ consists of maps $h'_n:C_n\rightarrow D_{n+1}$ satisfying $g_n-f_n=d_{n+1}h'_n+h'_{n-1}d_n$ for each $n\in \mathbb{Z}$, as pictured in the following diagram:
\begin{equation*}
\xymatrix{\cdots\ar[rr]&&C_{n+1} \ar[rr]^-{\scriptscriptstyle d_{n+1}}\ar@<0.5ex>[dd]^-{\scriptscriptstyle g_{n+1}}\ar@<-0.5ex>[dd]_-{\scriptscriptstyle f_{n+1}} && C_{n} \ar[ddll]_{\scriptscriptstyle h'_n} \ar[rr]^-{\scriptscriptstyle d_{n}}\ar@<0.5ex>[dd]^-{\scriptscriptstyle g_n}\ar@<-0.5ex>[dd]_-{\scriptscriptstyle f_n}&&C_{n-1} \ar[ddll]_{\scriptscriptstyle h'_{n-1}}\ar[rr]^-{\scriptscriptstyle d_{n-1}}\ar@<0.5ex>[dd]^-{\scriptscriptstyle g_{n-1}}\ar@<-0.5ex>[dd]_-{\scriptscriptstyle f_{n-1}}   &&\cdots
\\
\\
\cdots \ar[rr]&&D_{n+1}\ar[rr]_{\scriptscriptstyle d_{n+1}} && D_{n} \ar[rr]_{\scriptscriptstyle d_{n}} \ar[rr] &&D_{n-1} \ar[rr]_{\scriptscriptstyle d_{n-1}} &&\cdots}
\end{equation*}
The category of all chain complexes is called $\mathbf{Ch}$ and is discussed in detail by Kamps and Porter in \cite{KP}. They have proven that $\mathbf{Ch}$ is a 2-groupoid enriched Gray category. In our work, it is sufficient to consider non-negative chain complexes in which the subscripts are non-negative integers.

The chain complexes of length-0 gives us the category of vector spaces over $K$. We can denote it by $\mathbf{Ch}^0_K$. Now, we concentrate the category $\mathbf{Ch}^1_K$. If $\delta:C_1\rightarrow C_0$ is a linear transformation, then $\mathfrak{C}:=(\xymatrix{C_1\ar[r]^{\delta}&C_0})$ is a length-1 chain complex of vector spaces. Then $\mathfrak{C}$ can be considered as a chain complex
$$
\begin{tikzcd}
              \mathfrak{C}:=\cdots\ar[r]&0\ar[r]&C_{1}\ar[r,"\delta"]&C_{0}\ar[r]&0\ar[r]&\cdots
\end{tikzcd}
$$
Therefore, we obtain the category $\mathbf{Ch}^1_K$ whose objects are chain complexes of length-1 of $K$-vector spaces and morphisms are the chain maps between them. Clearly, $\mathbf{Ch}^1_K$ is a subcategory of $\mathbf{Ch}$.
Let $f,g:\mathfrak{C}\rightarrow \mathfrak{D}$ be chain maps in $\mathbf{Ch}^1_K$. Then, the homotopy $H:(H',f):f\simeq g$ is given by the following diagram
$$
\xymatrix{C_1\ar@<0.5ex>[rr]^{f_1}\ar@<-0.5ex>[rr]_{g_1}\ar[dd]_{\delta^C}&&D_1 \ar[dd]^{\delta^D}\\
\\
C_0\ar@<0.5ex>[rr]^{f_0}\ar@<-0.5ex>[rr]_{g_0}\ar[uurr]^{H'}&&D_0}
$$
where $H'$ is called the chain homotopy component, satisfying the conditions  $\delta^D H'=g_0-f_0$ and $H'\delta^C=g_1-f_1$. We can summarize the structure of $\mathbf{Ch}^1_K$ as follows:\\
$\bullet$ the 0-cells are the chain complexes of length-1; $\mathfrak{C}:=(\xymatrix{C_1\ar[r]^{\delta^C}& C_0})$,\\
$\bullet$ the 1-cells between 0-cells are the chain maps; $f:=(f_1,f_0):\delta^C\rightarrow \delta^D$,\\
$\bullet$ the 2-cells between 1-cells are the 1-homotopies; $H:=(H',f):f\simeq g$.

A 2-cell $H:=(H',f):f\simeq g$ may be written briefly as $H:f\Rightarrow g$. Let us now explain the vertical and horizontal compositions of 1-cells and 2-cells. The composition of 1-cells is the usual composition of the chain maps as explained before. For 0-cells $\mathfrak{C}$ and $\mathfrak{D}$, and 1-cells $f,g,k:\mathfrak{C}\rightarrow \mathfrak{D}$, the vertical composition of 2-cells $H:f\Rightarrow g$ and $\hat{H}:g\Rightarrow k$ is given by $\hat{H}\#_2 H:f\Rightarrow k$, where the chain homotopy component is $(\hat{H}\#_2 H)'=\hat{H}'+H'$.

To define the horizontal composition of 2-cells, as a first step, we have to give the whiskerings of 1-cells on 2-cells on the right side and on the left side. For the notions of whiskerings and whiskered categories see also \cite{Brown}. Let $H:f\Rightarrow f':\mathfrak{C}\rightarrow \mathfrak{D}$ be a 2-cell and $g:=(g_1,g_0):\mathfrak{D}\rightarrow \mathcal{E}$ a 1-cell. The left whiskering of $g$ on $H$ is given by $g\natural_1 H:g\#_0 f\Rightarrow g\#_0 f'$ where the chain homotopy component is $(g\natural_1 H)'=g_1H'$. This terminology can be explained by the picture
$$
\begin{array}{c}
    \begin{tikzcd}[row sep=small,column sep=normal]
                 &\ar[dd,Rightarrow,"{H}"]&\\
                 {\mathfrak{C}}\ar[rr,bend left=50,"f"]\ar[rr,bend right=50,"{f'}"']& &{\mathfrak{D}}\ar[r,"g"]&\mathcal{E}\\
                 &\
    \end{tikzcd}
\end{array}
{:=}
\begin{array}{c}
    \begin{tikzcd}[row sep=small,column sep=normal]
                  &\ar[dd,Rightarrow,"{g\natural_1 H}"]&\\
                  \mathfrak{C} \ar[rr,bend left=50,"g\#_0 f"]\ar[rr,bend right=50,"g\#_0 f'"']& &\mathcal{E}.\\
                  &\
    \end{tikzcd}
\end{array}
$$
The left whiskering of $g$ on $H$ appears on the left in the notation $g\natural_1 H$, but on the right in the picture. Similarly, the right whiskering of a 1-cell $f:=(f_1,f_0):\mathfrak{C}\rightarrow \mathfrak{D}$ on a 2-cell $K:g\Rightarrow g':\mathfrak{D}\rightarrow \mathcal{E}$ is given by $K\natural_1 f:g\#_0 f\Rightarrow g'\#_0 f$ where the chain homotopy component is $(K\natural_1 f)'=K'f_0$. In this case, the 1-cell $f$ will appear on the right in the notation $K\natural_1 f$, but on the left in the following picture:
$$
\begin{array}{c}
    \begin{tikzcd}[row sep=small,column sep=normal]
                    & \ & \  \ar[dd,Rightarrow,"{K}"] & \\
                \mathfrak{C} \ar[r,"f"]&\mathfrak{D}\ar[rr,bend left=50,"g"] \ar[rr,bend right=50,"g'"'] &  & \mathcal{E} \\
                \ & \ & \
    \end{tikzcd}
\end{array}
{:=}
\begin{array}{c}
    \begin{tikzcd}[row sep=small,column sep=normal]
                 & \ar[dd, Rightarrow, "{K\natural_1 f}"] & \\
                 \mathfrak{C}\ar[rr,bend left=50,"g\#_0 f"] \ar[rr,bend right=50,"g'\#_0 f"'] &  &\mathcal{E}.\\
                 & \
    \end{tikzcd}
\end{array}
$$
Now we can give the definition of horizontal composition for 2-cells. Let $H:f:(f_1,f_0)\Rightarrow f':(f'_1,f'_0):\mathfrak{C}\rightarrow \mathfrak{D}$ and $K:g:(g_1,g_0)\Rightarrow g':(g'_1,g'_0):\mathfrak{D}\rightarrow \mathcal{E}$ be 2-cells. The horizontal composition of $K$ and $H$ is $$K\#_0H:=(K\natural_1 f')\#_1(g\natural_1 H):g\#_0 f\Rightarrow g'\#_0 f'$$ where the chain homotopy component is $(K\#_0 H)'=g_1H'+K'f'_0$. The following picture shows the definition of horizontal composition of 2-cells:
$$
\begin{array}{c}
    \begin{tikzcd}[row sep=small,column sep=normal]
                 & \ \ar[dd, Rightarrow, "H"]  \ & \ &  \ar[dd, Rightarrow, "K"] &  \\
                \mathfrak{C}\ar[rr,bend left=50,"f"] \ar[rr,bend right=50,"f'"'] & \ & \mathfrak{D}\ar[rr,bend left=50,"g"]
                 \ar[rr,bend right=50,"g'"'] & \ & \mathcal{E}\\
                \ & \ & \ & \ & \
    \end{tikzcd}
\end{array}
{:=}
\begin{array}{c}
    \begin{tikzcd}[row sep=small,column sep=normal]
                 & \ar[dd, Rightarrow, "{K\#_0 H}"] & \\
                 {\mathfrak{C}}\ar[rr,bend left=50,"g\#_0 f"] \ar[rr,bend right=50,"g'\#_0 f'"']  &    & {\mathcal{E}}.  \\
                 & \
    \end{tikzcd}
\end{array}
$$
We also note from \cite{Barker} that $$(g'\natural_1 H)\#_1(K\natural_1 f)=(K\natural_1f')\#_1 (g\natural_1 H). $$
The interchange law for the vertical and horizontal compositions was given in \cite{Barker}. Therefore, $\mathbf{Ch}^1_K$ is a 2-category and is a groupoid-enriched category. By restriction the chain maps to those which are invertible, the category of invertible chain maps $\mathbf{invCh}^1_K$ is a 2-groupoid.

\subsection{The Structure Aut($\delta$) in $\mathbf{Ch}^1_K$ as a Cat$^1$-Group}

To define a linear representation of a crossed module (or equivalently a cat$^1$-group) $\mathfrak{C}$ as a 2-functor $\phi: \mathfrak{C}\rightarrow \mathbf{Ch}^1_K$, Barker in \cite{Barker} has constructed an automorphism cat$^1$-group $\mathbf{Aut(\delta)}$ as a 2-subcategory of $\mathbf{Ch}^1_K$. In this case, the functorial image of $\mathfrak{C}$ under the 2-functor $\phi$ is $\mathbf{Aut(\delta)}$. This subcategory is a collection of all chain isomorphisms $\delta \rightarrow \delta$ and homotopies between them.
\begin{defn}\rm{ (\cite{Barker})
Let $\delta:C_1\rightarrow C_0$ be a linear transformation of vector spaces over the field $K$. The automorphism cat$^1$-group of $\delta$; $\mathbf{Aut(\delta)}$ consists of:\\
$(i)$ the group $\mathbf{Aut}(\delta)_1$ of all chain isomorphisms $\delta\rightarrow \delta$,\\
$(ii)$ the group $\mathbf{Aut}(\delta)_2$  of all homotopies on $\mathbf{Aut}(\delta)_1$ ,\\
$(iii)$ morphisms $s,t:\mathbf{Aut}(\delta)_2 \rightarrow \mathbf{Aut}(\delta)_1$, selecting the source and target of each homotopy,\\
$(iv)$ the morphism $i:\mathbf{Aut}(\delta)_1 \rightarrow \mathbf{Aut}(\delta)_2$, which provides the identity homotopy on each chain automorphism.}
\end{defn}
Then, the only 0-cell in  $\mathbf{Aut(\delta)}$ is $\delta:C_1\rightarrow C_0$ and so $\mathbf{Aut}(\delta)_0$ is a singleton. 1-cells are chain isomorphisms from $\delta$ to itself, $f:(f_1,f_0):\delta\rightarrow \delta$. 2-cells are homotopies between 1-cells:
$$
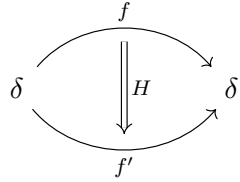
\begin{tikzcd}[row sep=small,column sep=normal]
                 &\ar[dd,Rightarrow,"{ H}"]\\
                  \delta \ar[rr,bend left=50,"f"] \ar[rr,bend right=50,"f'"']& &\delta\\
                 &\
\end{tikzcd}
$$
where $H=(H',f):f\Rightarrow f'$ and $H'$ is the chain homotopy component and it can be more clearly represented by the diagram
$$
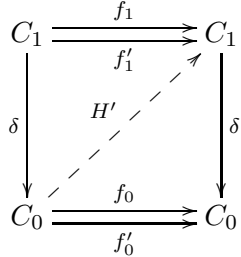
\xymatrix{C_1\ar@<0.5ex>[rr]^{f_1}\ar@<-0.5ex>[rr]_{f_1'}\ar[dd]_{\delta}&&C_1 \ar[dd]^{\delta}\\
\\
C_0\ar@<0.5ex>[rr]^{f_0}\ar@<-0.5ex>[rr]_{f_0'}\ar@{-->}[uurr]^{H'}&&C_0}
$$
with the homotopy conditions: $\delta H'=f'_0-f_0$ and $H'\delta =f'_1-f_1$. Since a 2-groupoid with a single object can be considered as a cat$^1$-group,
$$
\xymatrix{\mathbf{Aut(\delta)}:=\mathbf{Aut}(\delta)_1  \ar@<1ex>[r]^-{s,t}\ar@<0ex>[r]&\mathbf{Aut}(\delta)_0  \ar@<1ex>[l]^-{i}}
$$
is a cat$^1$-group.

Therefore, the 2-functor $\phi:\mathfrak{C}\rightarrow \mathbf{Aut(\delta)}\leqslant \mathbf{Ch}^1_K$ can be explained as follows:

$\bullet$ For the 0-cell $*$ in $\mathfrak{C}$, $\phi(*)=\delta$ is the chain complex of length-1 in $\mathbf{Ch}^1_K$ which is the 0-cell of $\mathbf{Aut(\delta)}$.

$\bullet$ For any 1-cell $x:* \longrightarrow *$, $x\in C_0$, $\phi(x)=f$ is the chain isomorphism $f:=(f_1,f_0):\delta\rightarrow \delta$ which is a 1-cell in $\mathbf{Aut}(\delta)_1$ .

$\bullet$ For any 2-cell $a:x\rightarrow y$ in $C_1$ with $s(a)=x$, $t(a)=y$, we have $\phi(a)=H:=(H',f):f\Rightarrow f'$  in $\mathbf{Aut}(\delta)_2$ is the homotopy from $\phi(x)=f$ to $\phi(y)=f'$.
\subsection{Regular Representation of Crossed Modules and Cat$^1$-Groups}

It is well-known that Cayley's theorem provides  regular representations of groups. For any group $G$, the right regular representation of $G$ is defined as $\lambda:G^{op}\rightarrow S_{|G|}$ with $\lambda_g(h)=hg$ (for $g,h\in G$) where $\lambda_g:G\rightarrow G$ is a permutation of the underlying set of $G$ and $G^{op}$ is the opposite category to $G$, since the right multiplication gives us $\lambda$ to be contravariant, i.e. $\lambda_{g_1g_2}(h)=h(g_1g_2)=(hg_1)g_2=\lambda_{g_2}(\lambda_{g_1})(h)$ for $g_1,g_2,h\in G$. Similarly, the left regular representation is defined by $\rho:G\rightarrow S_G$ with $\rho_g(h)=gh$. In this case, $\rho$ is covariant.

 Serre, \cite{Serre}, by reformulating Cayley's theorem for linear representations, defined the left regular representation as a linear representation of $G$. In this case, to the regular representations of $G$ there correspond regular linear representations $\lambda:G^{op}\rightarrow GL_{|G|}(K),\  \lambda:G\rightarrow GL_{|G|}(K)$ with $\lambda_g(\mathbf{e}_h):=\mathbf{e}_{hg}$ where $\mathbf{e}_g$ are the basis vectors in the group algebra $K(G)$. In this structure, $\lambda_g(\mathbf{e}_h):=\mathbf{e}_{hg}$ can be considered as $^{g}\mathbf{e}_{h}=\mathbf{e}_{hg}$ an action of $G$ on the group algebra $K(G)$ by multiplication on the right. Barker, \cite{Barker}, using the group algebra functor $K(.)$ from the category of groups to that of algebras has constructed regular representations of cat$^1$-groups. In this section, we will explain briefly this construction to use it in the next sections.
 \subsubsection{A Brief Description of the Group Algebra Functor $K(.)$}
 Let $G$ be a group and $K$ a fixed field. Suppose that $X_G$ is the underlying set of $G$. For any element $g\in G$, the notation $\mathbf{e}_g$ denotes the element in $X_G$. Then, $K(G)$ is a vector space over the field $K$ with basis $\{\mathbf{e}_g\}_{g\in G}$. Any element of $K(G)$ can be considered as in the form $\sum r_g\mathbf{e}_g$, with $r_g\in K$ and only finitely many $r_g\neq 0$ (cf. \cite{Passman}). The addition in $K(G)$ can be given by
 $$
\sum r_g\mathbf{e}_g +\sum s_g\mathbf{e}_g:=\sum (r_g+s_g)\mathbf{e}_g
 $$
 where $r_g+s_g$ is given by the addition in $K$, and the scalar multiplication can be given by
 $$
 s\sum r_g\mathbf{e}_g:=\sum (sr_g)\mathbf{e}_g,
 $$
  and for the elements $\sum r_g\mathbf{e}_g$ and $\sum s_h\mathbf{e}_h$ in $K(G)$ ($g,h\in G$), the multiplication in $K(G)$ is given by
$$
(\sum r_g\mathbf{e}_g)(\sum s_h\mathbf{e}_h):=\sum (r_gs_h)\mathbf{e}_{gh}
$$
where $(sr_g)$ and $(r_gs_h)$ are obtained by the multiplication in $K$ and $gh$ which appears in  $\mathbf{e}_{gh}$ is obtained by the operation in $G$.
Together with these operations, $K(G)$ is a $K$-algebra and it is called the group algebra of $G$. It is proven in \cite{Barker} that $K(.):\mathbf{Gr}\rightarrow \mathbf{Alg}_K$ is a functor from the category of groups to the category of $K$-algebras. For details see \cite{Barker} and \cite{Passman}. By ignoring the multiplication in the group algebra $K(G)$, we can consider the underlying vector space $K(G)$. We will also write $K(G)$ for this underlying vector space as well.

\subsubsection{The Regular Representation of Cat$^1$-Groups}
Using the group algebra functor for any cat$^1$-group, $\mathfrak{C}$, Barker has obtained a special cat$^1$-group algebra $\overline{K(\mathfrak{C})}$, and he defined from $\overline{K(\mathfrak{C})}$ a chain complex of length-1; $\overline{\delta}$ and then, by constructing a 2-functor $\lambda:\mathfrak{C}\rightarrow \mathbf{Aut}(\overline{\delta})$, he has given the regular representation of cat$^1$-groups. In this section, we recall this construction briefly. We will use this method to give the regular representation of 2-crossed modules and Gray 3-groupoids with a single 0-cell $*$.

Let $\partial:M\rightarrow N$ be a crossed module of groups. Recall that, using the action of $N$ on $M$, we can create the semi-direct product group $M\rtimes N$ together with the operation $(m,n)(m',n')=(m^{n}m',nn')$. It is known that
$$
\xymatrix{\mathfrak{C}:=M\rtimes N \ar@<1ex>[r]^-{s,t} \ar@<0ex>[r]&N \ar@<1ex>[l]^-{i}}
$$
is a cat$^1$-group with the structural maps given above. Applying the group algebra functor to this cat$^1$-group $\mathfrak{C}$, it is obtained that
$$
\xymatrix{K(\mathfrak{C}):=K(M\rtimes N) \ar@<1ex>[r]^-{\sigma,\tau} \ar@<0ex>[r]&K(N) \ar@<1ex>[l]^-{i}}
$$
is a pre-cat$^1$-algebra, where $K(M\rtimes N)$ and $K(N)$ are the group-algebras of the groups $M\rtimes N$ and $N$ respectively, and $K(N)$ has basis $\{\mathbf{e}_n:n\in N\}$, $K(M\rtimes N)$ has basis $\{\mathbf{e}_{m,n}:m\in M, n\in N\}$. The source, target and identity maps within this structure are defined by $\sigma(\mathbf{e}_{m,n})=\mathbf{e}_n$, $\tau(\mathbf{e}_{m,n})=\mathbf{e}_{\partial(m)n}$ and $i(\mathbf{e}_n)=\mathbf{e}_{1,n}$ respectively. It is clear that $\sigma i=\tau i =id.$ Any element of $\ke\sigma$ is in the form $\mathbf{v}_{m,n}=\mathbf{e}_{m,n}-\mathbf{e}_{1,n}$ since $\sigma(\mathbf{v}_{m,n})=0.$ Similarly, any element of $\ke\tau$ is in the form  $\mathbf{w}_{m,n}=\mathbf{e}_{m,n}-\mathbf{e}_{1,\partial(m)n}$. In this case, the set $\{\mathbf{v}_{m,n}:m\neq 1\}$ is a basis for $\ke\sigma$ and the set $\{\mathbf{w}_{m,n}:m\neq 1\}$ is a basis for $\ke\tau$. To satisfy the kernel condition $\ke\sigma.\ke\tau=0$, it would suffice that $\mathbf{v}_{m,n}.\mathbf{w}_{m',n'}=0$ and $\mathbf{w}_{m',n'}.\mathbf{v}_{m,n}=0$ for $\mathbf{w}_{m',n'}\in \ke\tau$ and $\mathbf{v}_{m,n}\in \ke\sigma$. However,
$$
\mathbf{v}_{m,n}.\mathbf{w}_{m',n'}=\mathbf{e}_{m^{n}m',nn'}-\mathbf{e}_{^{n}m',nn'}-\mathbf{e}_{m,n\partial(m')n'}+\mathbf{e}_{1,n\partial(m')n'}\neq 0.
$$
Since this expression is a linear combination of basis elements with non-zero coefficients, thus the kernel condition fails. Hence, in this form $K(\mathfrak{C})$ is not a cat$^1$-algebra but it is a pre-cat$^1$-algebra. In order to satisfy the kernel condition, Barker has found a suitable expression in $K(M\rtimes N)$ as
$$
\mathbf{e}_{m'm,n}-\mathbf{e}_{m,n}-\mathbf{e}_{m',\partial(m)n}+\mathbf{e}_{1,\partial(m)n}
$$
(for $m,m'\in M$ and $n\in N$) called \textit{cocycles}. $K(\mathfrak{C})$ can be factored by the ideal $J$  which is generated by these cocycles. This ideal is referred to as the cocycle ideal in $K(M\rtimes N)$. In this case, we can write
$$
J=\langle\mathbf{e}_{m'm,n}-\mathbf{e}_{m,n}-\mathbf{e}_{m',\partial(m)n}+\mathbf{e}_{1,\partial(m)n}:m,m'\in M\setminus\{1_M\}, n\in N\rangle.
$$

Let $\overline{K(M\rtimes N)}=K(M\rtimes N)/J$. Any  element in this factor algebra is the coset $\overline{\mathbf{e}}_{m,n}=\mathbf{e}_{m,n}+J$ for $\mathbf{e}_{m,n}\in K(M\rtimes N)$. Using the quotient map $q:K(M\rtimes N)\longrightarrow \overline{K(M\rtimes N)}$, the source, target and identity maps between $\overline{K(M\rtimes N)}$ and $K(N)$ are defined by
\begin{align*}
\overline{\sigma}(\overline{\mathbf{e}}_{m,n})&=\sigma(\mathbf{e}_{m,n})=\mathbf{e}_{n},\\
\overline{\tau}(\overline{\mathbf{e}}_{m,n})&=\tau(\mathbf{e}_{m,n})=\mathbf{e}_{\partial(m)n},\\
\overline{i}(\mathbf{e}_{n})&=\overline{\mathbf{e}}_{1,n}.
\end{align*}
Since
$$
\mathbf{v}_{m,n}.\mathbf{w}_{m',n'}=\mathbf{e}_{m^{n}m',nn'}-\mathbf{e}_{^{n}m',nn'}-\mathbf{e}_{m,n\partial(m')n'}+\mathbf{e}_{1,n\partial(m')n'}\in J,
$$
we obtain $\overline{\mathbf{v}}_{m,n}.\overline{\mathbf{w}}_{m',n'}=0+J=\overline{0}$. Thus, $\ke\overline{\sigma}.\ke\overline{\tau}=\overline{0}$. Hence, together with these structures
$$
\xymatrix{\overline{K(\mathfrak{C}}):=K(M\rtimes N)/J \ar@<1ex>[r]^-{\overline{\sigma},\overline{\tau}} \ar@<0ex>[r]&K(N) \ar@<1ex>[l]^-{\overline{i}}}
$$
is a cat$^1$-algebra. Thus, the construction of a regular representation can be summarised in the following definition (cf. \cite{Barker}).

\begin{defn}\rm{
The right regular representation of a cat$^1$-group $\mathfrak{C}$ is the 2-functor $$\mathbf{\lambda}:\mathfrak{C}^{op}\rightarrow \mathbf{Ch}^{1}_{K}$$ which is sending

$\bullet$ the 0-cell $*$ in  $\mathfrak{C}$ to the chain complex of vector spaces of length-1; $$\overline{\delta}=\overline{\tau}|_{\ke\overline{\sigma}}:\ke\overline{\sigma}\longrightarrow K(N),$$

$\bullet$ each $n\in N$, the 1-cell in $\mathfrak{C}$, to the chain automorphism $\mathbf{\lambda}_n:=(\lambda^{0}_n,\lambda^{1}_n):\overline{\delta}\longrightarrow \overline{\delta}$ defined on each level by
$$
\lambda^{0}_n(\mathbf{e}_{n'}):=\mathbf{e}_{n'n},  \hspace{1cm} \ \ \lambda^{1}_{n}(\mathbf{\overline{v}}_{m',n'}):=\mathbf{\overline{v}}_{m',n'n}
$$
where $\lambda^{0}_n:K(N)\rightarrow K(N)$ and $\lambda^{1}_n:\ke\overline{\sigma}\rightarrow \ke\overline{\sigma}$ are linear automorphisms,

$\bullet$ each $(m,n)\in M\rtimes N$, the 2-cell in $\mathfrak{C}$, to the homotopy
$$\mathbf{\lambda}_{m,n}:=(\lambda'_{(m,n)},\lambda_n):\lambda_n\Rightarrow \lambda_{\partial mn}$$ with the chain homotopy component $\lambda'_{(m,n)}:K(N)\longrightarrow \ke\overline{\sigma}$ defined by $\lambda'_{(m,n)}(\mathbf{e}_{n'}):=\mathbf{\overline{v}}_{^{n'}m,n'n}$ and where all chain automorphisms $\mathbf{\lambda}_n$ and homotopies $\mathbf{\lambda}_{m,n}$ resides in \textbf{Aut}$(\overline{\delta})$ for the linear transformation $\overline{\delta}:=\overline{\tau}|_{\ke\overline{\sigma}}$ obtained from the cat$^1$-algebra $\overline{K(\mathfrak{C})}$ of $\mathfrak{C}$.
}
\end{defn}
The construction can be applied to any cat$^1$-group, so it gives us a cat$^1$-group version of Cayley's theorem.

\section{The Category Ch$^2_K$ as a Gray Category}
Kamps and Porter, \cite{KP}, have constructed a Gray category structure on the category $\mathbf{Ch}^2_K$; of chain complexes of length-2. They have also given the relation between 2-crossed modules and Gray 3-groupoids. Using a slightly different language, Martins and Picken, \cite{Martins}, have proven this relationship between these structures. To define the linear representation of a 2-crossed module, Al-asady, in \cite{Jinan}, has considered these relationships and defined a 3-functor
$$
\phi:\mathfrak{C}^2\longrightarrow \mathbf{Ch}^2_K
$$
where $\mathfrak{C}^2$ is a Gray 3-groupoid with a single object set $\{*\}$ obtained from a 2-crossed module. In $\mathfrak{C}^2$, each object set is a group, so we can say that $\mathfrak{C}^2$ is a Gray 3-(group)-groupoid with a single object.

\subsection{Gray 3-Groupoids}

Kamps and Porter include in the Appendix of their work \cite{KP}, a sketch of Crans definition of Gray categories as algebraic structures (cf. \cite{Cr}). In this description, for convenience, they have inverted the interchange 3-cell. Martins and Picken, \cite{Martins}, gave a full definition of a Gray 3-groupoid. Their conventions are slightly different from the ones of \cite{Cr},\cite{KP}. The following definition is equivalent to that given by Martins and Picken in \cite{Martins}.

A (small) Gray 3-groupoid  $\mathfrak{C}^2$ is given by a set $C_0$ of objects (or 0-cells), a set $C_1$ of 1-cells, a set $C_2$ of 2-cells and a set $C_3$ of 3-cells, and maps $s_i,t_i:C_k\rightarrow C_{i-1}$, for $i=1,\cdots,k$ such that:
\begin{enumerate}
  \item $s_2\circ s_3=s_2$ and $t_2\circ t_3=t_2$, as maps $C_3\rightarrow C_1$.
  \item $s_1=s_1\circ s_2=s_1\circ s_3$ and $t_1=t_1\circ t_2=t_1\circ t_3$, as maps $C_3\rightarrow C_0$.
  \item $s_1=s_1\circ s_2$ and $t_1=t_1\circ t_2$, as maps $C_2\rightarrow C_0$.
  \item There exists a 2-vertical composition $J\#_3 J'$ of 3-cells if $s_3(J)=t_3(J')$, making it a groupoid with the set of objects is $C_2$ are implicit).
  \item There exists a vertical composition $$ \Gamma' \#_{2}\Gamma=\begin{bmatrix} \Gamma\\ \Gamma' \end{bmatrix}$$ of 2-cells if $s_2(\Gamma')=t_2(\Gamma)$, making it a groupoid whose set of objects is $C_1$ (identities are implicit).
  \item There exists a 1-vertical composition  $$ J' \#_{1}J=\begin{bmatrix} J\\ J' \end{bmatrix}$$
      of 3-cells if $s_2(J')=t_2(J)$ making the set of 3-cells $C_3$ a groupoid with set of objects $C_1$ and such that the boundaries $s_3,t_3:C_3\rightarrow C_2$ are functors (groupoid morphisms).
  \item The 1-vertical and 2-vertical compositions of 3-cells satisfy the interchange law:
$$
(J'\#_3 J)\#_1(J'_1\#_3 J_1)=(J'\#_1 J'_1)\#_3 (J\#_1 J_1).
$$
Combining with the previous axioms, this means that the 1-vertical and 2-vertical compositions of 3-cells and the vertical composition of 2-cells give $C_3$ the structure of a 2-groupoid (cf. \cite{HKK}), with set of objects being $C_1$, set of 1-cells $C_2$ and set of 2-cells $C_3$.
\item (\textbf{Whiskering by 1-cells}) For each $x,y\in C_0$, it can be defined a 2-groupoid $\mathfrak{C}(x,y)$ of all 1-, 2- and 3-cells $b$ such that $s_1(b)=x$ and $t_1(b)=y$. Given a 1-cell $\eta:y\rightarrow z$, there is a 2-groupoid map $\natural_1\eta:\mathfrak{C}(x,y)\rightarrow \mathfrak{C}(y,z)$. Similarly if $\eta':w\rightarrow x$, there is a 2-groupoid map $\eta'\natural_1:\mathfrak{C}(x,y)\rightarrow \mathfrak{C}(w,y)$.
\item There exists a horizontal composition $\eta \natural_1 \eta'$ of 1-cells if $s_1(\eta)=t_1(\eta')$, which is to be associative and to define a groupoid with set of objects $C_0$ and set of 1-cells $C_1$.

\item Given $\eta,\eta'\in C_1$;
\begin{align*}
\natural_1\eta \circ \natural_1 \eta'=&\natural_1(\eta'\eta)\\
\eta \natural_1\circ \eta'\natural_1=& (\eta \eta')\natural_1\\
\eta \natural_1\circ \natural_1 \eta'=&\natural_1 \eta'\circ \eta \natural_1,
\end{align*}
whenever these compositions make sense.
\item There are two horizontal compositions of 2-cells
\begin{align*}
               \begin{bmatrix} &\Gamma'\\ \Gamma& \end{bmatrix}
                                =&\left(\Gamma\natural_1t_1\left(\Gamma^\prime\right)\right)\#_2\left(s_1\left(\Gamma\right)\natural_1\Gamma^\prime\right)
\intertext{ and }
 \begin{bmatrix} \Gamma&\\ &\Gamma' \end{bmatrix}=&\left(t_1\left(\Gamma\right)\natural_{1}\Gamma^\prime\right)\#_2\left(\Gamma\natural_1s_1\left(\Gamma^\prime\right)\right)\\
\end{align*}
and of 3-cells:
\begin{align*}
               \begin{bmatrix} &J'\\ J& \end{bmatrix}
                                =&\left(J\natural_1t_1\left(J'\right)\right)\#_1\left(s_1\left(J\right)\natural_1 J'\right)
\intertext{ and }
 \begin{bmatrix} J&\\ &J' \end{bmatrix}=&\left(t_1\left(J\right)\natural_{1}J'\right)\#_1\left(J\natural_1s_1\left(J'\right)\right)\\
\end{align*}
It follows from the previous axioms that they are associative.
\item (\textbf{Interchange 3-cells}) For any 2-cells $\Gamma$ and $\Gamma'$, there is a 3-cell (called an interchange 3-cell)
$$
 \xymatrix{{\begin{bmatrix} &\Gamma'\\ \Gamma & \end{bmatrix}}=s_3(\Gamma\#\Gamma')\ar[rr]^{\scriptstyle(\Gamma\#\Gamma')} & &t_3(\Gamma\#\Gamma')={\begin{bmatrix} \Gamma &\\ & \Gamma' \end{bmatrix}}}
$$
\item(\textbf{2-functoriality}) For any 3-cells $\xymatrix{ \Gamma_1=s_3(J)\ar[r]^{\scriptstyle J} &t_3(J)=\Gamma_2}$ and $\xymatrix{ \Gamma'_1=s_3(J')\ar[r]^{\scriptstyle J'} &t_3(J')=\Gamma'_2}$, with $s_1(J')=t_1(J)$ the following upwards compositions (1-vertical compositions) of 3-cells coincide:
$$
 \xymatrix{{\begin{bmatrix} &\Gamma'_1\\ \Gamma_1 \end{bmatrix}}\ar[rr]^{\scriptstyle(\Gamma_1\#\Gamma'_1)} & &{\begin{bmatrix} \Gamma_1&\\ &\Gamma'_1 \end{bmatrix}}\ar[rr]^{\begin{bsmallmatrix} J&\\ &J' \end{bsmallmatrix}} & &{\begin{bmatrix} \Gamma_2&\\ &\Gamma'_2 \end{bmatrix}}}
$$
and
$$
 \xymatrix{{\begin{bmatrix} &\Gamma'_1\\ \Gamma_1 \end{bmatrix}}\ar[rr]^{\begin{bsmallmatrix} &J'\\ J& \end{bsmallmatrix}} & &{\begin{bmatrix} &\Gamma'_2\\ \Gamma_2& \end{bmatrix}}\ar[rr]^{\scriptstyle(\Gamma_2\#\Gamma'_2)} & &{\begin{bmatrix} \Gamma_2&\\ &\Gamma'_2 \end{bmatrix}}}
$$
This of course means that the collection $\Gamma\#\Gamma'$, for arbitrary 2-cells $\Gamma$ and $\Gamma'$ with $s_1(\Gamma')=t_1(\Gamma)$ defines a natural transformation between the 2-functors of 11. Note that by using the interchange condition for the vertical and upwards compositions, we only need to verify this condition for the case when either $J$ or $J'$ is an identity.(This is the way this axiom appears written in \cite{KP,Cr,Be}

\item(\textbf{1-functoriality}) For any three 2-cells $\gamma \xrightarrow{\Gamma \ } \phi \xrightarrow{\Gamma'}\psi$ and $\gamma'' \xrightarrow{\Gamma''} \phi''$ with $s_2(\Gamma')=t_2(\Gamma)$ and $t_1(\Gamma)=t_1(\Gamma')=s_1(\Gamma'')$ the following 1-vertical compositions of 3-cells coincide:
$$
 \xymatrix{{\begin{bmatrix} \psi \natural_1&\Gamma''\\ \Gamma'& \natural_1\gamma''\\ \Gamma&\natural_1\gamma'' \end{bmatrix}}\ar[rr]^{\begin{bsmallmatrix} \Gamma'\#\Gamma''\\ \Gamma\natural_1\gamma'' \end{bsmallmatrix}} & &{\begin{bmatrix} \Gamma'&\natural_1\phi''\\ \phi\natural_1& \Gamma''\\ \Gamma&\natural_1\gamma'' \end{bmatrix}}\ar[rr]^{\begin{bsmallmatrix} \Gamma'\natural_1\phi''\\ \Gamma\#\Gamma'' \end{bsmallmatrix}} & &{\begin{bmatrix} \Gamma'&\natural_1\phi''\\\Gamma& \natural_1\phi''\\  \gamma\natural_1&\Gamma'' \end{bmatrix}}}
$$
and
$$
 \xymatrix{{\begin{bmatrix} \psi \natural_1&\Gamma''\\ \Gamma'& \natural_1\gamma''\\ \Gamma&\natural_1\gamma'' \end{bmatrix}}\ar[rr]^{\begin{bsmallmatrix} \Gamma'\\ \Gamma \end{bsmallmatrix}\#\Gamma''} & &{\begin{bmatrix} \Gamma'&\natural_1\phi''\\\Gamma& \natural_1\phi''\\  \gamma\natural_1&\Gamma'' \end{bmatrix}}}
$$
An analogous identity obtained by exchanging the roles of the first and second columns should also hold.
A Gray 3-(group)-groupoid is a Gray 3-groupoid in which each objects set, $C_0,C_1,C_2$ and $C_3$ are groups and the structural maps $s_i,t_i$ are homomorphisms of groups.
We can show a Gray 3-(group)-groupoid with a single 0-cell $*$ pictorially as;
$$
\begin{array}{c}
  \boldsymbol {\mathfrak{C}^2}
\end{array}
{:=}
\begin{array}{c}
\xymatrix@C+=1.5cm{
& C_3\ar[ldd]_{\scriptscriptstyle s_3,t_3}\ar[rdd]^{\scriptscriptstyle s_1,t_1}\ar[d] & \\
& {*}  & \\
 C_2 \ar[rr]_{\scriptscriptstyle s_2,t_2}\ar[ur]  & &  C_1. \ar[ul] }
\end{array}
$$
\end{enumerate}
\begin{defn}
\rm{(\cite{Martins}) A (strict) Gray functor $\mathcal{F}:\mathfrak{C}^2\rightarrow\mathfrak{C'}^2$ between Gray 3-groupoids $\mathfrak{C}^2$ and $\mathfrak{C'}^2$ is given by maps $C_i\rightarrow C'_i$ preserving all compositions, identities and boundaries, strictly.
}
\end{defn}

We now recall the construction from \cite{Jinan} and \cite{KP}  of a Gray category structure from chain complexes length-2 of vector spaces over $K$. Let $C_2$, $C_1$ and $C_0$ be vector spaces over $K$ and $\delta^C_2:C_2\rightarrow C_1$,\ \ $\delta^C_1:C_1\rightarrow C_0$ be linear transformations.
Suppose that
$$
\begin{tikzcd}
             \mathcal{C}:=C_{2}\ar[r,"{\delta^C_2}"]&C_{1}\ar[r,"{\delta^C_{1}}"]& C_{0}
\end{tikzcd}
$$
is any chain complex of length-2. Then, 0-cells of  $ \mathbf{Ch}^2_K$ will be considered as chain complexes of length-2.

A chain map $F=(f_2,f_1,f_0)$ from a chain complex $\mathcal{C}$ to a chain complex $\mathcal{D}$ is given by the following commutative diagram
$$
\begin{array}{c}
 \xymatrix{\mathcal{C}\ar[d]_{F}\\
 \mathcal{D}}
 \end{array}
 {:=}
 \begin{array}{c}
   \xymatrix{C_{2} \ar[r]^-{\delta^C_2}\ar[d]^-{f_{2}}& C_{1} \ar[r]^-{\delta^C_{1}}\ar[d]^-{f_1}&C_{0}\ar[d]^-{f_{0}}\\
D_{2} \ar[r]_-{\delta^D_2} & D_{1} \ar[r]_-{\delta^D_1}  &D_{0}  }
\end{array}
$$
where $f_2,f_1$ and $f_0$ are linear transformations. Then, 1-cells of  $ \mathbf{Ch}^2_K$ will be considered as all chain maps between chain complexes of length-2.
Let $F$ and $G$ be chain maps between $\mathcal{C}$ and $\mathcal{D}$. A 1-homotopy $(H,F):F\Rightarrow G$ with the chain homotopy components $H'_1$ and $H'_2$, can be represented by the diagram
\begin{equation*}
\xymatrix{C_{2} \ar[rr]^-{\delta^C_2}\ar@<0.5ex>[dd]^-{g_{2}}\ar@<-0.5ex>[dd]_-{f_{2}} && C_{1} \ar[ddll]|{H'_2} \ar[rr]^-{\delta^C_1}\ar@<0.5ex>[dd]^-{g_1}\ar@<-0.5ex>[dd]_-{f_1}&&C_{0} \ar[ddll]|{H'_{1}}\ar@<0.5ex>[dd]^-{g_{0}}\ar@<-0.5ex>[dd]_-{f_{0}}
\\
\\
D_{2}\ar[rr]_{\delta^D_2} && D_{1} \ar[rr]_{\delta^D_1} \ar[rr] &&D_{0}}
\end{equation*}
where the chain homotopy components $H'_1:C_0\rightarrow D_1$ and $H'_2:C_1\rightarrow D_2$ are linear transformations, satisfying the following conditions:
\begin{enumerate}
\item $\delta^D_1 H'_1=g_0-f_0$,
\item $H'_1\delta^C_1+\delta^D_2H'_2=g_1-f_1,$
\item $H'_2\delta^C_2=g_2-f_2.$
\end{enumerate}
Then, 2-cells of  $ \mathbf{Ch}^2_K$ will be considered as all 1-homotopies between chain maps.

Let $(H,F),(K,F):F\Rightarrow G$ be 1-homotopies between $F$ and $G$. A 2-homotopy $\alpha:=(\alpha',H,F):(H,F)\Rrightarrow (K,F)$ with the chain homotopy component $\alpha':C_0\rightarrow D_2$ can be given by the diagram
\begin{equation*}
\xymatrix{C_{2} \ar[rr]^-{\delta^C_2}\ar@<0.5ex>[dd]\ar@<-0.5ex>[dd] && C_{1} \ar@<-0.6ex>[ddll]|<<<<<<<<<<<<<<<<{\scriptscriptstyle{H'_2}} \ar@<0.6ex>[ddll]|<<<<<<<<{\scriptscriptstyle{K'_2}} \ar[rr]^-{\delta^C_1}\ar@<0.5ex>[dd]\ar@<-0.5ex>[dd]&&C_{0}\ar@{-->}@/^{0.5pc}/[ddllll]_(.4){\large\boldsymbol\alpha'} \ar@<-0.6ex>[ddll]|<<<<<<<<<<<<<<<<{\scriptscriptstyle{H'_1}} \ar@<0.6ex>[ddll]|<<<<<<<<{\scriptscriptstyle{K'_1}}\ar@<0.5ex>[dd]\ar@<-0.5ex>[dd]
\\
\\
D_{2}\ar[rr]_{\delta^D_2} && D_{1} \ar[rr]_{\delta^D_1} \ar[rr] &&D_{0}}
\end{equation*}
where for the homotopy component $\alpha'$, the conditions $\delta^D_2\alpha'=K'_1-H'_1$ and $\alpha'\delta^C_1=K'_2-H'_2$ are satisfied. Therefore, we can illustrate 0-,1-,2- and 3-cells in $ \mathbf{Ch}^2_K$ by a diagram
$$
  \begin{tikzcd}[row sep=0.3cm,column sep=scriptsize]
                & \ar[dd, Rightarrow, "{(H,F)}"{swap,name=f,description}]& &\ar[dd, Rightarrow, "{(K,F)}"'{swap,name=g,description}]& \\
              {\mathcal{C}\ }\ar[rrrr,bend left=40,"F"] \ar[rrrr,bend right=40,"{G}"']& \tarrow["{(\alpha', H,F)}" ,from=f,to=g, shorten >= -1pt,shorten <= 1pt ]{rrr}&  & & {\ \mathcal{D}}. \\
                \ & \ & \ & \ & \
    \end{tikzcd}
$$
A 3-cell $\alpha:=(\alpha',H,F):(H,F)\Rrightarrow (K,F)$ may be written briefly as $\alpha:H\Rrightarrow K$. The source and target maps are given by
$s_3(\alpha)=(H,F)=\left((H'_1,H'_2),(f_2,f_1,f_0)\right)$ and $t_3(\alpha)=(K,F)=\left((K'_1,K'_2),(f_2,f_1,f_0)\right)$ where $K'_1=H'_1+\delta^D_2\alpha'$ and  $K'_2=H'_2+\alpha'\delta^C_1$. The vertical composition of 3-cells $\alpha:=(\alpha',H,F)$ and $\beta:=(\beta',K,F)$ is defined by
$$
\begin{bmatrix} \alpha\\ \beta \end{bmatrix}=\beta\#_3\alpha:=(\beta'+\alpha',H,F)
$$
with $t_3(\alpha)=s_3(\beta).$ The other source and target maps for a 3-cell $\alpha$ are given by $s_2(\alpha)=F, \ \ t_2(\alpha)=G$, and  $s_1(\alpha)=\mathcal{C},\ \ t_1(\alpha)=\mathcal{D}$. The composition of 1-cells is the usual composition of the chain maps. For  1-cells $F,G,T:\mathcal{C}\rightarrow \mathcal{D}$, the vertical composition of 2-cells $(H,F):F\Rightarrow G$ and $(K,G):G\Rightarrow T$ is given by $K\#_2 H:F\Rightarrow T$, where the chain homotopy component is $(K\#_2 H)'=K'+H'$, where $K'=(K'_1,K'_2)$ and $H'=(H'_1,H'_2).$

For any 2-cells
$$\Gamma=(K,G)=\left((K'_1,K'_2),(G_2,G_1,G_0)\right):G\Rightarrow G'$$
and
$$\Gamma'=(H,F)=\left((H'_1,H'_2),(F_2,F_1,F_0)\right):F\Rightarrow F'$$
the whiskering of $F'$ on $(K,G)$ is given by
$$(K,G)\natural_1 F'=(K'_1F'_0,K'_2F'_1,G\#_0F')$$
and this can be represented pictorially as;
$$
\begin{array}{c}
    \begin{tikzcd}[row sep=small,column sep=0.8cm]
                    &  &    \ar[dd,Rightarrow,"\scriptscriptstyle{(K,G)}"{description}] & \\
                 \mathfrak{C} \ar[r,"\scriptscriptstyle F'"] &\mathfrak{D}\ar[rr,bend left=50,"\scriptscriptstyle G"] \ar[rr,bend right=50,"\scriptscriptstyle G'"'] &  & \mathcal{E} \\
                \ & \ & \  & \
    \end{tikzcd}
\end{array}
{:=}
\begin{array}{c}
    \begin{tikzcd}[row sep=small,column sep=0.8cm]
                 & \ar[dd, Rightarrow, "\scriptscriptstyle{(K,G)\natural_1 F'}"{description}] & \\
                 \mathfrak{C}\ar[rr,bend left=50,"\scriptscriptstyle G\#_0 F'"] \ar[rr,bend right=50,"\scriptscriptstyle G'\#_0 F'"'] &  &\mathcal{E}.\\
                 & \
    \end{tikzcd}
\end{array}
$$
The whiskering of $G$ on $(H,F)$ is given by
$$G\natural_1(H,F)=(G_1H'_1,G_2H'_2,G\#_0F)$$
as illustrated in the following diagram;
$$
\begin{array}{c}
    \begin{tikzcd}[row sep=small,column sep=0.8cm]
                 &\ar[dd,Rightarrow,"\scriptscriptstyle{(H,F)}"description]&\\
                 \mathfrak{C}\ar[rr,bend left=50,"\scriptscriptstyle F"]\ar[rr,bend right=50,"\scriptscriptstyle {F'}"']& & \mathfrak{D}\ar[r,"\scriptscriptstyle G"]& \mathcal{E}\\
                 &\
    \end{tikzcd}
\end{array}
{:=}
\begin{array}{c}
    \begin{tikzcd}[row sep=small,column sep=0.8cm]
                  &\ar[dd,Rightarrow,"{\scriptscriptstyle G\natural_1 (H,F)}"description]&\\
                   \mathfrak{C} \ar[rr,bend left=50,"\scriptscriptstyle G\#_0 F"]\ar[rr,bend right=50,"\scriptscriptstyle G\#_0 F'"']& &  \mathcal{E}.\\
                  &\
    \end{tikzcd}
\end{array}
$$
The horizontal compositions of $\Gamma$ and $\Gamma'$ can be given by the vertical composition of these two 2-cells;
$$
\begin{bmatrix} &\Gamma'\\ \Gamma& \end{bmatrix}=(\Gamma \natural_1 t_1(\Gamma'))\#_2 (s_1(\Gamma) \natural_1 \Gamma')=(K'_1F'_0 + G_1H'_1,K'_2F'_1 + G_2H'_2, G\#_0 F).
$$
and
$$
 \begin{bmatrix} \Gamma&\\ &\Gamma' \end{bmatrix}=(t_1(\Gamma)\natural_1\Gamma')\#_2 (\Gamma\natural_1 s_1(\Gamma'))=(K'_1F_0 + G'_1H'_1,K'_2F_1 + G'_2H'_2, G\#_0 F).
$$
For any 3-cells:
$$J=(\beta',K,G):\Gamma_1=(K,G)\Rrightarrow\Gamma_2=(K',G)$$ and;
$$J'=(\alpha',H,F):\Gamma'_1=(H,F)\Rrightarrow\Gamma'_2=(H',F)$$
The horizontal composition of $J$ and $J'$
$$
\begin{bmatrix} &J' \\ J& \end{bmatrix}:\begin{bmatrix} &\Gamma'_1 \\ \Gamma_1& \end{bmatrix}\longrightarrow\begin{bmatrix} &\Gamma'_2 \\ \Gamma_2& \end{bmatrix};
$$
with $s_1(J')=D=t_1(J)$ can be defined by
$$
\begin{bmatrix} &J' \\ J& \end{bmatrix}=\left(G_2\alpha'+\beta'F'_0,\left(K'_1F'_0+G_1H'_1,K'_2F'_1+G_2H'_2\right),G\#_0F\right).
$$
In these morphisms, we have
$$
\begin{bmatrix} &\Gamma'_1 \\ \Gamma_1& \end{bmatrix}=\left(\left(K'_1F'_0+G_1H'_1,K'_2F'_1+G_2H'_2\right),G\#_0F \right)=s_3\left(\begin{bmatrix} &J' \\ J& \end{bmatrix}\right)\\
$$
and
$$
\begin{bmatrix} &\Gamma'_2 \\ \Gamma_2& \end{bmatrix}=\left(\left(K''_1F'_0+G_1H''_1,K''_2F'_1+G_2H''_2\right),G\#_0F \right)=t_3\left(\begin{bmatrix} &J' \\ J& \end{bmatrix}\right)\\
$$
where
\begin{align*}
K'_1F'_0+G_1H'_1+\delta_2^E(G_2\alpha'+\beta'F'_0)=&K'_1F'_0+G_1H'_1+\delta_2^EG_2\alpha'+\delta_2^E\beta'F'_0\\
                                                  =&K''_1F'_0+G_1H''_1 \hspace{1.0cm} (\because K''_1=\delta_2^E\beta'+K'_1 \ ,\ H''_1=\delta_2^D\alpha'+H'_1)
\end{align*}
and
\begin{align*}
K'_2F'_1+G_2H'_2+(G_2\alpha'+\beta'F'_0)\delta_1^C=&K'_2F'_1+G_2H'_2+G_2\alpha'\delta_1^C+\beta'F'_0\delta_1^C\\
                                                  =&K''_2F'_1+G_2H''_2
\end{align*}
Similarly, the horizontal composition
$$
\begin{bmatrix} J& \\ &J' \end{bmatrix}:\begin{bmatrix} \Gamma_1& \\ &\Gamma'_1\end{bmatrix}\longrightarrow\begin{bmatrix} \Gamma_2& \\ &\Gamma'_2 \end{bmatrix} $$
can be defined. Therefore, $\mathbf{Ch}^2_K$ has a Gray category structure.

For the calculations about 1-homotopies, along with chain complexes of length-2, and 2-homotopies between 1-homotopies and other stuff, the following diagram can be regarded as our \textit{black box} including the necessary tools to see the details.
$$
\xymatrix{C_2\ar@<0.6ex>[rr]|{\scriptscriptstyle F_2}\ar@<-0.6ex>[rr]|<<<<<{\scriptscriptstyle F'_2}\ar[dd]_{\scriptscriptstyle \delta^C_2}&& D_2\ar@<0.6ex>[rr]|{\scriptscriptstyle G_2}\ar@<-0.6ex>[rr]|<<<<<{\scriptscriptstyle G'_2}\ar[dd]^{\scriptscriptstyle\delta^D_2}&& E_2\ar[dd]^{\scriptscriptstyle \delta^E_2}\\
\\
C_1\ar@<0.6ex>[rr]|{\scriptscriptstyle F_1}\ar@<-0.6ex>[rr]|<<<<<{\scriptscriptstyle F'_1}\ar[dd]_{\delta^C_1}\ar@<0.6ex>[uurr]|{\scriptscriptstyle H'_2}\ar@<-0.6ex>[uurr]|<<<<<<<{\scriptscriptstyle H''_2}&& D_1\ar[dd]^{\delta^D_1}\ar@<0.6ex>[rr]|{\scriptscriptstyle G_1}\ar@<-0.6ex>[rr]|<<<<<{\scriptscriptstyle G'_1}\ar@<0.6ex>[uurr]|{\scriptscriptstyle K'_2}\ar@<-0.6ex>[uurr]|<<<<<<<{\scriptscriptstyle K''_2} &&E_1\ar[dd]^{\delta^E_1}\\
\\
C_0\ar@<0.6ex>[rr]|{\scriptscriptstyle F_0}\ar@<-0.6ex>[rr]|<<<<<{\scriptscriptstyle F'_0}\ar@<0.6ex>[uurr]|{\scriptscriptstyle H'_1}\ar@<-0.6ex>[uurr]|<<<<<<<{\scriptscriptstyle H''_1}\ar@{-->}@/^{-0.5pc}/[uuuurr]_(.7){\large\boldsymbol\alpha'}&&D_0\ar@<0.6ex>[rr]|{\scriptscriptstyle G_0}\ar@<-0.6ex>[rr]|<<<<<{\scriptscriptstyle G'_0}\ar@<0.6ex>[uurr]|{\scriptscriptstyle K'_1}\ar@<-0.6ex>[uurr]|<<<<<<<{\scriptscriptstyle K''_1}\ar@{-->}@/^{-0.5pc}/[uuuurr]_(.7){\large\boldsymbol\beta'}&&E_0.}
$$

\subsection{The Structure Aut($\delta$) in $\mathbf{Ch}^2_K$ as a Cat$^2$-Group}

To define linear and matrix representations of a cat$^2$-group $\mathfrak{C}^2$ as a 3-functor $\phi: \mathfrak{C}^2\rightarrow \mathbf{Ch}^2_K$, Al-asady, \cite{Jinan}, has constructed an automorphism cat$^2$-group $\mathbf{Aut(\delta)}$ as a Gray 3-groupoid with a single 0-cell or a cat$^2$-group in $\mathbf{Ch}^2_K$. In this case, the functorial image of $\mathfrak{C}^2$ under $\phi$ will be $\mathbf{Aut(\delta)}$. This structure is a collection of all chain automorphisms $\delta \rightarrow \delta$ and 1-homotopies between them and 2-homotopies between 1-homotopies. Since an isomorphism from an object to itself is known as automorphism, the structure $\mathbf{Aut(\delta)}$ is called an automorphism cat$^2$-group. The following definition is due to \cite{Jinan}.

\begin{defn}\rm{
  Let $$
\begin{tikzcd}
             \delta:=C_{2}\ar[r,"{\delta_2}"]&C_{1}\ar[r,"{\delta_{1}}"]& C_{0}
\end{tikzcd}
$$
be a chain complex of vector spaces over the field $K$ of length-2. The automorphism cat$^2$-group $\mathbf{Aut}(\delta)$ consists of:

$1.$ the group $\mathbf{Aut}(\delta)_0=\{\delta\}$, where $\delta:=(\delta_2,\delta_1)$ is the chain complex of length-2 given above,

$2.$ the group $\mathbf{Aut}(\delta)_1$ of chain automorphisms $F:=(F_2,F_1,F_0):\delta\rightarrow\delta$,

$3.$ the group $\mathbf{Aut}(\delta)_2$ of all 1-homotopies $(H,F):=((H'_1,H'_2),F):F\Rightarrow G$ between chain automorphisms,

$4.$ the group $\mathbf{Aut}(\delta)_3$ of all 2-homotopies $\alpha:=(\alpha',H,F):(H,F)\Rrightarrow (K,F)$ between 1-homotopies together with source,target and identity maps between these groups given in \cite{Jinan}.
}
\end{defn}


\section{2-Crossed Modules and Gray 3-Groupoids with a Single 0-cell}
In this section,  we first recall the construction of  a Gray 3-groupoid with a single object from a 2-crossed module (cf. \cite{KP}, \cite{Martins}, \cite{Jinan}) to use it for giving the regular representation in the next section.

Cat$^2$-groups, \cite{WL}, are higher dimensional cat$^1$-groups. They can be regarded as cat$^1$-objects in the category of cat$^1$-groups. A cat$^2$-group $\mathfrak{C}^2$ is a 5-tuple $(G,s_1,t_1,s_2,t_2)$ where $(G,s_i,t_i)$ \ $i=1,2$ are cat$^1$-groups and
\begin{enumerate}
\item $s_is_j=s_js_i,\ t_it_j=t_jt_i, \ s_it_j=t_js_i \ \text{for} \ i,j=1,2 \ \,i\neq j$
\item $[\ke s_i,\ke t_i]=1\ \text{for} \ i=1,2.$
\end{enumerate}

Recall from \cite{Con} that \emph{a 2-crossed module} of groups consists of a complex of groups
$$
\xymatrix{\mathfrak{X}:=(L\ar[r]^-{\partial_2}&M\ar[r]^{\partial_1}&N)}
$$
together with (a) actions of $N$ on $M$ and $L$ so that $\partial
_{2},\partial _{1}$ are morphisms of $N$-groups,  and (b) an
$N$-equivariant function
\begin{equation*}
\{\quad ,\quad \}:M\times M\longrightarrow L
\end{equation*}%
called a Peiffer lifting. This data must  satisfy the following
axioms:
\begin{equation*}
\begin{array}{lrrll}
\mathbf{2CM1)} &  & \partial _{2}\{m,m^{\prime }\} & = & \left( ^{\partial
_{1}m}m^{\prime }\right) mm^{\prime }{}^{-1}m^{-1}\newline
\\
\mathbf{2CM2)} &  & \{\partial _{2}l,\partial _{2}l^{\prime }\} & = &
[l^{\prime },l]\newline
\\
\mathbf{2CM3)} &  & (i)\quad \{mm^{\prime },m^{\prime \prime }\} & = &
^{\partial _{1}m}\{m^{\prime },m^{\prime \prime }\}\{m,m^{\prime }m^{\prime
\prime }m^{\prime }{}^{-1}\}\newline
\\
&  & (ii)\quad \{m,m^{\prime }m^{\prime \prime }\} & = & \{m,m^{\prime
}\}^{mm^{\prime }m^{-1}}\{m,m^{\prime \prime }\}\newline
\\
\mathbf{2CM4)} &  & \{m,\partial _{2}l\}\{\partial _{2}l,m\} & = &
^{\partial _{1}m}ll^{-1}\newline
\\
\mathbf{2CM5)} &  & ^{n}\{m,m^{\prime }\} & = & \{^{n}m,^{n}m^{\prime }\}%
\newline
\end{array}%
\end{equation*}%
\newline
for all $l,l^{\prime }\in L$, $m,m^{\prime },m^{\prime \prime }\in M$ and $%
n\in N$.
Now suppose that $
\xymatrix{\mathfrak{X}:=(L\ar[r]^-{\partial_2}&M\ar[r]^{\partial_1}&N)}
$ is a 2-crossed module of groups with the Peiffer lifting $\{\quad ,\quad \}:M\times M\longrightarrow L$. We will remind ourselves of the connection between 2-crossed modules and Gray 3-groupoids with  a single 0-cell and cat$^2$-groups.

First, we define the groups of 0-cells, 1-cells, 2-cells and 3-cells. The group of 0-cells is $C_0=\{*\}$. The group of 1-cells is $C_1=N$. Then, a 1-cell in $\mathfrak{C}^2$ is an element $n\in N$. It will be considered as a morphism over $*$. The composition of 1-cells $n$ and $n'$ in $C_1$ is given by the operation of the group $N$. We can show a 1-cell $n\in C_1$ by $n:* \rightarrow *$. Using the group-action of $N$ on $M$, we can create the semi-direct product group $C_2=M \rtimes N$ together with the operation $ (m,n)(m',n')=(m^{n}m^{\prime},nn') $
for $m,m'\in M$ and $n,n'\in N$. An element $(m,n)$ of $C_2$ can be considered as a 2-cell from $n$ to $\partial_1 mn$, so we can define source, target maps between $C_2$ and $C_1$ as follows: for $(m,n)\in (M \rtimes N)=C_2$, the 1-source of this 2-cell is $n$ and so $s_2(m,n)=n$ and 1-target of this 2-cell is $t_2(m,n)=\partial_1 mn$. The 0-source and 0-target of $(m,n)$ is $*$. We can represent a 2-cell $(m,n)$ in $\mathfrak{C}^2$ pictorially as:
$$
\begin{array}{c}
    \xymatrix{ {*} \ar[r]\ar[d]_{\scriptscriptstyle n}^{\:\: \scriptscriptstyle(m,n)} &{*}\ar[d]^{ \scriptscriptstyle\partial_1 mn}\\
  {*}   \ar[r]  &  {*} }
  \end{array}
{:=}
  \begin{array}{c}
    \begin{tikzcd}[row sep=tiny,column sep=small]
                & \ar[dd, Rightarrow, "{\scriptscriptstyle(m,n)}"] \\
              {*}\ar[rr,bend left=50,"\scriptscriptstyle n"] \ar[rr,bend right=50,"\scriptscriptstyle \partial_1 mn"']  &  & \ {*} \\
                 & \
    \end{tikzcd}
  \end{array}
$$\\
The vertical composition of $\Gamma=\left(m,n\right)$ and $\Gamma'=\left(m',\partial_1 mn\right)$ in $C_2$ is given by
$$ \Gamma' \#_{2}\Gamma=\begin{bmatrix} \Gamma\\ \Gamma' \end{bmatrix}=\left(m',\partial_1 mn\right)\#_2(m,n)=(m'm,n).$$
To define the horizontal composition of 2-cells, we need to give the whiskering of a 1-cell on  a 2-cell on the right and left sides. The whiskering $n'\in C_1$ on $(m,n)\in C_2$ on the left side is $n'\natural_{1}(m,n)=(^{n'}m,n'n)$. Similarly the right whiskering of $n^\prime$ on $(m,n)$ is given by  $\left(m,n\right)\natural_1n^\prime=\left(m,nn^\prime\right).$ Therefore, the horizontal compositions of $\Gamma$ and $\Gamma^\prime$ are:
$$
               \begin{bmatrix} &\Gamma'\\ \Gamma& \end{bmatrix}
                                =(\Gamma\natural_1t_1(\Gamma^\prime))\#_2(s_1(\Gamma)\natural_1\Gamma^\prime)
                                =(m^nm^\prime,nn^\prime)
$$
and
$$
                \begin{bmatrix} \Gamma&\\ &\Gamma' \end{bmatrix}=(t_1(\Gamma)\natural_{1}\Gamma^\prime)\#_2(\Gamma\natural_1s_1(\Gamma^\prime))
                                      =(^{\partial_1m}(^nm^\prime)m,nn^\prime).
$$

Note that
$\begin{bmatrix} &\Gamma'\\ \Gamma& \end{bmatrix}\neq \begin{bmatrix} \Gamma&\\ &\Gamma' \end{bmatrix}$
since $\partial_1$ is not a crossed module.

We can show easily that $s_2$ and $t_2$ are  homomorphisms of groups from $C_1$ to $C_0$. For $\Gamma=(m,n)$ and $\Gamma^\prime=(m^\prime,n^\prime)$, we have;
$$s_2\left(\Gamma\#_0\Gamma^\prime\right)=s_2\left(m^nm^\prime,nn^\prime\right)=nn^\prime=s_2(\Gamma)s_2(\Gamma')$$
and
$$t_2\left(\Gamma\#_0\Gamma^\prime\right)=t_2\left(m^nm^\prime,nn^\prime\right)=\partial_1\left(m^nm^\prime\right)nn^\prime=\partial_1 mn\partial_1 m' n'=t_2(\Gamma)t_2(\Gamma').$$
Now, define the group of 3-cells in $\mathfrak{C}^2$. Using the group action of $M$ and of  $N$ on $L$, we can create the semi-direct product group  $C_3=L\rtimes M\rtimes N$ with the multiplication
$$(l,m,n)(l',m',n')=(l\{\partial_2(^{n}l'),m\}^{n}(l')^{-1},m^{n}m',nn')$$
where $\{-,-\}:M\times M\rightarrow L$ is the Peiffer lifting of the 2-crossed module $\mathfrak{X}$. Using the equality $\{\partial_2(l),m\}l^{-1}=^ml$, we can rewrite it
$$(l,m,n)(l',m',n')=(l^{m}(^{n}l'),m^{n}m',nn').$$

Any 3-cell in $C_3$ can be represented by an element $(l,m,n)$ in $L\rtimes M\rtimes N$ for $l\in L,\ m\in M, \  n\in N$. The 2-source of a 3-cell $(l,m,n)$ is given by $s_3(l,m,n)=(m,n)$ and 2-target is given by $t_3(l,m,n)=(\partial_2 lm,n).$ We can show a 3-cell in $C_3$ by a diagram;
$$
\begin{array}{c}
    \xymatrix@C-=0.5cm{{*} \ar@{=>}[rr]^{\scriptscriptstyle(m,n)}\ar[dd]_{\scriptscriptstyle n} & &{*}\ar[dd]^{\scriptscriptstyle \partial_1 mn}
\\&\scriptscriptstyle(l,m,n)& \\
  {*}  \ar@{=>}[rr]_{\scriptscriptstyle(\partial_2 lm,n)}  & & {*} }
  \end{array}
{:=}
  \begin{array}{c}
   \begin{tikzcd}[row sep=small,column sep=scriptsize]
                & \ar[dd, Rightarrow, "{\scriptscriptstyle (m,n)}" {swap,name=f,description}]& &\ar[dd, Rightarrow, "{\scriptscriptstyle (\partial_2 lm,n)}"'{swap,name=g,description}]& \\
              {*} \ar[rrrr,bend left=40,"\scriptscriptstyle n"] \ar[rrrr,bend right=40,"{\scriptscriptstyle \partial_1 mn}"']& \tarrow["\scriptscriptstyle {(l,m,n)}" ,from=f,to=g, shorten >= -1pt,shorten <= 1pt ]{rrr}&  & & {*} \\
                \ & \ & \ & \ & \
    \end{tikzcd}
  \end{array}
$$
The 2-vertical composition of 3-cells $J'=(l,m,n)$ and $J=(l^\prime,\partial_2lm,n)$ is \\
$$ J\#_{3}J' =\begin{bmatrix} J'\\ J \end{bmatrix}=(l^\prime,\partial_2lm,n)\#_3(l,m,n)=(l^\prime l,m,n).
$$
The right whiskering of a 2-cell $\Gamma=(m,n)$ on a 3-cell $J=\left(l,m^\prime,\partial_1 mn\right)$ is given by
$$\left(l,m^\prime,\partial_1 mn\right)\natural_2(m,n)=(l,m^\prime m,n).$$
This can be represented pictorially as
$$
\begin{array}{c}
    \begin{tikzcd}[row sep=small,column sep=normal]
                    &  & \tarrow["\scriptscriptstyle{(l,m',n')}" {description} ]{dd} & \\
                \scriptstyle n \ar[r,Rightarrow,"\scriptscriptstyle {(m,n)}"] &\scriptstyle {\partial_1 mn \ }\ar[rr,Rightarrow,bend left=40,"\scriptscriptstyle{ \left(m',n' \right)}"] \ar[rr,Rightarrow,bend right=50,"\scriptscriptstyle {\left(\partial_2 lm',n' \right)}"'] &  &\ \ \scriptstyle {\partial_1 m'n'}\\
                \ & \ & \
    \end{tikzcd}
\end{array}
{:=}
\begin{array}{c}
    \begin{tikzcd}[row sep=small,column sep=normal]
                 & \tarrow["\scriptscriptstyle{(l,m'm,n)}" {description} ]{dd}&\\
                 \scriptstyle {n} \ \ \ \ar[rr,Rightarrow,bend left=40,"\scriptscriptstyle { \left(m'm,n \right)}"] \ar[rr,Rightarrow,bend right=50,"\scriptscriptstyle {\left(\partial_2 lm'm,n \right)}"'] &  &  \scriptstyle {\partial_1 m'n'}\\
                 \ & \ &\
    \end{tikzcd}
\end{array}
$$
where $\partial_1 mn=n'.$

Similarly, 1-vertical composition in $\mathfrak{C}^2$ of 3-cells $(l,m,n)$ and $(l',m',\partial_1mn)$ is given by
$$(l,m,n)\#_1(l',m',\partial_1mn)=(l'^{m'}l,m'm,n).$$

The left whiskering of a 2-cell $\Gamma=(m',\partial_1 mn)$ on a 3-cell $J=(l,m,n)$ is given by
$$\left(m',\partial_1 mn\right)\natural_2\left(l,m,n\right)=(^{m'}l,m' m,n)$$
where \  $^{m^\prime}l=\{\partial_2l,m^\prime\}l$. This can be represented pictorially  as
$$
\begin{array}{c}
    \begin{tikzcd}[row sep=small,column sep=normal]
                 &\tarrow["\scriptscriptstyle{(l,m,n)}"{description}]{dd}&\\
                 \scriptstyle {\ \ n} \ \ \  \ar[rr,Rightarrow,bend left=40,"\scriptscriptstyle {(m,n)}"]\ar[rr,Rightarrow,bend right=50,"\scriptscriptstyle {(\partial_{2}lm,n)}"']& & \scriptstyle {\partial_{1} mn}\ar[r,Rightarrow,"\scriptscriptstyle {(m',n')}"]& \scriptstyle {\partial_{1}m'n'}\\
                 &\
    \end{tikzcd}
\end{array}
{:=}
\begin{array}{c}
    \begin{tikzcd}[row sep=small,column sep=normal]
                  &\tarrow["{\scriptscriptstyle(^{m^\prime}l,m^\prime m,n)}"{description}]{dd}&\\
                   \scriptstyle{n}\ \ \ \ar[rr,Rightarrow,bend left=40,"\scriptscriptstyle {(m'm,n)}"]\ar[rr,Rightarrow,bend right=50,"\scriptscriptstyle {(m'\partial_{2}lm,n)}"']& & \scriptstyle {\partial_{1}m'n'}\\
                  &\
    \end{tikzcd}
\end{array}
$$
where $\partial_1 mn=n'.$ The right whiskering of a 1-cell $n'$ on a 3-cell $(l,m,n)$ is $(l,m,n)\natural_1 n^\prime=(l,m,nn^\prime)$ and the left whiskering is $n^\prime\natural_1(l,m,n)=(^{n'}l,^{n'}m,n'n).$

The horizontal compositions of 3-cells  $J=(l,m,n):\Gamma_1\rightarrow \Gamma_2$ and $J^\prime=(l^\prime,m^\prime,n^\prime):\Gamma'_1\rightarrow \Gamma'_2$ in $C_3$  are given by
$$
\begin{bmatrix} &J' \\ J& \end{bmatrix}=((J\natural_1t_1(J^\prime))\#_1(s_1(J)\natural_1J^\prime)
                                       =(l{^m(^nl^\prime)},m^nm^\prime,nn^\prime)
$$
where
$${s_3}\left(\begin{bmatrix}  &J'\\ J& \end{bmatrix}\right)=\begin{bmatrix}  &s_3(J')\\ s_3(J)&  \end{bmatrix}=\begin{bmatrix}   &(m',n')\\ (m,n)&\end{bmatrix}=\begin{bmatrix}  &\Gamma'_1\\ \Gamma_1& \end{bmatrix}=\left(m^nm^\prime,nn^\prime\right)$$
and
$${t_3}\left(\begin{bmatrix}  &J'\\ J& \end{bmatrix}\right)=\begin{bmatrix}  &\Gamma'_2\\ \Gamma_2& \end{bmatrix}$$
similarly
$$
                       \begin{bmatrix} J & \\ & J' \end{bmatrix}
                                 =(t_1(J)\natural_1J^\prime)\#_1(J\natural_1s_1(J^\prime))
                                 =(^{\partial_1m}(^nl^\prime)^{^{\partial_1m}(^nm^\prime)}l,^{\partial_1m}(^nm^\prime)m,nn^\prime)
$$
where
$${s_3}\left(\begin{bmatrix}  J& \\ &J' \end{bmatrix}\right)=\begin{bmatrix}  s_3(J)&\\ &s_3(J')  \end{bmatrix}=\begin{bmatrix}   (m,n)&\\ &(m',n')\end{bmatrix}=\begin{bmatrix}  \Gamma_1&\\ &\Gamma'_1 \end{bmatrix}=\left(^{\partial_1m}(^nm^\prime)m,nn^\prime\right)$$
and
$${t_3}\left(\begin{bmatrix}  J& \\ &J' \end{bmatrix}\right)=\begin{bmatrix}  \Gamma_2&\\ &\Gamma'_2 \end{bmatrix}$$
Using the multiplication in $C_3$ given above, we can show that $s_3,t_3:C_3\longrightarrow C_2$ are homomorphisms of groups:
\begin{align*}
s_3\left((l,m,n)(l^\prime,m^\prime,n^\prime)\right)=&s_3\left(l \ {^m(^nl^\prime)},m^nm^\prime,nn^\prime\right)\\
                                                   =&(m^nm^\prime,nn^\prime)\\
                                                   =&(m,n)(m^\prime,n^\prime)\\
                                                   =&s_3(l,m,n) s_3(l^\prime,m^\prime,n^\prime)
\end{align*}
and
\begin{align*}
t_3\left((l,m,n)(l^\prime,m^\prime,n^\prime)\right)=&t_3\left(l \ {^m(^nl^\prime)},m^nm^\prime,nn^\prime\right)\\
                                                   =&\left(\partial_2\left(l \ {^m(^nl^\prime)}\right) m^nm^\prime,nn^\prime\right)\\
                                                   =&\left(\partial_2lm \partial_2(^nl^\prime)^nm^\prime,nn^\prime\right))\ \ (\partial_2 \text{ crossed module})\\
                                                   =&\left(\partial_2lm ^n(\partial_2l^\prime m^\prime),nn^\prime\right)\\
                                                   =&(\partial_2lm,n)(\partial_2l^\prime m^\prime,n^\prime)\\
                                                   =&t_3(l,m,n) t_3(l^\prime,m^\prime,n^\prime)
\end{align*}
for all  $(l,m,n),(l^\prime,m^\prime,n^\prime)\in C_3$. The interchange law for $\#_3$ and the semidirect product of 3-cells in $L\rtimes M\rtimes N$ can be found in Section \ref{interchangelaw} of Appendix. For any 2-cells $\Gamma=(m,n)$ and $\Gamma^\prime=(m^\prime,n^\prime)$, the interchange 3-cell is $\Gamma\#\Gamma^\prime=\left(\{m,^nm^\prime \},m^nm^\prime,nn^\prime\right)$ as given in \cite{KP}. For this interchange 3-cell, we have
$$s_3\left(\Gamma\#\Gamma^\prime\right)=\left(m^nm^\prime,nn^\prime\right)=\begin{bmatrix} &\Gamma'\\ \Gamma & \end{bmatrix}$$
and
\begin{align*}
t_3\left(\Gamma\#\Gamma^\prime\right)=&\left(\partial_2\{m,^nm^\prime\}m^nm^\prime,nn^\prime\right)\\
                                     =&\left(^{\partial_1(m)n}(m^\prime)m(^nm^\prime)^{-1}m^{-1}m^nm^\prime,nn^\prime\right)\\
                                     =&\left(^{\partial_1m}(^nm^\prime)m,nn^\prime\right)\\
                                     =&\begin{bmatrix} \Gamma&\\ &\Gamma' \end{bmatrix}
\end{align*}
Then we obtain
$$
    \xymatrix{ {\begin{bmatrix} &\Gamma'\\ \Gamma & \end{bmatrix}} \ar@3{->}[rr]^{\scriptstyle(\Gamma\#\Gamma')} & &{\begin{bmatrix} \Gamma&\\ &\Gamma' \end{bmatrix}}}
$$
Of course, this is functorial and hence defines a functor from the category of 2-crossed modules of groups to the category of Gray 3-(group)-groupoids with a single object: $\mathbf{\Theta}:\mathbf{X_2 Mod}\longrightarrow \mathbf{Gray}.$
\subsection{From Gray 3-(Group)-Groupoids to 2-Crossed Modules of Groups}
In this section, using a slightly different method, we give a construction of a 2-crossed module of groups from a Gray 3-(group)-groupoid with a single object. Similar construction appears in \cite{Martins}.

Suppose that $(C_3,C_2,C_1,s_1,t_1,s_2,t_2,*)$ is a Gray 3-(group)-groupoid with a single object. Using the properties of $s_i$,$t_i$, we have a complex of homomorphisms of groups:
\begin{equation*}
\xymatrix{\ke s_3\ar[r]^-{\overline{t}_3}&\ke s_2 \ar[r]^-{\overline{t}_2}&C_1 }
\end{equation*}
where $\overline{t}_2=t_2|_{\ke s_1}$ and $\overline{t}_3=t_3|_{\ke s_3 }$. The multiplication in $\ke s_2$ is taken to be the horizontal composition of 2-morphisms;
$$ \Gamma \Gamma'=\begin{bmatrix} &\Gamma'\\ \Gamma& \end{bmatrix}.$$
Similarly, the multiplication in $\ke s_3 $ is given by $ JJ'=\begin{bmatrix} &J'\\ J& \end{bmatrix}$ for $J,J'\in \ke s_3 $.
Using the whiskering $\natural_1$ of a 1-cell on 2-cell, an action of $\zeta\in C_1$ on a 2-cell \ $\Gamma\in \ke s_1$ is given by
$$^{\zeta}(\Gamma)=\zeta\natural_1\Gamma\natural_1\zeta^{-1}.$$
We have  $\overline{t}_2(\zeta \natural_1\Gamma)=\zeta\overline{t}_2(\Gamma)\zeta^{-1}$. Thus, $\overline{t}_2$ is a pre-crossed module. For any 3-cell $J \in \ke{s_3}$ and a 2-cell $\eta\in \ke s_1$, the action of $\eta$ on $J$ is given by
$$^{\eta}J=\eta\natural_2J\natural_2\eta^{-1}$$
where $\natural_2$ is the whiskering of the 2-cell $\eta$ on the 3-cell $J$. Together with this action $\overline{t}_3=t_3|_{\ke s_3 }$ is a crossed module.  The Peiffer Lifting map
 $$\{-,-\}:\ke s_2\times \ke s_2\longrightarrow \ke s_3$$
 is given by
 $$\{\Gamma,\Gamma^\prime\}=\begin{bmatrix} &e_3\left(s_3(\Gamma\#\Gamma^\prime)\right)^{-1}\\ \Gamma\#\Gamma^\prime & \end{bmatrix}$$
 where $\Gamma\#\Gamma^\prime$ is the interchange 3-cell for $\Gamma, \Gamma'\in \ke s_2$.

 This defines a functor from the category of Gray 3-(group)-groupoids with a single object to the category of 2-crossed modules of groups; $$\mathbf{\Delta}:\mathbf{Gray}\longrightarrow\mathbf{X_2 Mod}.$$

 Suppose that
 \begin{equation*}
\xymatrix{\mathfrak{X}:=L \ar[r]^-{\partial_{2}}&M\ar[r]^-{\partial_{1}} & N}
\end{equation*}
is a 2-crossed module of groups and its associated Gray 3-(group)-groupoid is $\mathbf{\Theta}(\mathfrak{X})=\mathfrak{C}^2$. By applying the functor $\mathbf{\Delta}$ to $\mathfrak{C}^2$, we will obtain a 2-crossed module
\begin{equation*}
\xymatrix{L' \ar[r]^-{\overline{\partial}_{2}}&M'\ar[r]^-{\overline{\partial}_{1}} & N'}
\end{equation*}
which is isomorphic to $\mathfrak{X}$ on each level. We obtain $N^\prime=C_1=N$ and  $M'=\ke s_2=M\rtimes\{1\}\cong M$ and $\overline{\partial}_1=t_2|_{\ke s_2}$, $\overline{\partial}_1(m,1)=\partial_1m$ for all $m\in M$.

By the definition of $s_3$, for $C_3=L\rtimes M\rtimes N$, we obtain
\begin{align*}
L'=\ke s_3=&\{\alpha=(l,m,n)\in C_3:s_3(\alpha)=(m,n)=(1,1)\}\\
=&\{(l,1,1):l\in L\}\cong L
\end{align*}
 The map $\overline{\partial}_2$ is given by $\overline{\partial}_2=t_3|_{\ke s_3}$ and then $\overline{\partial}_2(l,1,1)=t_3(l,1,1)=(\partial_2l,1)$. Therefore, we obtain a chain complex of homomorphisms of groups:
\begin{equation*}
\xymatrix{L\rtimes \{1\}\rtimes \{1\}\ar[r]^-{\overline{\partial}_{2}}&M\rtimes \{1\}\ar[r]^-{\overline{\partial}_{1}} & N}
\end{equation*}
which is isomorphic to $\mathfrak{X}$ on each level.\\
For any 2-cells $\Gamma=(m,1)\in \ke s_2$ and $\Gamma'=(m',1)\in \ke s_2$, We can define the Peiffer Lifting
$$\{-,-\}^\prime:M^\prime \times M^\prime\rightarrow L^\prime $$
by
\begin{align*}
\{\Gamma,\Gamma^\prime\}=&{\begin{bmatrix} &e_3\left(s_3(\Gamma\#\Gamma^\prime)\right)^{-1}\\ \Gamma\#\Gamma^\prime & \end{bmatrix}}\\
                        =&{\begin{bmatrix} &e_3s_3\left(\{m^\prime,m\},m^\prime m,1\right)^{-1}\\ \left(\{m^\prime,m\},m^\prime m,1\right) & \end{bmatrix}}\\
                        =&{\begin{bmatrix} &\left(1, (m^\prime m)^{-1},1\right)\\ \left(\{m^\prime,m\},m^\prime m,1\right) & \end{bmatrix}}\\
                        =&\left(\{m^\prime,m\},1,1\right)\in L'.
\end{align*}
\begin{remark}\rm{
Consider a 2-crossed module $\xymatrix{\mathfrak{X}:=(L\ar[r]^-{\partial_2}&M\ar[r]^{\partial_1}&N)}$ and its associated  Gray-3-(group)-groupoid $\mathfrak{C}^2$. From \cite{Jinan}, we have the construction of a linear representation of any cat$^2$-group or any Gray 3-(group)-groupoid with a single object, and this may also be considered as a linear representation of the corresponding 2-crossed module. Thus, we might define a linear representation of the 2-crossed module $\mathfrak{X}$ to be a representation of its associated Gray 3-(group)-groupoid $\mathfrak{C}^2$.\textit{ Thus we can say that a possible way of obtaining with a direct definition of a linear representation of the 2-crossed module $\mathfrak{X}$ would be to pass to the associated Gray 3-(group)-groupoid $\mathfrak{C}^2$ with a single object and find a representation for this.
}

 Therefore,  a linear representation of $\mathfrak{C}^2$ which is obtained from the 2-crossed module $\mathfrak{X}$ is a lax 3-functor $\phi:\mathfrak{C}^2\rightarrow \mathbf{Aut(\delta)}\leqslant \mathbf{Ch}^2_K$. This functor can be summarised as follows:

$\bullet$ For the 0-cell $*$ in $\mathfrak{C}^2$, $\phi(*)=\delta:=(\delta_2,\delta_1)$ is the chain complex of length-2, the 0-cell of $\mathbf{Aut(\delta)}$.

$\bullet$ For any 1-cell $n:* \longrightarrow *$, $n\in N$, $\phi(n)$  in $\mathbf{Aut}(\delta)_1$  is the chain automorphism $\delta\rightarrow \delta$

$\bullet$ For any 2-cell $(m,n):n\Rightarrow \partial_1 mn$ in $M\rtimes N$, \  $\phi(m,n)$  in $\mathbf{Aut}(\delta)_2$  is a 1-homotopy from $\phi(n)$ to $\phi(\partial_1 mn)$.

$\bullet$ For any 3-cell $(l,m,n):(m,n)\Rrightarrow (\partial_2 lm,n)$ in $L\rtimes M\rtimes N$,\  $\phi(l,m,n)$  in $\mathbf{Aut}(\delta)_3$ is a 2-homotopy between 1-homotopies $\phi(m,n)$ and $\phi(\partial_2 lm,n)$.}
\end{remark}
\section{Regular Representation of 2-Crossed Modules and Cat$^2$-groups}
 Consider the Gray 3-(group)-groupoid with a single object $*$;
 \begin{equation*}
\xymatrix{\mathfrak{C}^2:=L\rtimes M \rtimes N \ar@<0.5ex>[r]\ar@<-0.5ex>[r]& M\rtimes N\ar@<0.5ex>[r]\ar@<-0.5ex>[r] & N\ar@<0.5ex>[r]\ar@<-0.5ex>[r]& \{*\}}
\end{equation*}
which is obtained from the 2-crossed module $\mathfrak{X}$. We can consider $\mathfrak{C}^2$  as a cat$^2$-group. Then the big group is $L\rtimes M\rtimes N$. The structural homomorphisms are $s_1(l,m,n)=t_1(l,m,n)=*$ and $s_2(l,m,n)=(1,m,n)$, \ $t_2(l,m,n)=(1,\partial_2lm,n)$. Recall that the multiplication in $L\rtimes M\rtimes N$ is given by
 $$(l,m,n)(l^\prime,m^\prime,n^\prime)=\left(l^m(^{n}l^\prime),m^nm^\prime, nn^\prime\right).$$
Using this multiplication, we can show that $[\ke s_2,\ke t_2]=1$ as follows:

 For $(l,m,n)\in L\rtimes M\rtimes N,$  we obtain $s_2(l,m,n)=(1,m,n)=(1,1,1),$ so $(l,1,1)\in \ke s_2.$ Similarly,
 $t_2(l,m,n)=(1,\partial_2lm,n)=(1,1,1),$ then $m=\partial_2 l^{-1},\ n=1,$ and so $(l,\partial_2 l^{-1},1)\in \ke t_2.$  The inverse of $(l,m,n)$ is
 $$(l,m,n)^{-1}=(^{n^{-1}}{(^{m^{-1}}l^{-1})},^{n^{-1}}m^{-1},n^{-1}).$$
For $(l,1,1)\in \ke s_2$ and $(l^\prime,\partial_2 l^{\prime{-1}},1)\in \ke t_2,$ we obtain
\begin{align*}
  [(l,1,1),(l^\prime,\partial_2 l^{\prime{-1}},1)]=&(l,1,1)(l^\prime,\partial_2 l^{\prime{-1}},1)(l^{-1},1,1)(^{\partial_2l^{\prime}}(l^\prime)^{-1},\partial_2 l^{\prime},1)\\
  =&(ll^\prime,\partial_2 l^{\prime{-1}},1)(l^{-1} \ ^{\partial_2 l^{\prime}}(l^{\prime{-1}}),\partial_2 l^{\prime},1) \\
  =&(ll^\prime,\partial_2 l^{\prime{-1}},1)(l^{-1}l^\prime l^{\prime{-1}}l^{\prime{-1}},\partial_2 l^{\prime},1)\\
  =&(ll^\prime,\partial_2 l^{\prime{-1}},1)(l^{-1}l^{\prime{-1}},\partial_2 l^{\prime},1)\\
  =&(ll^\prime \ ^{\partial_2 l^{\prime{-1}}}(l^{-1}l^{\prime{-1}}),\partial_2 l^{\prime{-1}}\partial_2 l^{\prime},1)\\
  =&(1,1,1)
\end{align*}
  Thus, $(L\rtimes M\rtimes N , s_1,t_1 , s_2,t_2)$ is a cat$^2$-group.

\subsection{The Construction of Cat$^2$-Group Algebra}
 The notion of cat$^1$-algebra is well-known, at least as on anologue of a cat$^1$-group in the category of algebras. The equivalence between cat$^1$-algebras and crossed modules of algebras appears in \cite{Nizar}. More general expositions of cat$^1$-objects were introduced by Ellis \cite{Ellis} and Porter \cite{Porter}. Recall from that a  cat$^1$-algebra $\mathcal{A}$ consists of $K$-algebras $\mathcal{A}_0$,\ $\mathcal{A}_1$ and $K$-algebra morphisms $\sigma,\tau:\mathcal{A}_1\rightarrow\mathcal{A}_0,\ i:\mathcal{A}_0\rightarrow \mathcal{A}_1$ (called structural morphisms) satisfying $1.$ $\sigma_i=\tau_i=id_{\mathcal{A}_0}$,\ $2.$ $\ke \sigma \cdot \ke \tau=0$,\  $\ke \tau \cdot \ke \sigma=0$. Condition 2 is called the kernel condition.

In this section using the method given by Barker in \cite{Barker}, we shall give a construction Gray 3-(group)-algebra groupoid with a single object $*$ or a cat$^2$-group algebra from  $\mathfrak{C}^2$.

 Recall from \cite{Arvasi} that a cat$^2$-algebra $\mathcal{A}^2$ is a 5-tuple $(\mathcal{A},\sigma_1,\tau_1,\sigma_2,\tau_2)$ where $(\mathcal{A},\sigma_i,\tau_i)$ \ $i=1,2$ are  cat$^1$-algebras and
\begin{enumerate}
\item $ \sigma_i\tau_j=\tau_j\sigma_i,\ \tau_j\tau_i=\tau_i\tau_j,\ \sigma_i\tau_j=\tau_j\sigma_i \ \text{for} \ i,j=1,2,\ i\neq j  \ \text{and} $
\item $ \ke \tau_i\cdot\ke \sigma_i=0  \text{ and }\ke \sigma_i\cdot\ke \tau_i=0 \ \text{for} \ i=1,2.$
\end{enumerate}

The cat$^2$-group $\mathfrak{C}^2$ can be regarded as:
$$
\xymatrix{\mathfrak{C}^2:=L\rtimes M \rtimes N \ar@<1ex>[r]^-{s_3,t_3} \ar@<0ex>[r]&1\rtimes M \rtimes N \ar@<1ex>[r]^-{s_2,t_2} \ar@<0ex>[r]\ar@<1ex>[l]^-{e_3}& 1\rtimes 1 \rtimes N\ar@<1ex>[l]^-{e_2}\ar@<0.5ex>[r]\ar@<-0.5ex>[r]& \{*\}}
$$
where\\
$s_3(l,m,n)=(1,m,n)$,\ \ $t_3(l,m,n)=(1,\partial_2 lm,n)$ and
$s_2(1,m,n)=(1,1,n)$,\ \ $t_2(1,m,n)=(1,1,\partial_1 mn)$.\\
In this structure;

1.$$\xymatrix{1\rtimes M \rtimes N \ar@<1ex>[r]^-{s_2,t_2} \ar@<0ex>[r]& 1\rtimes 1\rtimes N}$$
 is a pre-cat$^1$-group with the group multiplication in the big group $1\rtimes M\rtimes N$ given by
  $$(1,m,n)(1,m^\prime,n^\prime)=(1,m^nm^\prime,nn^\prime)$$
and

2. $$\xymatrix{L\rtimes M \rtimes N \ar@<1ex>[r]^-{s_3,t_3} \ar@<0ex>[r]&1\rtimes M \rtimes N }$$
 is a cat$^1$-group with the group multiplication given by
$$(l,m,n)(l^\prime,m^\prime,n^\prime)=(l^m(^{n}l^\prime),m^nm^\prime,nn^\prime).$$
If we apply the group algebra functor $K(.)$ to this structure for each level, we obtain  the following structure:
$$
\xymatrix{K(\mathfrak{C}^2):=K(L\rtimes M \rtimes N) \ar@<1ex>[r]^-{\sigma_3, \tau_3} \ar@<0ex>[r]&K(1\rtimes M \rtimes N) \ar@<1ex>[r]^-{\sigma_2, \tau_2} \ar@<0ex>[r]& K(1\rtimes 1 \rtimes N)\ar@<0.5ex>[r]\ar@<-0.5ex>[r]& \{*\}.}
$$
In this structure, $K(1\rtimes 1\rtimes N)$ has basis $\{\mathbf{e}_{1,1,n}:n\in N\}$, $K(1\rtimes M\rtimes N)$ has basis $\{{\mathbf{e}_{1,m,n}:m\in M,n\in N}\}$ and $K(L\rtimes M\rtimes N)$ has basis $\{\mathbf{e}_{l,m,n}:l\in L,m\in M,n\in N\}$.

 The structural maps are between them;
 $$\sigma_3(\mathbf{e}_{l,m,n})=\mathbf{e}_{1,m,n}\ ,\ \ \tau_3(\mathbf{e}_{l,m,n})=\mathbf{e}_{1,\partial_2lm,n}$$
 $$\sigma_2(\mathbf{e}_{1,m,n})=\mathbf{e}_{1,1,n}\ ,\ \ \tau_2(\mathbf{e}_{1,m,n})=\mathbf{e}_{1,1,\partial_1mn}$$
 $$\sigma_1(\mathbf{e}_{l,m,n})=\mathbf{e}_{1,1,n}\ ,\ \ \tau_1(\mathbf{e}_{l,m,n})=\mathbf{e}_{1,1,\partial_1mn}.$$
 We can picture a 3-cell in this structure by
$$
   \begin{tikzcd}[row sep=0.3cm,column sep=scriptsize]
                & \ar[dd, Rightarrow, "{\scriptscriptstyle \mathbf{e}_{1,m,n}}"{swap,name=f,description}]& &\ar[dd, Rightarrow, "{\scriptscriptstyle \mathbf{e}_{1,\partial_2lm,n}}"'{swap,name=g,description}]& \\
               {*} \ar[rrrr,bend left=40,"{\scriptscriptstyle \mathbf{e}_{1,1,n}}"] \ar[rrrr,bend right=40,"{\scriptscriptstyle \mathbf{e}_{1,1,\partial_1mn}}"']& \tarrow["\scriptscriptstyle {\mathbf{e}_{l,m,n}}" ,from=f,to=g, shorten >= -1pt,shorten <= 1pt ]{rrr}&  & & {*} \\
                \ & \ & \ & \ & \
    \end{tikzcd}
$$
Since $K(\cdot)$ is a functor, the first condition of cat$^2$-group algebra is induced from the equivalent condition on $\mathfrak C^2$.

We have to show that the kernel condition, $\ke\sigma_3\cdot \ke\tau_3=0$.
To show this equality we need to find bases for $\ke \sigma_3$ and $\ke\tau_3$.  For any $\alpha=\mathbf{e}_{l,m,n}\in K(L\rtimes M\rtimes N)$; we obtain
\begin{align*}
\alpha-i_3\sigma_3(\alpha)=&\mathbf{e}_{l,m,n}-i_3\sigma_3(\mathbf{e}_{l,m,n})\\
=&\mathbf{e}_{l,m,n}-\mathbf{e}_{1,m,n}\\
=&\mathbf{v}^{\scriptscriptstyle22}_{l,m,n}\in \ke\sigma_3.
\end{align*}
Similarly $\alpha-i_3\tau_3(\alpha)={\mathbf{e}_{l,m,n}}-{\mathbf{e}_{1,\partial_2 lm,n}}=\mathbf{w}^{\scriptscriptstyle22}_{l,m,n} \in \ke \tau_3. $ We can give the following result.
\begin{lem}
The set $\{\mathbf{v}^{\scriptscriptstyle22}_{l,m,n}:l\neq1\}$ is as basis for $ \ke \sigma_3$ in $K(L\rtimes M\rtimes N).$
\end{lem}
\begin{pf}
 We can say easily that $\mathbf{v}^{\scriptscriptstyle22}_{l,m,n}\in \ke \sigma_3,$ since
$$\sigma_3(\mathbf{v}^{\scriptscriptstyle22}_{l,m,n})=\sigma_3(\mathbf{e}_{l,m,n}-\mathbf{e}_{1,m,n})=\mathbf{e}_{1,m,n}-\mathbf{e}_{1,m,n}=0$$
We need to show that the elements of the set $\{\mathbf{v}^{\scriptscriptstyle22}_{l,m,n}:l\neq1\}$ span $\ke \sigma_3$ and are linearly independent.
Let $\mathbf{v}\in \ke \sigma_3.$ We can write it as\\
\begin{equation*} \mathbf{v}=\sum_{n \in N}\left(\sum_{m \in M}\left(\sum_{l \in L} r_{l,m,n}\mathbf{e}_{l,m,n}\right)\right) \end{equation*} with $\sigma_3(\mathbf{v})=0$. But
\begin{equation*} \sigma_3 \mathbf{(\mathbf{v})}=\sum_{n \in N}\sum_{m \in M}\left(\sum_{l \in L} r_{l,m,n}\right)\mathbf{e}_{1,m,n}
\end{equation*}
and since $\mathbf{e}_{1,m,n}$ is a basis for $K(L\rtimes M\rtimes N)$, this is zero if $\sum_{l \in L}r_{l,m,n}=0$ for each $m \in M$ and $n \in N$.

Now,
\begin{align*}
\mathbf{v}=&\sum_{n}\sum_{m}\sum_{l\neq 1} r_{l,m,n}\left(\mathbf{e}_{l,m,n}-\mathbf{e}_{1,m,n}+\mathbf{e}_{1,m,n}\right)+\sum_{n}\sum_{m}r_{1,m,n}\mathbf{e}_{1,m,n} \\
          =&\sum_{n}\sum_{m}\sum_{l\neq 1} r_{l,m,n}\mathbf{v}^{\scriptscriptstyle22}_{l,m,n}+\sum_{n}\sum_{m}\sum_{l}r_{l,m,n}\mathbf{e}_{1,m,n}
\end{align*}
but
\begin{equation*}
\sum_{n}\sum_{m}\sum_{l\neq 1} r_{l,m,n}\mathbf{e}_{1,m,n}=\sum_{n}\sum_{m}\left(\sum_{l}r_{l,m,n}\right)\mathbf{e}_{1,m,n}=0
\end{equation*}
because $\sum_{l}r_{l,m,n}=0$ for every $(1,m,n) \in 1\rtimes M\rtimes N$.
Hence,
\begin{equation*}
\mathbf{v}=\sum_{n}\sum_{m}\sum_{l\neq 1} r_{l,m,n}\mathbf{v}^{\scriptscriptstyle22}_{l,m,n}
\end{equation*}
for every $\mathbf{v}\in \ke \sigma_3$, so the $\mathbf{v}^{\scriptscriptstyle22}_{l,m,n}$ do indeed span $\ke \sigma_3$.
Now suppose that
\begin{equation*}
\mathbf{v}=\sum_{n}\sum_{m}\sum_{l\neq 1} r_{l,m,n}\mathbf{v}^{\scriptscriptstyle22}_{l,m,n}=0
\end{equation*}
Then,
\begin{align*}
\sum_{n}\sum_{m}\sum_{l\neq 1} r_{l,m,n}\left(\mathbf{e}_{l,m,n}-\mathbf{e}_{1,m,n}\right)=0\Leftrightarrow&\sum_{n}\sum_{m}\sum_{l\neq 1}r_{l,m,n}\mathbf{e}_{l,m,n}-\sum_{n}\sum_{m}\sum_{l\neq 1} r_{l,m,n}\mathbf{e}_{1,m,n}=0\\
\Leftrightarrow&\sum_{n}\sum_{m}\sum_{l}r'_{l,m,n}\mathbf{e}_{l,m,n}=0
\end{align*}
where $r'_{l,m,n}=r_{l,m,n}$ when $l\neq 1$ and $r'_{1,m,n}=-\sum_{l\neq1}r_{l,m,n}$. This is a linear combination of basis vectors $\mathbf{e}_{l,m,n}$ in $K(L\rtimes M\rtimes N)$ so $r'_{l,m,n}=0$ for each $(l,m,n) \in L\rtimes M\rtimes N$. In particular, this is true for $l \neq 1$, so every $r_{l,m,n}=0$ and  $\mathbf{v}^{\scriptscriptstyle22}_{l,m,n}$ are linearly independent.
\end{pf}

Consider the elements  $$\mathbf{e}^{*_2}_{l,m,n}={i_3\sigma_3(\mathbf{e}_{l,m,n})}-{\mathbf{e}_{l,m,n}+i_3\tau_3\mathbf{e}_{l,m,n}}=\mathbf{e}_{1,m,n}-\mathbf{e}_{l,m,n}+\mathbf{e}_{(1,\partial_2 lm,n)}\in K(L\rtimes M \rtimes N)$$
while
$$\mathbf{e}^{*_2}_{1,m,n}=i_3\sigma_3(\mathbf{e}_{1,m,n})-\mathbf{e}_{1,m,n}+i_3\tau_3(\mathbf{e}_{1,m,n})=\mathbf{e}_{1,m,n}-\mathbf{e}_{1,m,n}+\mathbf{e}_{1,m,n}
=\mathbf{e}_{1,m,n}\in K(1\rtimes M \rtimes N)$$ Therefore,
\begin{align*}
(\mathbf{v}^{\scriptscriptstyle22}_{l,m,n})^{*_2}=&(\mathbf{e}_{l,m,n}-\mathbf{e}_{1,m,n})^{*_2}\\
                               =&\mathbf{e}_{1,m,n}-\mathbf{e}_{l,m,n}+\mathbf{e}_{1,\partial_2 lm,n}-\mathbf{e}_{1,m,n}\\
                               =&-(\mathbf{e}_{l,m,n}-\mathbf{e}_{1,\partial_2 lm,n})\\
                               =&-\mathbf{w}^{\scriptscriptstyle22}_{l,m,n}\in \ke \tau_3.
\end{align*}
Hence $\{{-\mathbf{w}^{\scriptscriptstyle22}_{l,m,n}:l\neq 1}\}$ is a basis for $\ke \tau_3$
and so $\{{\mathbf{w}^{\scriptscriptstyle22}_{l,m,n}:l\neq 1}\}$ is also basis for $\ke \tau_3.$
Therefore we have proven the following lemma.
\begin{lem}
The set $\{{\mathbf{w}^{\scriptscriptstyle22}_{l,m,n}:l\neq 1}\}$ is a basis for $\ke \tau_3$.\\
\end{lem}
\begin{lem}
$\xymatrix{K(1\rtimes M \rtimes N) \ar@<1ex>[r]^-{\sigma_2,\tau_2} \ar@<0ex>[r]& K(1\rtimes 1\rtimes N)\ar@<1ex>[l]^-{i_2}}$ is a pre-cat$^1$-group algebra.
\end{lem}
\begin{pf}
According to the multiplication in $K(1\rtimes M\rtimes N)$  given by
$$\mathbf{e}_{1,m,n}\mathbf{e}_{1,m^\prime, n^\prime}=\mathbf{e}_{1,m^nm^\prime,nn^\prime},$$
we can easily see that $\sigma_2$, $\tau_2$ and $i_2$ are homomorphisms of algebras and $\sigma_2i_2=\tau_2i_2=id.$
\end{pf}

For any elements $\mathbf{e}_{l,m,n}$ and $\mathbf{e}_{l^\prime,m^\prime,n^\prime}$ in $ K(L\rtimes M\rtimes N),$ we obtain
 $$\mathbf{e}_{l,m,n}-\mathbf{e}_{1,m,n}=\mathbf{v}^{\scriptscriptstyle22}_{l,m,n}\in \ke \sigma_3$$
 and
 $$\mathbf{w}^{\scriptscriptstyle22}_{l^\prime,m^\prime,n^\prime}=\mathbf{e}_{l^\prime,m^\prime,n^\prime}-\mathbf{e}_{1,\partial_2 l^\prime m^\prime,n^\prime} \in \ke \tau_3.$$

 The multiplication of these elements in $ K(L\rtimes M\rtimes N)$  is
 \begin{align*}
 \mathbf{v}^{\scriptscriptstyle22}_{l,m,n}\cdot\mathbf{w}^{\scriptscriptstyle22}_{l^\prime,m^\prime,n^\prime}=&(\mathbf{e}_{l,m,n}-\mathbf{e}_{1,m,n})\cdot(\mathbf{e}_{l^\prime,m^\prime,n^\prime}-\mathbf{e}_{1,\partial_2 l^\prime m^\prime,n^\prime})\\
 =&\mathbf{e}_{l^m (^nl^\prime),m^nm^\prime,nn^\prime}-\mathbf{e}_{l,m^n(\partial_2l^\prime m^\prime),nn^\prime}-\mathbf{e}_{^m(^nl^\prime),m^nm^\prime,nn^\prime}+\mathbf{e}_{1,m^n(\partial_2l^\prime m^\prime),nn^\prime}\neq 0
 \end{align*}
Thus, in the structure,
$$
\xymatrix{K(L\rtimes M \rtimes N)\ar@<1ex>[r]^-{\sigma_3,\tau_3} \ar@<0ex>[r]&K(1\rtimes M \rtimes N) \ar@<1ex>[l]^-{i_3}}
$$
we see that $\ke\sigma_3\cdot\ke \tau_3\neq 0.$  In order to construct a cat$^1$-algebra from $K(\mathfrak C^2)$, it is necessary to impose some relations in $K(L\rtimes M\rtimes N)$, so that the kernel condition is satisfied. Suitable expressions are formed $K(L\rtimes M\rtimes N)$ can be factored by the ideal generated by these relations.
Consider the expressions of the form
 \begin{align*}
{\mathbf{u}_{\scriptscriptstyle1}}=&\mathbf{e}_{l^\prime l,m,n}-\mathbf{e}_{l,m,n}-\mathbf{e}_{l^\prime,\partial_2lm,n}+\mathbf{e}_{1,\partial_2lm,n}\\
{\mathbf{v}_{\scriptscriptstyle1}}=&\mathbf{e}_{l^{\prime m^\prime} l,m^\prime m,n}-\mathbf{e}_{l,m,n}-\mathbf{e}_{l^\prime,m^\prime,\partial_1mn}+\mathbf{e}_{1,1,\partial_1mn}\\
{\mathbf{v}_2}=&\mathbf{e}_{1,m^\prime m,n}-\mathbf{e}_{1,m,n}-\mathbf{e}_{1,m^\prime,\partial_1mn}+\mathbf{e}_{1,1,\partial_1mn}
\end{align*}
in $K(L\rtimes M\rtimes N)$.  These elements can be represented diagramatically as
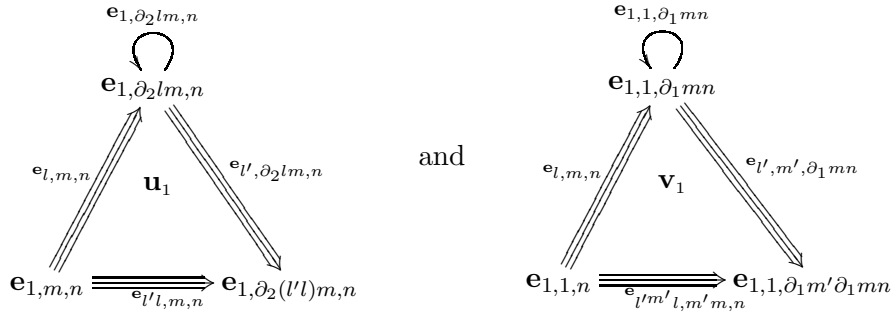
$$
\begin{array}{c}
\xymatrix@C-=0.0cm{
& \mathbf{e}_{1,\partial_2lm,n}\ar@3{->}[ddr]^{\scriptscriptstyle \mathbf{e}_{l',\partial_2 lm,n}}\ar@(ur,ul)_{\mathbf{e}_{1,\partial_2lm,n}} & \\
 &  \ {\mathbf{u}_{\scriptscriptstyle1}} &  \\
\mathbf{e}_{1,m,n} \ar@3{->}[rr]_{\scriptscriptstyle \mathbf{e}_{l'l,m,n}}\ar@3{->}[uur]^{\scriptscriptstyle \mathbf{e}_{l,m,n}}   &  & \mathbf{e}_{1,\partial_2(l'l)m,n} }
\end{array}
\text{ \ \ \ and  \ \ \ }
\begin{array}{c}
\xymatrix@C-=0.0cm{
& \mathbf{e}_{1,1,\partial_1 mn}\ar@3{->}[ddr]^{\scriptscriptstyle \mathbf{e}_{l',m',\partial_1 mn}}\ar@(ur,ul)_{\mathbf{e}_{1,1,\partial_1 mn}} & \\
 &  \ \ {\mathbf{v}_{\scriptscriptstyle1}} &  \\
\mathbf{e}_{1,1,n} \ar@3{->}[rr]_{\scriptscriptstyle \mathbf{e}_{l^{\prime m^\prime}l,m'm,n}}\ar@3{->}[uur]^{\scriptscriptstyle \mathbf{e}_{l,m,n}}   & & \mathbf{e}_{1,1,\partial_1m'\partial_1 mn} }
\end{array}
$$
In these pictures, each arrow represents a 3-cell in $K(L\rtimes M\rtimes N)$. For the generating element
 $$ {\mathbf{u}_{\scriptscriptstyle1}}= \mathbf{e}_{l^\prime l,m,n}-\mathbf{e}_{l,m,n}-\mathbf{e}_{l^\prime,\partial_2lm,n}+\mathbf{e}_{1,\partial_2lm,n},$$
if we take $(l,m,n)=(1,m,n)$, we obtain
 $$\mathbf{e}_{l^\prime,m,n}-\mathbf{e}_{1,m,n}-\mathbf{e}_{l^\prime,m,n}+\mathbf{e}_{1,m,n}=0,$$
 and if we take $(l^\prime,m^\prime,n^\prime)=(1,m^\prime,n^\prime)$ we obtain
 $$\mathbf{e}_{l,m,n}-\mathbf{e}_{l,m,n}-\mathbf{e}_{1,\partial_2lmn}+\mathbf{e}_{1,\partial_2lmn}=0.$$
 Therefore, we can rewrite this expression as follows:
 \begin{equation}\label{1}
  \mathbf{e}_{l'l,m,n}=\mathbf{e}_{l,m,n}+\mathbf{e}_{l',\partial_2lm,n}-\mathbf{e}_{1,\partial_2lm,n}.
 \end{equation}
 Similarly, for the expressions
  \begin{align*}
 {\mathbf{v}_{\scriptscriptstyle1}}=&\mathbf{e}_{l^{\prime m^\prime}l,m^\prime m,n}-\mathbf{e}_{l,m,n}-\mathbf{e}_{l^\prime,m^\prime,\partial_1mn}+\mathbf{e}_{1,1,\partial_1mn}\\
 {\mathbf{v}_{\scriptscriptstyle2}}=&\mathbf{e}_{1,m^\prime m,n}-\mathbf{e}_{1,m,n}-\mathbf{e}_{1,m^\prime,\partial_1mn}+\mathbf{e}_{1,1,\partial_1mn}
  \end{align*}
 if we take $(l,m,n)=(1,1,n)$, we obtain
 $$\mathbf{e}_{l^\prime,m^\prime,n}-\mathbf{e}_{1,1,n}-\mathbf{e}_{l^\prime,m^\prime,n}+\mathbf{e}_{1,1,n}=0$$
 and if we take $(l^\prime,m^\prime,n^\prime)=(1,1,n^\prime)$, we obtain
 $$\mathbf{e}_{l,m,n}-\mathbf{e}_{l,m,n}-\mathbf{e}_{1,1,\partial_1mn}+\mathbf{e}_{1,1,\partial_1mn}=0.$$
 We can rewrite these expression as follows:
\begin{equation}\label{2}
   \mathbf{e}_{l^{\prime m^\prime} l,m^\prime m,n}=\mathbf{e}_{l,m,n}+\mathbf{e}_{l^\prime,m^\prime,\partial_1mn}-\mathbf{e}_{1,1,\partial_1mn}
 \end{equation}
 and
 \begin{equation}\label{3}
 \mathbf{e}_{1,m^\prime m,n}=\mathbf{e}_{1,m,n}+\mathbf{e}_{1,m^\prime,\partial_1mn}-\mathbf{e}_{1,1,\partial_1mn}.
 \end{equation}
 Let $J_2$ be ideal of $K(L\rtimes M\rtimes N)$ generated by the elements ${\mathbf{u}_{\scriptscriptstyle1}},{\mathbf{v}_{\scriptscriptstyle1}},{\mathbf{v}_{\scriptscriptstyle2}}$ and define
 $$\overline{K(L\rtimes M\rtimes N)}=K(L\rtimes M\rtimes N)/ J_2.$$

We can also consider the expression
$${\mathbf{v}_{\scriptscriptstyle2}}=\mathbf{e}_{1,m^\prime m,n}-\mathbf{e}_{1,m,n}-\mathbf{e}_{1,m^\prime,\partial_1mn}+\mathbf{e}_{1,1,\partial_1mn}$$
 in  $K(1\rtimes M\rtimes N)$.  Let $J_1$ be the ideal of $K(1\rtimes M\rtimes N)$ generated by the elements ${\mathbf{v}_{\scriptscriptstyle2}}$. Since
 \begin{align*}
 \sigma_3({\mathbf{v}_{\scriptscriptstyle1}})=&\mathbf{e}_{1,m^\prime m,n}-\mathbf{e}_{1,m,n}-\mathbf{e}_{1,m^\prime,\partial_1mn}+\mathbf{e}_{1,1,\partial_1mn}\in J_1\\
 \tau_3({\mathbf{v}_{\scriptscriptstyle1}})=&\mathbf{e}_{1,\partial_2l^\prime m^\prime \partial_2lm,n}-\mathbf{e}_{1,\partial_2lm,n}-\mathbf{e}_{1,\partial_2l^\prime m^\prime,\partial_1mn}+\mathbf{e}_{1,1,\partial_1mn}\in J_1\\
 \sigma_3({\mathbf{u}_{\scriptscriptstyle1}})=&\mathbf{e}_{1,m,n}-\mathbf{e}_{1,m,n}-\mathbf{e}_{1,\partial_2lm,n}+\mathbf{e}_{1,\partial_2lm,n}=0\\
 \tau_3({\mathbf{u}_{\scriptscriptstyle1}})=&\mathbf{e}_{1,\partial_2l^\prime \partial_2lm,n}-\mathbf{e}_{1,\partial_2lm,n}-\mathbf{e}_{1,\partial_2l^\prime \partial_2lm,n}+\mathbf{e}_{1,\partial_2lm,n}=0
\end{align*}
 We obtain $\sigma_3(J_2)\subseteq J_1$ and $\tau_3(J_2)\subseteq J_1.$ Therefore, we can give the following commutative diagram;
$$
\xymatrix@C-=0.4cm{K(L\rtimes M\rtimes N)/J_2 \ar@<0.5ex>[rr]^{\overline{\sigma}_{3}}\ar@<-0.5ex>[rr]_{\overline{\tau}_{3}}&&K(1\rtimes M\rtimes N)/J_{1}\\ \\ K(L\rtimes M\rtimes N)\ar@<0.5ex>[rr]^{\sigma_{3}}\ar@<-0.5ex>[rr]_{\tau_{3}}\ar@{->>}[uu]^{q_{2}}&&K(1\rtimes M \rtimes N)\ar@{->>}[uu]_{q_{1}}}
$$where $q_1$ and $q_2$ are quotient maps. Since $\sigma_3(J_2)\subseteq J_1$, the map $\overline {\sigma}_3$ given  by
$$\overline {\sigma}_3(\overline {\mathbf{e}}_{l,m,n})=\overline {\sigma}_3(\mathbf{e}_{l,m,n}+J_2)=\sigma_3(\mathbf{e}_{l,m,n})+J_1=\mathbf{e}_{1,m,n}+J_1=\overline {\mathbf{e}}_{1,m,n}$$
is a well-defined homomorphism. Since $\tau_3(J_2)\subseteq J_1$, the map $\overline{\tau}_3$ given by
$$\overline{\tau}_3(\overline {\mathbf{e}}_{l,m,n})=\overline{\tau}_3(\mathbf{e}_{l,m,n}+J_2)=\tau_3(\mathbf{e}_{l,m,n})+J_1=\mathbf{e}_{1,\partial_2lm,n}+J_1=\overline{\mathbf{e}}_{1,\partial_2lm,n}$$
is a well-defined homomorphism of algebras.

Similarly, $$\overline{i}_3(\overline{\mathbf{e}}_{1,m,n})=\overline{{i}_3(\mathbf{e}_{1,m,n})}=i_3(\mathbf{e}_{1,m,n})+J_2=\overline{\mathbf{e}}_{1,m,n}.$$

Let $\overline{K(L\rtimes M \rtimes N)}=K(L\rtimes M \rtimes N)/J_2$ and $\overline{K(1\rtimes M \rtimes N)}=K(1\rtimes M \rtimes N)/J_1.$
We obtain the following structure;
$$
\xymatrix{\overline{K(\mathfrak{C}^2)}:=\overline{K(L\rtimes M \rtimes N)}\ar@<1ex>[r]^-{\overline{\sigma}_3,\overline{\tau}_3} \ar@<0ex>[r]&\overline{K(1\rtimes M \rtimes N)} \ar@<1ex>[r]^-{\overline{\sigma}_2,\overline{\tau}_2} \ar@<0ex>[r]\ar@<1ex>[l]^-{\overline{i}_3}& K(1\rtimes 1\rtimes N)\ar@<1ex>[l]^-{\overline{i}_2}\ar@<0.5ex>[r]\ar@<-0.5ex>[r]& \{*\}}
$$
where
\begin{align*}
\overline{\sigma}_2(\overline{\mathbf{e}}_{1,m,n})=&\mathbf{e}_{1,1,n}, \\
\overline{\tau}_2(\overline{\mathbf{e}}_{1,m,n})=&\mathbf{e}_{1,1,\partial_1mn},\\
\overline{\sigma}_3(\overline{\mathbf{e}}_{l,m,n})=&\overline{\mathbf{e}}_{1,m,n},\\
\overline{\tau}_3(\overline{\mathbf{e}}_{l,m,n})=&\overline{\mathbf{e}}_{1,\partial_2lm,n}.
\end{align*}
Since;
$$\sigma_2({\mathbf{v}_{\scriptscriptstyle2}})=\mathbf{e}_{1,1,n}-\mathbf{e}_{1,1,n}-\mathbf{e}_{1,1,\partial_1mn}+\mathbf{e}_{1,1,\partial_1mn}=0$$
and
$$\tau_2({\mathbf{v}_{\scriptscriptstyle2}})=\mathbf{e}_{1,1,\partial_1 m^\prime \partial_1mn}-\mathbf{e}_{1,1,\partial_1mn}-\mathbf{e}_{1,1,\partial_1 m^\prime \partial_1mn}+\mathbf{e}_{1,1,\partial_1mn}=0,$$
we have $\sigma_2(J_1)=0$ and $\tau_2(J_1)=0.$ Therefore, the maps $\overline{\sigma}_2$ and $\overline{\tau}_2$ are well-defined homomorphisms.

In $\overline{K(\mathfrak{C}^2)}$ , it can be easily seen that
$$\xymatrix{\overline{K(1\rtimes M \rtimes N)} \ar@<1ex>[r]^-{\overline{\sigma}_2,\overline{\tau}_2} \ar@<0ex>[r]& K(1\rtimes 1\rtimes N)\ar@<1ex>[l]^-{\overline{i}_2}}$$
 is a pre-cat$^1$-group algebra. Now we will show that
 $$\xymatrix{\overline{K(L\rtimes M \rtimes N)}\ar@<1ex>[r]^-{\overline{\sigma}_3,\overline{\tau}_3} \ar@<0ex>[r]&\overline{K(1\rtimes M \rtimes N)}\ar@<1ex>[l]^-{\overline{i}_3}}$$
 is a cat$^1$-group algebra. To show this, the kernel condition for $\overline{\sigma}_3,\overline{\tau}_3$ must be satisfied.
Let
\begin{align*}
\overline{\mathbf{v}^{\scriptscriptstyle22}}_{l,m,n}=&\overline{\mathbf{v}^{\scriptscriptstyle22}}_{l,m,n}+J_2\in \ke \overline{\sigma}_3\\
\overline{\mathbf{w}^{\scriptscriptstyle22}}_{l',m',n'}=&\overline{\mathbf{w}^{\scriptscriptstyle22}}_{l',m',n'}+J_2\in \ke \overline{\tau}_3.
\end{align*}
Then, we obtain
$$\overline{\mathbf{v}^{\scriptscriptstyle22}}_{l,m,n}\cdot\overline{\mathbf{w}^{\scriptscriptstyle22}}_{l',m',n'}=\mathbf{e}_{l^m(^nl'),m^n m',nn'}-\mathbf{e}_{^m(^nl'),m^n m^\prime,nn^\prime}-\mathbf{e}_{l,m^n(\partial_2l^\prime m^\prime),nn^\prime}+\mathbf{e}_{1,m^n(\partial_2l^\prime m^\prime),nn^\prime}+J_2$$
and by considering the relations (\ref{1}),\ (\ref{2}) and (\ref{3}), we obtain
$$\overline{\mathbf{v}^{\scriptscriptstyle22}}_{l,m,n}\cdot\overline{\mathbf{w}^{\scriptscriptstyle22}}_{l^\prime,m^\prime,n^\prime}=\overline{0}=0+J_2.$$
Thus, we obtain that $\ke \overline{\sigma}_3. \ke \overline{\tau}_3=\overline{0}.$

Relation (\ref{1}) can be rewritten as
\begin{multline*}
(\mathbf{e}_{l^\prime l,m,n}-\mathbf{e}_{1,m,n})-(\mathbf{e}_{l,m,n}-\mathbf{e}_{1,m,n})-(\mathbf{e}_{l^\prime,\partial_2 lm,n}-\mathbf{e}_{1,\partial_2 lm,n})\\
\begin{aligned}
&=\overline{\mathbf{v}^{\scriptscriptstyle22}}_{l^\prime l,m,n}-\overline{\mathbf{v}^{\scriptscriptstyle22}}_{l,m,n}-\overline{\mathbf{v}^{\scriptscriptstyle22}}_{l^\prime,\partial_2 lm,n}.
\end{aligned}
\end{multline*}
Since this is in the ideal $J_2$, it will be killed off by factorisation. Thus in $K(L\rtimes M\rtimes N)/ J_2$, we get
$$\overline{\mathbf{v}^{\scriptscriptstyle22}}_{l^\prime l,m,n}=\overline{\mathbf{v}^{\scriptscriptstyle22}}_{l,m,n}+\overline{\mathbf{v}^{\scriptscriptstyle22}}_{l^\prime,\partial_2 lm,n}.$$

Relation (\ref{2}) can be rewritten as
\begin{multline*}
(\mathbf{e}_{l^\prime{^{m^\prime}}l,m^\prime m,n}-\mathbf{e}_{1,m^\prime m,n}+\mathbf{e}_{1,m^\prime m,n})-(\mathbf{e}_{l,m,n}-\mathbf{e}_{1,m,n}+\mathbf{e}_{1,m,n})\\
\begin{aligned}
&\hspace{3.7cm} -(\mathbf{e}_{l^\prime,m^\prime,\partial_1mn}-\mathbf{e}_{1,m^\prime,\partial_1mn}+\mathbf{e}_{1,m^\prime,\partial_1mn})+\mathbf{e}_{1,1,\partial_1mn}\\
=&{\mathbf{v}^{\scriptscriptstyle22}}_{l^\prime{^{m^\prime}}l,m^\prime m,n}+\mathbf{e}_{1,m^\prime m,n}-{\mathbf{v}^{\scriptscriptstyle22}}_{l,m,n}-\mathbf{e}_{1,m,n}-{\mathbf{v}^{\scriptscriptstyle22}}_{l^\prime,m^\prime,\partial_1mn}-\mathbf{e}_{1,m^\prime,\partial_1mn}+\mathbf{e}_{1,1,\partial_1mn}\\
=&{\mathbf{v}^{\scriptscriptstyle22}}_{l^\prime{^{m^\prime}}l,m^\prime m,n}-{\mathbf{v}^{\scriptscriptstyle22}}_{l,m,n}-{\mathbf{v}^{\scriptscriptstyle22}}_{l^\prime,m^\prime,\partial_1 mn}+\underbrace{\mathbf{e}_{1,m^\prime m,n}-\mathbf{e}_{1,m,n}-\mathbf{e}_{1,m^\prime,\partial_1 mn}+\mathbf{e}_{1,1,\partial_1mn} }_{\mathbf{v'}},
\end{aligned}
\end{multline*}
where ${\mathbf{v}_{\scriptscriptstyle2}}=\mathbf{e}_{1,m^\prime m,n}-\mathbf{e}_{1,m,n}-\mathbf{e}_{1,m^\prime,\partial_1 mn}+\mathbf{e}_{1,1,\partial_1 mn} \in J_2$, we obtain
${\mathbf{v}_{\scriptscriptstyle2}}+J_2=0+J_2$ and thus
$$\overline{\mathbf{v}^{\scriptscriptstyle22}}_{l'{^{m^\prime}}l,m^\prime m,n}=\overline{\mathbf{v}^{\scriptscriptstyle22}}_{l,m,n}+\overline{\mathbf{v}^{\scriptscriptstyle22}}_{l^\prime,m^\prime,\partial_1 mn}.$$
There are redundancies among the $\overline{\mathbf{v}^{\scriptscriptstyle22}}_{l,m,n}$ so these do not form a basis for $\ke \overline{\sigma}_3$. Similar relations hold for the $\overline{\mathbf{w}^{\scriptscriptstyle22}}_{l',m',n'}$.
\subsubsection{The Chain Complex $\overline{\delta}$ from $\overline{K(\mathfrak{C}^2)}$} \label{complex}

For any  $\overline{\mathbf{e}}_{1,m,n}\in K(1\rtimes M\rtimes N)/ J_1$, we obtain
$$\overline{\mathbf{v}^{\scriptscriptstyle11}}_{1,m,n}=\overline{\mathbf{e}}_{1,m,n}-\overline{\mathbf{e}}_{1,1,n}\in \ke \overline{\sigma}_2=\mathrm{K}_2$$ and for any $\overline{\mathbf{e}}_{l,m,n}\in K(L\rtimes M\rtimes N) /J_2$,  we obtain
$$\overline{\mathbf{v}^{\scriptscriptstyle22}}_{l,m,n}=\overline{\mathbf{e}}_{l,m,n}-\overline{\mathbf{e}}_{1,m,n}\in \ke \overline{\sigma}_3=\mathrm{K}_3.$$
For any $\overline{\mathbf{v}^{\scriptscriptstyle22}}_{l,m,n}\in \mathrm{K}_3 $; we obtain
\begin{align*}
\overline{\tau}_2 \overline{\tau}_3(\overline{\mathbf{v}^{\scriptscriptstyle22}}_{l,m,n})=&\overline{\tau}_2 \left(\overline{\tau}_3(\overline{\mathbf{e}}_{l,m,n}-\overline{\mathbf{e}}_{1,m,n})\right) \\
=&\overline{\tau}_2(\overline{\mathbf{e}}_{1,\partial_2 lm,n}-\overline{\mathbf{e}}_{1,m,n}) \\
=&\overline{\tau}_2(\left(\overline{\mathbf{e}}_{1,\partial_2 lm,n}-\overline{\mathbf{e}}_{1,1,n})-(\overline{\mathbf{e}}_{1,m,n}-\overline{\mathbf{e}}_{1,1,n})\right)\\
=&\overline{\tau}_2(\overline{\mathbf{v}^{\scriptscriptstyle11}}_{1,\partial_2 lm,n}-\overline{\mathbf{v}^{\scriptscriptstyle11}}_{1,m,n})\\
=&\mathbf{e}_{1,1,\partial_1\partial_2(l)\partial_1mn}-\mathbf{e}_{1,1,\partial_1mn}\\
=&\mathbf{e}_{1,1,\partial_1mn}-\mathbf{e}_{1,1,\partial_1mn}  \hspace{4cm}  (\because \partial_1\partial_2l=1)\\
=&0.
\end{align*}
Thus, we have the following chain complex of length-2:
$$
\begin{tikzcd}
             \overline{\delta}:=\mathrm{K}_3\ar[r,"\overline{\tau}_3"]&\mathrm{K}_2\ar[r,"\overline{\tau}_2"]& \mathrm{K}_1
\end{tikzcd}
$$
where $\mathrm{K}_1=K(1\rtimes 1\rtimes N)$. Thus, $\overline{\delta}$ can be considered as an object of $\mathbf{Ch}^2_K$ by ignoring the algebra multiplications in each level, so the components of this chain complex are vector spaces and the boundaries are linear transformations.

Now, we show that the construction of $\overline{K(\mathfrak{C}^2)}$ from $\mathfrak{C}^2$ is functorial. For a 2-crossed module
$$\xymatrix{\mathfrak{X}:=(L\ar[r]^-{\partial_2}&M\ar[r]^{\partial_1}&N),}$$
we obtained a Gray 3-(group)-groupoid with a single 0-cell;
$$
\xymatrix{\mathfrak{C}^2:=L\rtimes M \rtimes N \ar@<1ex>[r]^-{s_3,t_3} \ar@<0ex>[r]&1\rtimes M \rtimes N \ar@<1ex>[r]^-{s_2,t_2} \ar@<0ex>[r]\ar@<1ex>[l]^-{e_3}&1\rtimes 1\rtimes N\ar@<1ex>[l]^-{e_2}\ar@<0.5ex>[r]\ar@<-0.5ex>[r]& \{*\}}
$$
and
$$
\xymatrix{ \overline{K(\mathfrak{C}^2)}:=K(L\rtimes M \rtimes N)/ J_2  \ar@<1ex>[r]^-{\overline{\sigma}_3,\overline{\tau}_3} \ar@<0ex>[r]& K(1\rtimes M \rtimes N)/J_1 \ar@<1ex>[r]^-{\overline{\sigma}_2,\overline{\tau}_2} \ar@<0ex>[r]\ar@<1ex>[l]^-{\overline{i}_3}& K(1\rtimes 1\rtimes N)\ar@<1ex>[l]^-{\overline{i}_2}\ar@<0.5ex>[r]\ar@<-0.5ex>[r]& \{*\}}
$$
given above.

Suppose now that $\mathfrak{C}^2_i:=((L_i\rtimes M_i\rtimes N_i),(1\rtimes M_i\rtimes N_i),(1\rtimes 1 \rtimes N_i),s_i,t_i,i_i)$ are Gray 3-(group)-groupoids with a single object $*$, for $i=1,2,3$ and $\phi:\mathfrak{C}^2_1\rightarrow \mathfrak{C}^2_2$ and $\psi:\mathfrak{C}^2_2\rightarrow \mathfrak{C}^2_3$ are (strict) Gray functor between them.
Applying the functor $K(\cdot)$ for each level gives us the following diagram of group algebras;
$$
\xymatrix@C-=0.7cm{K(L_1\rtimes M_1\rtimes N_1)\ar@<1.0ex>[dd]|{\tau_3}\ar@<-1.0ex>[dd]|<<<<<{\scriptstyle \sigma_3}\ar[r]^-{\scriptscriptstyle K(\phi_3)}& K(L_2\rtimes M_2\rtimes N_2)\ar@<1.0ex>[dd]|{\tau_3}\ar@<-1.0ex>[dd]|<<<<<{\sigma_3}\ar[r]^-{\scriptscriptstyle K(\psi_3)}& K(L_3\rtimes M_3\rtimes N_3)\ar@<1.0ex>[dd]|{\scriptstyle \tau_3}\ar@<-1.0ex>[dd]|<<<<<{\scriptstyle \sigma_3}\\
\\
K(1\rtimes M_1\rtimes N_1)\ar@<1.0ex>[dd]|{\scriptstyle \tau_2}\ar@<-1.0ex>[dd]|<<<<<{\scriptstyle \sigma_2}\ar[r]^-{\scriptscriptstyle K(\phi_2)}&K(1\rtimes M_2\rtimes N_2)\ar[r]^-{\scriptscriptstyle K(\psi_2)}\ar@<1.0ex>[dd]|{\scriptstyle \tau_2}\ar@<-1.0ex>[dd]|<<<<<{\scriptstyle \sigma_2}&K(1\rtimes M_3\rtimes N_3)\ar@<1.0ex>[dd]|{\scriptstyle \tau_2}\ar@<-1.0ex>[dd]|<<<<<{\scriptstyle \sigma_2}\\
\\
K(1\rtimes 1\rtimes N_1)\ar[r]_-{\scriptscriptstyle K(\scriptstyle \phi_1)}&K(1\rtimes 1\rtimes N_2)\ar[r]_-{\scriptscriptstyle K(\psi_1)}&K(1\rtimes 1\rtimes N_3)}
$$
where
$$K(\phi_1)(\mathbf{\mathbf{e}}_{1,1,n_1})=\mathbf{e}_{1,1,{\phi_1}(n_1)}, \  \ K(\phi_2)(\mathbf{e}_{1,m_1,n_1})=\mathbf{e}_{{\phi_2}(1,m_1,n_1)}$$ and
$$K(\phi_3)(\mathbf{e}_{l_1,m_1,n_1})=\mathbf{e}_{{\phi_3}(l_1,m_1,n_1)}$$ with similar conditions for $K(\psi_i)$,  $i=1,2,3.$ We form the ideals ${J}^i_2$ of $K(L_i\rtimes M_i\rtimes N_i)$ generated by elements of the forms
\begin{align*}
\mathbf{u}_i=&\mathbf{e}_{{l}^\prime_i l_i,m_i,n_i}-\mathbf{e}_{l_i,m_i,n_i}-\mathbf{e}_{{l}^\prime_i,\partial_2l_im_i,n_i}+\mathbf{e}_{1,\partial_2l_im_i,n_i}\\
\mathbf{v}_i=&\mathbf{e}_{{l'_i}^{{m'_i}}{l_i},{m}^\prime_im_i,n_i}-\mathbf{e}_{l_i,m_i,n_i}-\mathbf{e}_{{l}^\prime_i,{m}^\prime_i,\partial_1m_in_i}+\mathbf{e}_{1,1,\partial_1m_in_i}\\
\mathbf{v'}_i=&\mathbf{e}_{1,{m}^\prime_im_i,n_i}-\mathbf{e}_{1,m_i,n_i}-\mathbf{e}_{1,{m}^\prime_i,\partial_1m_in_i}+\mathbf{e}_{1,1,\partial_1m_in_i}\\
\intertext{and ${J}^i_1$ of $K(1\rtimes M_i\rtimes N_i)$ generated by the elements of the form}
\mathbf{v'}_i=&\mathbf{e}_{1,{m}^\prime_im_i,n_i}-\mathbf{e}_{1,m_i,n_i}-\mathbf{e}_{1,{m}^\prime_i,\partial_1m_in_i}+\mathbf{e}_{1,1,\partial_1m_in_i}
\end{align*}
for $l_i,{l}^\prime_i\in L_i$ ,  $m_i,{m}^\prime_i\in M_i$ ,  $n_i\in N_i$ ,  $i=1,2,3.$

Next factor each $K(L_i\rtimes M_i\rtimes N_i)$ by the corresponding ideals ${J}^i_2$ and replace ${\sigma}^i_3,{\tau}^i_3$ by induced $\overline{\sigma^i_3},\overline{\tau^i_3}$ and factor each $K(1\rtimes M_i\rtimes N_i)$ by the corresponding ideals ${J}^i_1$ and replace ${\sigma}^i_2, \ {\tau}^i_2$ by the induced $\overline{\sigma^i_2},\   \overline{\tau^i_2}$ for $i=1,2,3$.

This gives us the following diagram;
$$
\xymatrix{K(L_1\rtimes M_1\rtimes N_1)/ {J^1_2}\ar@<1.0ex>[dd]^{\overline{\tau^1_3}}\ar@<-1.0ex>[dd]_{\scriptstyle \overline{\sigma^1_3}}\ar[r]^-{\scriptscriptstyle {\overline K(\phi_3)}}& K(L_2\rtimes M_2\rtimes N_2)/ {J^2_2}\ar@<1.0ex>[dd]^{\overline{\tau^2_3}}\ar@<-1.0ex>[dd]_{\overline{\sigma^2_3}}\ar[r]^-{\scriptscriptstyle {\overline K(\psi_3)}}& K(L_3\rtimes M_3\rtimes N_3)/ {J^3_2}\ar@<1.0ex>[dd]^{\overline{\tau^3_3}}\ar@<-1.0ex>[dd]_{\overline{\sigma^3_3}}\\
\\
K(1\rtimes M_1\rtimes N_1)/ {J^1_1}\ar@<1.0ex>[dd]^{\overline{\tau^1_2}}\ar@<-1.0ex>[dd]_{\overline{\sigma^1_2}}\ar[r]^-{\scriptscriptstyle {\overline K(\phi_2)}}&K(1\rtimes M_2\rtimes N_2)/ {J^2_1}\ar[r]^-{\scriptscriptstyle {\overline K(\psi_2)}}\ar@<1.0ex>[dd]^{\overline{\tau^2_2}}\ar@<-1.0ex>[dd]_{\overline{\sigma^2_2}}&K(1\rtimes M_3\rtimes N_3)/ {J^3_1}\ar@<1.0ex>[dd]^{\overline{\tau^2_3}}\ar@<-1.0ex>[dd]_{\overline{\sigma^3_2}}\\
\\
K(1\rtimes1\rtimes N_1)\ar[r]_-{\scriptscriptstyle \overline{K}(\scriptscriptstyle \phi_1)}&K(1\rtimes1\rtimes N_2)\ar[r]_-{\scriptscriptstyle {\overline K(\psi_1)}}&K(1\rtimes1\rtimes N_3)}
$$
where $\overline{K}(\phi_i)$ and $\overline{K}(\psi_i)$ are induced by the quotient maps. Again, the commutativity is satisfied. Define $$\overline{K}(\mathfrak{C}^2):=\overline{K(\mathfrak{C}^2)}$$
and $\overline{K}(\phi)$ be the Gray 3-(group) algebra groupoid map with $\overline{K}(\phi_i) \ (i=2,3)$ defined as above and
$$\overline{K}(\phi_1):=K(\phi_1) \ , \ \overline{K}(\psi_1):=K(\psi_1) .$$
Therefore, we have
$$\overline{K}(\psi\phi):=\overline{K}(\psi)\overline{K}(\phi)$$
and it is easy to see that $\overline{K}$ also preserves the trivial morphisms on $\mathfrak{C}^2$, and so this gives a functor $\overline{K}$ from the category of Gray 3-(group)-groupoids to that of  Gray 3-(group) algebra groupoids with a single object $*$.

\section{The Construction of Regular Representation as a 3-functor}
In this section, we give the right regular representation (Cayley's theorem) for the 2-crossed module $\xymatrix{\mathfrak{X}:=(L\ar[r]^-{\partial_2}&M\ar[r]^{\partial_1}&N)}$ and so its associated structure $\mathfrak{C}^2$.
Consider the associated Gray 3-(group)-groupoid with a single object or equivalently cat$^2$-group:
$$
\xymatrix{\mathfrak{C}^2:=L\rtimes M \rtimes N \ar@<1ex>[r]^-{s_3,t_3} \ar@<0ex>[r]& 1\rtimes M \rtimes N \ar@<1ex>[r]^-{s_2,t_2} \ar@<0ex>[r]\ar@<1ex>[l]^-{e_3}& 1\rtimes1\rtimes N\ar@<1ex>[l]^-{e_2}\ar@<0.5ex>[r]\ar@<-0.5ex>[r]& \{*\}.}
$$
In the previous section we obtain a chain complex of algebras
$$
\begin{tikzcd}
             \overline{\delta}:=\mathrm{K}_3\ar[r,"\overline{\tau}_3"]&\mathrm{K}_2\ar[r,"\overline{\tau}_2"]& \mathrm{K}_1
\end{tikzcd}
$$
 over the field $K$ from $\overline{K(\mathfrak{C}^2)}$. Consider the cat$^2$-group $\mathbf{Aut(\overline{\delta})}$ from \cite{Jinan}. This can be considered as a subcategory of $\mathbf{Ch}^2_K$. We will define a right regular representation as a lax 3-(contravariant) functor $$\lambda:\mathfrak{C}^2 \rightarrow \mathbf{Aut(\overline{\delta})}.$$
 This functor $\lambda$ will be organised as follows:
\begin{align*}
 \lambda:\mathfrak{C}^2& \xrightarrow{\makebox[2cm]{}}  \mathbf{Aut(\overline{ \delta})}\\
 *&\xmapsto{\hspace{2cm}}\overline{\delta}\ (\text{single object})\\
 (1,1,n)&\xmapsto{\hspace{2cm}} \lambda_n=({\lambda}^0_n,{\lambda}^1_n,{\lambda}^2_n)\ (\text{chain map})\\
 (1,m,n)&\xmapsto{\hspace{2cm}} \lambda_{m,n}=(({\lambda}^\prime_{m,n},{\lambda}^{\prime \prime}_{m,n}),\lambda_n)\ (\text{1-homotopy})\\
 (l,m,n)&\xmapsto{\hspace{2cm}}\lambda_{l,m,n}=(\alpha_{l,m,n},({\lambda}^\prime_{m,n},{\lambda}^{\prime \prime}_{m,n}),\lambda_n) \ (\text{2-homotopy}).
\end{align*}

\textbf{Remark:} We can summarise that any representation of $\mathfrak{C}^2$ maps elements of $1\rtimes 1 \rtimes N$ to chain automorphisms in $\mathbf{(Aut\overline{ \delta})}_{1}$, elements of $1\rtimes M \rtimes N$ to 1-homotopies in $\mathbf{(Aut\overline{ \delta})}_{2}$ and elements of $L\rtimes M \rtimes N$ to 2-homotopies in $\mathbf{(Aut\overline{ \delta})}_{3}$ for some representation complex on length-2 of algebras (or vector spaces by ignoring the multiplication) $\overline{\delta}$. We can picture an action of $\mathfrak{C}^2$ on $\overline{K(\mathfrak{C}^2)}$ by right multiplication. Of course, this is left action, and the elements of $\mathfrak{C}^2$ appear in the left side on the picture, while they appear on the right in the algebraic notation. In the following pictures, the broken arrows show elements of $\mathfrak{C}^2$ and the unbroken arrows are as in $\overline{K(\mathfrak{C}^2)}$.

A 1-cell in $\mathfrak{C}^2$ is an element $(1,1,n)\in 1\rtimes 1\rtimes N$. This can act both on the 1-cells, 2-cells and 3-cells of $\overline{K(\mathfrak{C}^2)}$. The action of $(1,1,n)$ on a 1-cell is;
$$
\begin{tikzcd}
              \ar[rr,densely dotted,"{\scriptscriptstyle{(1,1,n)}}"]&&{}\ar[rr,"{\mathbf{e}_{1,1,n'}}"]&& := \ar[rr,"{\mathbf{e}_{1,1,n'n}}"]&&{}
\end{tikzcd}
$$
The action on a 2-cell is similar:
$$
\begin{array}{c}
    \begin{tikzcd}[row sep=small,column sep=normal]
                    &  &  \ar[dd,Rightarrow,"{\overline{\mathbf{e}}_{1,m',n'}}"{description}] & \\
                \ar[r,"{\scriptscriptstyle{ (1,1,n)}}",densely dotted]& \ar[rr,bend left=40,"{\mathbf{e}_{1,1,n'}}"] \ar[rr,bend right=50,"{\mathbf{e}_{1,1,\partial_1m'n'}}"'] &  & {} \\
                    \ &  \ & \
    \end{tikzcd}
\end{array}
{:=}
\begin{array}{c}
    \begin{tikzcd}[row sep=small,column sep=normal]
                 & \ar[dd, Rightarrow, "{\overline{\mathbf{e}}_{1,m',n'n}}"{description}] & \\
                 \ar[rr,bend left=40,"{\mathbf{e}_{1,1,n'n}}"] \ar[rr,bend right=50,"{\mathbf{e}_{1,1,\partial_1m'n'n}}"'] &  &{}\\
                 & \
    \end{tikzcd}
\end{array}
$$
The action of $(1,m,n)$ on a 1-cell is given by
$$
\begin{array}{c}
    \begin{tikzcd}[row sep=small,column sep=normal]
                    &  \ar[Rightarrow,"{\scriptscriptstyle (1,m,n)}",dotted]{dd} & &  \\
                 \ar[rr,bend left=40,"{\scriptscriptstyle{(1,1,n)}}",dotted] \ar[rr,bend right=50,"{\scriptscriptstyle{(1,1,\partial_1 mn)}}"',dotted] &  &\ar[r,"{\mathbf{e}_{1,1,n'}}"] & {}  \\
                \ & \ &   \ &
    \end{tikzcd}
\end{array}
{:=}
\begin{array}{c}
    \begin{tikzcd}[row sep=small,column sep=normal]
                 &   \ar[Rightarrow,"{\overline{\mathbf{e}}_{1,^{n'}m,n'n}}"{description}]{dd}&\\
                  \ar[rr,bend left=40,"{\mathbf{e}_{1,1,n'n}}"] \ar[rr,bend right=50,"{\mathbf{e}_{1,1,n'\partial_1mn}}"'] &  & {} \\
                 \ & \ &
    \end{tikzcd}
\end{array}
$$
The action of $(1,1,n)$ on a 3-cell is given by:
$$
\begin{array}{c}
    \begin{tikzcd}[row sep=small,column sep=normal]
                    &  &    \tarrow["\scriptscriptstyle{\overline{\mathbf{e}}_{l',m',n'}}"]{dd} & \\
                \ar[r,"\scriptscriptstyle {(1,1,n)}",densely dotted]  &\ar[rr,Rightarrow,bend left=40,"{\mathbf{e}_{1,m',n'}}"] \ar[rr,Rightarrow,bend right=50,"{\mathbf{e}_{1,\partial_2l'm',n'}}"'] &  &{} \\
                \ & \ & \
    \end{tikzcd}
\end{array}
{:=}
\begin{array}{c}
    \begin{tikzcd}[row sep=small,column sep=normal]
                 &   \tarrow["{\scriptscriptstyle{\overline{\mathbf{e}}_{l',m',n'n}}}" {description} ]{dd}&\\
                  \ar[rr,Rightarrow,bend left=40,"{\mathbf{e}_{1,m',n'n}}"] \ar[rr,Rightarrow,bend right=50,"{\mathbf{e}_{1,\partial_2l'm',n'n}}"'] &  & {} \\
                 \ & \ &
    \end{tikzcd}
\end{array}
$$
The action of $(1,m,n)$ on a 3-cell is given by:
$$
\begin{array}{c}
    \begin{tikzcd}[row sep=small,column sep=normal]
                    &  &    \tarrow["{\scriptscriptstyle{\overline{\mathbf{e}}_{l',m',n'}}}"]{dd} & \\
                \ar[r,Rightarrow,"\scriptscriptstyle {(1,m,n)}",dotted] & \ar[rr,Rightarrow,bend left=40,"{\mathbf{e}_{1,m',n'}}"] \ar[rr,Rightarrow,bend right=50,"{\mathbf{e}_{1,\partial_2l'm',n'}}"'] &  &{} \\
                \ & \ & \
    \end{tikzcd}
\end{array}
{:=}
\begin{array}{c}
    \begin{tikzcd}[row sep=small,column sep=normal]
                 &   \tarrow["{\scriptscriptstyle{\overline{\mathbf{e}}_{l',m'^{n'}m,n'n}}}"{description}]{dd}&\\
                  \ar[rr,Rightarrow,bend left=40,"{\mathbf{e}_{1,m'^{n'}m,n'n}}"] \ar[rr,Rightarrow,bend right=50,"{\mathbf{e}_{1,\partial_2l'm'^{n'}m,n'n}}"'] &  & {} \\
                 \ & \ &
    \end{tikzcd}
\end{array}
$$
The action of $(l,m,n)$ on a 1-cell is given by
$$
\begin{array}{c}
    \begin{tikzcd}[row sep=small,column sep=normal]
                    &   \tarrow["{\scriptscriptstyle(l,m,n)}",dotted]{dd} & &  \\
                 \ar[rr,Rightarrow,bend left=40,"{\scriptscriptstyle{(1,m,n)}}",dotted] \ar[rr,Rightarrow,bend right=50,"{\scriptscriptstyle{(1,\partial_2lm,n)}}"',dotted] &  &\ar[r,"{\mathbf{e}_{1,1,n'}}"] & {} \\
                \ & \ & \  &
    \end{tikzcd}
\end{array}
{:=}
\begin{array}{c}
    \begin{tikzcd}[row sep=small,column sep=normal]
                 &   \tarrow["{\scriptscriptstyle{\overline{\mathbf{e}}_{^{n'}l,^{n'}m,n'n}}}"{description}]{dd}&\\
                  \ar[rr,Rightarrow,bend left=40,"{\mathbf{e}_{1,^{n'}m,n'n}}"] \ar[rr,Rightarrow,bend right=50,"{\mathbf{e}_{1,^{n'}(\partial_2lm),n'n}}"'] &  & {} \\
                 \ & \ &
    \end{tikzcd}
\end{array}
$$
Using these actions,  we can give the construction of the 3-functor $\lambda$, on each level as follows:
\subsection{$\lambda$ over 0-cells}
Since the 0-cell of $\mathfrak{C}^2$ is $*$, and the 0-cell of  $\mathbf{Aut(\overline{ \delta})}$ is $\overline{\delta}$, so  we can define $\lambda(*)=\overline{\delta}$.
Thus, we have $$\lambda_{*}=\overline{\delta} \in \mathbf{(Aut \overline{ \delta})}_{0}.$$
\subsection{$\lambda$ over 1-cells}
For any 1-cell of $\mathfrak{C}^2$; $(1,1,n):*\rightarrow *$ for $n\in N$, using the action of $(1,1,n)$ on the cells in $\overline{K(\mathfrak{C}^2)}$ given in the above pictures and considering only $\lambda_n$ for the image of $(1,1,n)$ under $\lambda$, we can define $\lambda_n:\overline{\delta}\rightarrow \overline{\delta}$ as follows:
$$
\begin{array}{c}
\xymatrix{\overline{\delta}\ar[d]_-{\lambda_n}\\
\overline{\delta}}
\end{array}
 {:=}
\begin{array}{c}
\xymatrix{\mathrm{K}_3 \ar[r]^-{\overline{\tau}_3}\ar[d]^{{\lambda}^2_n}& \mathrm{K}_2 \ar[r]^-{\overline{\tau}_2}\ar[d]^-{{\lambda}^1_n}&\mathrm{K}_1\ar[d]^-{{\lambda}^0_n}\\
\mathrm{K}_3\ar[r]_-{\overline{\tau}_3} & \mathrm{K}_2 \ar[r]_-{\overline{\tau}_2}  &\mathrm{K}_1 }
\end{array}
$$
where
$${\lambda}^0_n(\mathbf{e}_{1,1,n^\prime})=\mathbf{e}_{1,1,n^\prime n} \ ,\  \ {\lambda}^1_n (\overline{\mathbf{v}^{\scriptscriptstyle11}}_{1,m^\prime,n^\prime})= \overline{\mathbf{v}^{\scriptscriptstyle11}}_{1,m^\prime,n^\prime n}\ ,\ \ {\lambda}^2_n(\overline{\mathbf{v}^{\scriptscriptstyle22}}_{l',m',n'})=\overline{\mathbf{v}^{\scriptscriptstyle22}}_{l',m',n'n}$$
for each $\mathbf{e}_{1,1,n'}\in \mathrm{K}_1$,\ \ $\overline{\mathbf{v}^{\scriptscriptstyle11}}_{1,m',n'}\in\mathrm{K}_2$ and $\overline{\mathbf{v}^{\scriptscriptstyle22}}_{l',m',n'}\in \mathrm{K}_3$.

 The maps ${\lambda}^0_n,{\lambda}^1_n$ and ${\lambda}^2_n$ are linear automorphisms from $\mathrm{K}_i$ to $\mathrm{K}_i$ for $i=1,2,3$. We will show that $\lambda_n=({\lambda}^0_n,{\lambda}^1_n,{\lambda}^2_n)$ is a chain automorphism in  $\mathbf{Aut(\overline{ \delta})}$. It must be that the last diagram is commutative.
For any element $\overline{\mathbf{v}^{\scriptscriptstyle11}}_{1,m',n'}\in\mathrm{K}_2$, we obtain
\begin{align*}
\overline\tau_2 {\lambda}^1_n(\overline{\mathbf{v}^{\scriptscriptstyle11}}_{1,m',n'})=&\overline\tau_2(\overline{\mathbf{v}^{\scriptscriptstyle11}}_{1,m',n'n})\\
=&\overline\tau_2(\overline{\mathbf{e}}_{1,m',n'n}-\overline{\mathbf{e}}_{1,1,n'n})\\
=&\mathbf{e}_{1,1,\partial_1m'n'n}-\mathbf{e}_{1,1,n'n}\in \mathrm{K}_1
\end{align*}
and
\begin{align*}
{\lambda}^0_n \overline\tau_2(\overline{\mathbf{v}^{\scriptscriptstyle11}}_{1,m',n'})=&{\lambda}^0_n \overline\tau_2(\overline{\mathbf{e}}_{1,m',n'}-\overline{\mathbf{e}}_{1,1,n'})\\
=&{\lambda}^0_n(\mathbf{e}_{1,1,\partial_1m'n'}-\mathbf{e}_{1,1,n'})\\
=&\mathbf{e}_{1,1,\partial_1m'n'n}-\mathbf{e}_{1,1,n'n}\in \mathrm{K}_1
\end{align*}
and so $$\overline\tau_2 \lambda^1_n=\lambda^0_n \overline\tau_2.$$
For any $\overline{\mathbf{v}^{\scriptscriptstyle22}}_{l',m',n'}\in \mathrm{K}_3$, we obtain
\begin{align*}
\overline\tau_3 \lambda^2_n(\overline{\mathbf{v}^{\scriptscriptstyle22}}_{l',m',n'})=&\overline\tau_3(\lambda^2_n(\overline{\mathbf{e}}_{l',m',n'}-\overline{\mathbf{e}}_{1,m',n'}))\\
=&\overline\tau_3(\overline{\mathbf{e}}_{l',m',n'n}-\overline{\mathbf{e}}_{1,m',n'n})\\
=&\overline{\mathbf{e}}_{1,\partial_2l'm',n'n}-\overline{\mathbf{e}}_{1,m',n'n})\\
=&(\overline{\mathbf{e}}_{1,\partial_2l'm',n'n}-\overline{\mathbf{e}}_{1,1,n'n})-(\overline{\mathbf{e}}_{1,m',n'n}-\overline{\mathbf{e}}_{1,1,n'n})\\
=&\overline{\mathbf{v}^{\scriptscriptstyle11}}_{1,\partial_2l'm',n'n}-\overline{\mathbf{v}^{\scriptscriptstyle11}}_{1,m',n'n}\in \mathrm{K}_2
\end{align*}
and;
\begin{align*}
\lambda^1_n  \overline\tau_3(\overline{\mathbf{v}^{\scriptscriptstyle22}}_{l',m',n'})=&\lambda^1_n(\overline\tau_3(\overline{\mathbf{e}}_{l',m',n'}-\overline{\mathbf{e}}_{1,m',n'}))\\
=&\lambda^1_n (\overline{\mathbf{e}}_{1,\partial_2l'm',n'}-\overline{\mathbf{e}}_{1,m',n'})\\
=&\lambda^1_n ((\overline{\mathbf{e}}_{1,\partial_2l'm',n'}-\overline{\mathbf{e}}_{1,1,n'})-(\overline{\mathbf{e}}_{1,m',n'}-\overline{\mathbf{e}}_{1,1,n'}))\\
=&\lambda^1_n (\overline{\mathbf{v}^{\scriptscriptstyle11}}_{1,\partial_2l'm',n'}-\overline{\mathbf{v}^{\scriptscriptstyle11}}_{1,m',n'})\\
=&\overline{\mathbf{v}^{\scriptscriptstyle11}}_{1,\partial_2l'm',n'n}-\overline{\mathbf{v}^{\scriptscriptstyle11}}_{1,m',n'n}\in \mathrm{K}_2
\end{align*}
and so;
$$\overline\tau_3\lambda^2_n=\lambda^1_n \overline\tau_3.$$
Therefore, $\lambda_n=({\lambda}^0_n,{\lambda}^1_n,{\lambda}^2_n)$ is a chain automorphism over $\overline{\delta}$ in $\mathbf{(Aut \overline{ \delta})}_{1}.$
\subsubsection{$\lambda_n$ preserves the composition of 1-cells}
For any 1-cells  $n,n':{*} \longrightarrow {*} $ in $\mathfrak{C}^2$, and for any ${\mathbf{e}}_{1,1,p}\in \mathrm{K}_1$,
 we have $$\lambda^0_{nn'}({\mathbf{e}}_{1,1,p})={\mathbf{e}}_{1,1,p(nn')}={\mathbf{e}}_{1,1,(pn)n'}$$ and
$$({\lambda}^0_{n'}\circ{\lambda}^0_n)({\mathbf{e}}_{1,1,p})={\lambda}^0_{n'}({\lambda}^0_n({\mathbf{e}}_{1,1,p}))={\lambda}^0_{n'}({\mathbf{e}}_{1,1,pn})={\mathbf{e}}_{1,1,(pn)n'}$$ and so $${\lambda}^0_{nn'}={\lambda}^0_{n'}\circ{\lambda}^0_n.$$
For any $\overline{\mathbf{v}^{\scriptscriptstyle11}}_{1,k,p}\in\mathrm{K}_2$, we have
\begin{align*}
\lambda'_{nn'}({\mathbf{v}^{\scriptscriptstyle11}}_{1,k,p})=&\overline{\mathbf{v}^{\scriptscriptstyle11}}_{1,k,p(nn')}\\
=&\overline{\mathbf{v}^{\scriptscriptstyle11}}_{1,k,(pn)n'}\\
=&\lambda'_{n'}(\overline{\mathbf{v}^{\scriptscriptstyle11}}_{1,k,pn})\\
=&\lambda'_{n'}(\lambda'_{n}(\overline{\mathbf{v}^{\scriptscriptstyle11}}_{1,k,p}))\\
=&({\lambda}^1_{n'}\circ{\lambda}^1_n)({\mathbf{v}^{\scriptscriptstyle11}}_{1,k,p})
\end{align*}
and so ${\lambda}^1_{nn'}={\lambda}^1_{n'}\circ{\lambda}^1_n $. Similarly, we have ${\lambda}^2_{nn'}={\lambda}^2_{n'}\circ{\lambda}^2_n$.
\subsection{$\lambda$ over 2-cells}
For any 2-cell $(1,m,n)\in 1\rtimes  M\rtimes N$, there must be a corresponding 2-cell $\lambda_{m,n}\in \mathbf{(Aut\overline{ \delta})}_{2}$ as a 1-homotopy. $\lambda$  need to be a functor, so it must preserve the source and target of each 2-cell and so, $\lambda_{m,n}$ will be a homotopy in $\mathbf{Ch}^2_K$ from $\lambda_n$ to $\lambda_{\partial_1mn}$. We can show it by its source and chain homotopy component as $\lambda_{m,n}:=((\lambda'_{m,n},\lambda''_{m,n}),\lambda_n)$, where $\lambda'_{m,n}$,\ $\lambda''_{m,n}$ are the chain homotopy components and $\lambda_n$ is the source of $\lambda_{m,n}$.

Here, the chain homotopy components are $\lambda'_{m,n}:\mathrm{K}_1\longrightarrow \mathrm{K}_2$ and  $\lambda''_{m,n}:\mathrm{K}_2\longrightarrow \mathrm{K}_3$. These maps can be shown in the following diagram:
\begin{equation*}
\xymatrix{\mathrm{K}_3 \ar[rr]^-{\overline{\tau}_3}\ar@<0.5ex>[dd]^-{\scriptscriptstyle\lambda^2_{\partial_1mn}}\ar@<-0.5ex>[dd]_-{{\scriptscriptstyle\lambda}^2_n} &&\mathrm{K}_2 \ar[ddll]|<<<<<<<{\Large\boldsymbol\lambda''_{m,n}} \ar[rr]^-{\overline{\tau}_2}\ar@<0.5ex>[dd]^-{\scriptscriptstyle\lambda^1_{\partial_1mn}}\ar@<-0.5ex>[dd]_-{\scriptscriptstyle{\lambda}^1_n}&&\mathrm{K}_1\ar[ddll]|<<<<<<<{\Large\boldsymbol\lambda'_{m,n}}\ar@<0.5ex>[dd]^-{\scriptscriptstyle\lambda^0_{\partial_1mn}}\ar@<-0.5ex>[dd]_-{\scriptscriptstyle{\lambda}^0_n}
\\
\\
\mathrm{K}_3\ar[rr]_{\overline{\tau}_3} && \mathrm{K}_2 \ar[rr]_{\overline{\tau}_2} \ar[rr] &&\mathrm{K}_1.}
\end{equation*}
We must check that $\lambda_{m,n}$ satisfies the chain homotopy conditions given below:
\begin{enumerate}
\item $\overline\tau_2\lambda'_{m,n}=\lambda^0_{\partial_1mn}-\lambda^0_{n}$,
\item $\lambda'_{m,n}\overline\tau_2+\overline\tau_3 \lambda''_{m,n}=\lambda^1_{\partial_1mn}-\lambda^1_{n},$
\item $\lambda''_{m,n}\overline\tau_3=\lambda^2_{\partial_1mn}-\lambda^2_{n}.$
\end{enumerate}
We define, for any ${\mathbf{e}}_{1,1,n'}\in \mathrm{K}_1$, the map $\lambda'_{m,n}:\mathrm{K}_1\longrightarrow \mathrm{K}_2$ by
$$\lambda'_{m,n}({\mathbf{e}}_{1,1,n'})=\overline{\mathbf{v}^{\scriptscriptstyle11}}_{(1,1,n')(1,m,n)}=\overline{\mathbf{v}^{\scriptscriptstyle11}}_{1,^{n'}m,n'n}$$
and for any $\overline{\mathbf{v}^{\scriptscriptstyle11}}_{1,m',n'}\in \mathrm{K}_2$, the map $\lambda''_{m,n}:\mathrm{K}_2\longrightarrow \mathrm{K}_3$ by
$$\lambda''_{m,n}(\overline{\mathbf{v}^{\scriptscriptstyle11}}_{1,m',n'})=\overline{\mathbf{v}^{\scriptscriptstyle22}}_{(1,m',n')\#(1,m,n)}
=\overline{\mathbf{v}^{\scriptscriptstyle22}}_{\{m',^{n'}m\},m'^{n'}m,n'n},$$
where $(1,m',n')\#(1,m,n)=({\{m',^{n'}m\},m'^{n'}m,n'n})$ is the interchange 3-cell of 2-cells $(1,m,n)$ and $(1,m',n')$ in $\mathfrak{C}^2$ and
$\{-,-\}:M\times M\longrightarrow L$ is the Peiffer lifting of 2-crossed module $\mathfrak{X}$.

We shall check the first condition;

$1.$ For any ${\mathbf{e}}_{1,1,n'}\in \mathrm{K}_1;$ we obtain
\begin{align*}
\overline\tau_2 \lambda'_{m,n}({\mathbf{e}}_{1,1,n'})=&\overline\tau_2(\overline{{\mathbf{v}}^{\scriptscriptstyle11}}_{1,^{n'}m,n'n})\\
=&\overline\tau_2(\overline {\mathbf{e}}_{1,^{n'}m,n'n}-\overline {\mathbf{e}}_{1,1,n'n})\\
=&{\mathbf{e}}_{1,1,\partial_1(^{n'}m)n'n}-{\mathbf{e}}_{1,1,n'n}\\
=&{\mathbf{e}}_{1,1,n'\partial_1mn}-{\mathbf{e}}_{1,1,n'n}\\
=&\lambda^0_{\partial_1mn}({\mathbf{e}}_{1,1,n'})-\lambda^0_n({\mathbf{e}}_{1,1,n'})\\
=&(\lambda^0_{\partial_1mn}-\lambda^0_n)({\mathbf{e}}_{1,1,n'}).
\end{align*}
Thus, we have $\overline\tau_2 \lambda'_{m,n}=\lambda^0_{\partial_1mn}-\lambda^0_n$.

$2.$
We must show that
$$\lambda'_{m,n} \overline\tau_2+\overline\tau_3 \lambda''_{m,n}=\lambda'_{\partial_1mn}-\lambda'_n.$$
For any element $\overline{\mathbf{v}^{\scriptscriptstyle11}}_{1,m',n'}\in \mathrm{K}_2$, we obtain
$$\overline\tau_3 \lambda''_{m,n}(\overline{\mathbf{v}^{\scriptscriptstyle11}}_{1,m',n'})=\overline\tau_3(\overline{\mathbf{v}^{\scriptscriptstyle22}}_{\{m',^{n'}m\},m'^{n'}m,n'n})=\overline{\mathbf{v}^{\scriptscriptstyle11}}_{1,m'^{n'}m,n'n}-\overline{\mathbf{v}^{\scriptscriptstyle11}}_{1,\partial_2 \{m',^{n'}m\}m'^{n'}m,n'n}.$$
Since $\partial_2\{m',^{n'}m\}=^{\partial_1m'}(^{n'}m)m'(^{n'}m)^{-1}{m'}^{-1}$, we obtain
$$\overline{\mathbf{v}^{\scriptscriptstyle11}}_{1,\partial_2 \{m',^{n'}m\}m'(^{n'}m),n'n}=\overline{\mathbf{v}^{\scriptscriptstyle11}}_{1,^{\partial_1m'}(^{n'}m)m',n'n}$$ and then
$$\overline\tau_3\lambda''_{m,n}(\overline{\mathbf{v}^{\scriptscriptstyle11}}_{1,m',n'})=\overline{\mathbf{v}^{\scriptscriptstyle11}}_{1,m'^{n'}m,n'n}-\overline{\mathbf{v}^{\scriptscriptstyle11}}_{1,^{\partial_1m'}(^{n'}m)m',n'n}$$
On the other hand,
$$(\lambda'_{\partial_1mn}-\lambda'_n)(\overline{\mathbf{v}^{\scriptscriptstyle11}}_{1,m',n'})=\lambda'_{\partial_1mn}(\overline{\mathbf{v}^{\scriptscriptstyle11}}_{1,m',n'})-\lambda'_n(\overline{\mathbf{v}^{\scriptscriptstyle11}}_{1,m',n'})=\overline{\mathbf{v}^{\scriptscriptstyle11}}_{1,m',n'\partial_1mn}-\overline{\mathbf{v}^{\scriptscriptstyle11}}_{1,m',n'n}$$
and
\begin{align*}
  \lambda'_{m,n}\overline\tau_2(\overline{\mathbf{v}^{\scriptscriptstyle11}}_{1,m',n'})=&\lambda'_{m,n}\overline\tau_2(\overline {\mathbf{e}}_{1,m',n'}-\overline {\mathbf{e}}_{1,1,n'})\\
  =&\lambda'_{m,n}(\mathbf{e}_{1,1,\partial_1m'n'}-\mathbf{e}_{1,1,n'})\\
  =&\overline{\mathbf{v}^{\scriptscriptstyle11}}_{1,^{\partial_1m'}(^{n'}m),\partial_1m'n'n}-\overline{\mathbf{v}^{\scriptscriptstyle11}}_{1,^{n'}m,n'n}
\end{align*}
Thus we obtain
\begin{align*}
(\lambda'_{\partial_1mn}-\lambda'_{n}-\lambda'_{m,n}\overline\tau_2)(\overline{\mathbf{v}^{\scriptscriptstyle11}}_{1,m',n'})=&\overline{\mathbf{v}^{\scriptscriptstyle11}}_{1,m',n'\partial_1mn}-\overline{\mathbf{v}^{\scriptscriptstyle11}}_{1,m',n'n}+\overline{\mathbf{v}^{\scriptscriptstyle11}}_{1,^{n'}m,n'n}-\overline{\mathbf{v}^{\scriptscriptstyle11}}_{1,^{\partial_1m'}(^{n'}m),\partial_1m'n'n}\\
=&(\overline{\mathbf{v}^{\scriptscriptstyle11}}_{1,m',n'\partial_1mn}+\overline{\mathbf{v}^{\scriptscriptstyle11}}_{1,^{n'}m,n'n})-(\overline{\mathbf{v}^{\scriptscriptstyle11}}_{1,^{\partial_1m'}(^{n'}m),\partial_1m'n'n}+\overline{\mathbf{v}^{\scriptscriptstyle11}}_{1,m',n'n}).
\end{align*}
We know that
$${\mathbf{u}_{\scriptscriptstyle2}}=\mathbf{e}_{1,m'm,n}-\mathbf{e}_{1,m,n}-\mathbf{e}_{1,m',\partial_1mn}+\mathbf{e}_{1,1,\partial_1mn}\in J_1.$$
We can rewrite this expression as follows:
$$(\mathbf{e}_{1,m'm,n}-\mathbf{e}_{1,1,n})-(\mathbf{e}_{1,m,n}-\mathbf{e}_{1,1,n})-(\mathbf{e}_{1,m',\partial_1mn}-\mathbf{e}_{1,1,\partial_1mn})=
{\mathbf{v}^{\scriptscriptstyle11}}_{1,m'm,n}-{\mathbf{v}^{\scriptscriptstyle11}}_{1,m,n}-{\mathbf{v}^{\scriptscriptstyle11}}_{1,m',\partial_1mn}\in J_1$$
and thus, we have $$\overline{\mathbf{v}^{\scriptscriptstyle11}}_{1,m'm,n}=\overline{\mathbf{v}^{\scriptscriptstyle11}}_{1,m,n}+\overline{\mathbf{v}^{\scriptscriptstyle11}}_{1,m',\partial_1mn}.$$
Using this equality, we obtain
$$\overline{\mathbf{v}^{\scriptscriptstyle11}}_{1,m',n'\partial_1mn}+\overline{\mathbf{v}^{\scriptscriptstyle11}}_{1,^{n'}m,n'n}=\overline{\mathbf{v}^{\scriptscriptstyle11}}_{1,m'^{n'}m,n'n}$$
and
$$\overline{\mathbf{v}^{\scriptscriptstyle11}}_{1,m',n'n}+\overline{\mathbf{v}^{\scriptscriptstyle11}}_{1,^{\partial_1m'}(^{n'}m),\partial_1m'n'n}=\overline{\mathbf{v}^{\scriptscriptstyle11}}_{1,^{\partial_1m'}(^{n'}m)m',n'n}.$$
Thus, we obtain
$$(\lambda^1_{\partial_1mn}-\lambda^1_{n}-\lambda'_{m,n} \overline\tau_2)(\overline{\mathbf{v}^{\scriptscriptstyle11}}_{1,m',n'})=\overline{\mathbf{v}^{\scriptscriptstyle11}}_{1,m'^{n'}m,n'n}-\overline{\mathbf{v}^{\scriptscriptstyle11}}_{1,^{\partial_1m'}(^{n'}m)m',n'n}=
\overline\tau_3(\lambda''_{m,n}(\overline{\mathbf{v}}^{11}_{1,m',n'})).$$
Therefore we have;
$\lambda^1_{\partial_1mn}-\lambda^1_{n}-\lambda'_{m,n} \overline\tau_2=\overline\tau_3 \lambda''_{m,n}$ and thus
$\overline\tau_3\lambda''_{m,n}+\lambda'_{m,n}\overline\tau_2=\lambda^1_{\partial_1mn}-\lambda^1_{n}.$
Thus, second chain homotopy condition is satisfied. In section \ref{sect:App5} of Appendix, we will show that third chain homotopy condition; $\lambda'' \overline\tau_3=\lambda^2_{\partial_1mn}-\lambda^2_{n}$ is satisfied.

\subsubsection{$\lambda_{m,n}$ preserves the vertical composition of 2-cells}
We will show that $\lambda_{m,n}$ preserves the vertical composition of 2-cells. We must show that $$\lambda_{(m',\partial_1mn)\#_2(m,n)}=\lambda_{m'm,n}=\lambda_{(m',\partial_1mn)}\#_2\lambda_{(m,n)}$$
for the vertical composition $\#_2$ of 2-cells in $\mathfrak{C}^2$.

Consider the following diagram;
$$
 \begin{tikzcd}
\mathrm{K}_3 \arrow[rrr,shift left=0.90ex,"{\scriptscriptstyle {{\lambda}^2_n}}"description,pos=.75]\arrow[rrr,shift left=-0.90ex,"{\scriptscriptstyle \lambda^2_{\partial_1mn}}"description,pos=.35]    \arrow[ddd,"{\overline{\tau}_3}"'] &&&\mathrm{K}_3\arrow[rrr,shift left=0.95ex,"{\scriptscriptstyle \lambda^2_{\partial_1mn}}"description,pos=.70]\arrow[rrr,shift left=-0.95ex,"{\scriptscriptstyle{\lambda}^2_{\partial_1m'n'}}"description,pos=.35]    \arrow[ddd,"{\overline{\tau}_3}"]                &&& \mathrm{K}_3 \arrow[ddd,,"{\overline{\tau}_3}"]  \\ \\ \\
\mathrm{K}_2 \arrow[rrr,shift left=0.90ex,"{\scriptscriptstyle {{\lambda}^1_n}}"description,pos=.75]\arrow[rrr,shift left=-0.90ex,"{\scriptscriptstyle \lambda^1_{\partial_1mn}}"description,pos=.35]  \arrow[ddd,"{\overline{\tau}_2}"']\arrow[uuurrr,"{\scriptscriptstyle{H'_2}}"{sloped,above=-0.3ex,xshift=0.0em}]&&&\mathrm{K}_2 \arrow[rrr,shift left=0.95ex,"{\scriptscriptstyle \lambda^1_{\partial_1mn}}"description,pos=.70]\arrow[rrr,shift left=-0.95ex,"{\scriptscriptstyle{\lambda}^1_{\partial_1m'n'}}"description,pos=.35]  \arrow[ddd,"{\overline{\tau}_2}"]\arrow[uuurrr,"{\scriptscriptstyle{K'_2}}"{sloped,above=-0.3ex,xshift=0.0em}]             &&&\mathrm{K}_2  \arrow[ddd,"{\overline{\tau}_2}"]\\ \\ \\
\mathrm{K}_1 \arrow[rrr,shift left=0.90ex,"{\scriptscriptstyle {{\lambda}^0_n}}"description,pos=.75]\arrow[rrr,shift left=-0.90ex,"{\scriptscriptstyle \lambda^0_{\partial_1mn}}"description,pos=.35] \arrow[uuurrr,"{\scriptscriptstyle {H'_1}}"{sloped,above=-0.3ex,xshift=0.0em}] &&&\mathrm{K}_1 \arrow[rrr,shift left=0.95ex,"{\scriptscriptstyle \lambda^{0}_{\partial_{1}mn}}"description,pos=.70]\arrow[rrr,shift left=-0.95ex,"{\scriptscriptstyle{\lambda}^{0}_{\partial_{1}m'n'}}"description,pos=.35]\arrow[uuurrr,"{\scriptscriptstyle{K'_1}}"{sloped,above=-0.3ex,xshift=0.0em}]   &&&\mathrm{K}_1
\end{tikzcd}
$$
where $n'=\partial_1 mn$. We can take
$$(H,F)=((H'_1,H'_2),F)=((\lambda'_{m,n},\lambda''_{m,n}),\lambda_n):F\Rightarrow G$$
and
$$F=\lambda_n=({\lambda}^0_{n},{\lambda}^1_{n},{\lambda}^2_{n}); \ \ G=\lambda_{\partial_1mn}=({\lambda}^0_{\partial_1mn},{\lambda}^1_{\partial_1mn},{\lambda}^2_{\partial_1mn})$$ and
$$(K,G)=((K'_1,K''_2),G)=((\lambda'_{m',\partial_1mn},\lambda''_{m',\partial_1mn}),G):G\Rightarrow T$$
where
$$T=\lambda_{\partial_1m'\partial_1mn}=({\lambda}^0_{\partial_1m'\partial_1mn},{\lambda}^1_{\partial_1m'\partial_1mn},{\lambda}^2_{\partial_1m'\partial_1mn}).$$
In $\mathbf{(Aut\overline{ \delta})}_{2}$, the vertical composition of 1-homotopies is given by
$$(K,G)\#_2(H,F)=(K+H,F)=((K'_1+H'_1,K'_2+H'_2),F).$$
For any $(\mathbf{e}_{1,1,p})\in \mathrm{K}_1$; we obtain
$\lambda'_{m'm,n}(\mathbf{e}_{1,1,p})=\overline{\mathbf{v}^{\scriptscriptstyle11}}_{1,^{p}{(m'm)},pn}\in\mathrm{K}_2$.
From the relation ${\mathbf{u}_{\scriptscriptstyle2}}\in J_1$, we can write
$$\overline{\mathbf{v}^{\scriptscriptstyle11}}_{1,^{p}{(m'm)},pn}
=\overline{\mathbf{v}^{\scriptscriptstyle11}}_{1,^{p}m',p\partial_1mn}+\overline{\mathbf{v}^{\scriptscriptstyle11}}_{1,^{p}m,pn}$$
and then
\begin{align*}
\lambda'_{(m',\partial_1mn)\#_2(m,n)}(\mathbf{e}_{1,1,p})=&\lambda'_{m'm,n}(\mathbf{e}_{1,1,p})\\
=&\overline{\mathbf{v}^{\scriptscriptstyle11}}_{1,^{p}{(m'm)},pn}\\
=&\overline{\mathbf{v}^{\scriptscriptstyle11}}_{1,^{p}m',p\partial_1mn}+\overline{\mathbf{v}^{\scriptscriptstyle11}}_{1,^{p}m,pn}\\ =&\lambda'_{m',\partial_1mn}(\mathbf{e}_{1,1,p})+\lambda'_{m,n}(\mathbf{e}_{1,1,p})\\
=&(K'_1+H'_1)(\mathbf{e}_{1,1,p})
\end{align*}
thus $\lambda'_{(m',\partial_1mn)\#_2(m,n)}=K'_1+H'_1$. We  show the equality $\lambda''_{(m',\partial_1mn)\#_2(m,n)}=K'_2+H'_2$ in section \ref{sect:App1} of Appendix.
Therefore, $\lambda_{m,n}$ preserves the vertical composition of 2-cells.

\subsubsection{$\lambda_{m,n}$ preserves the horizontal composition of 2-cells }
For any 2-cells $\Gamma=(m,n)$ and $\Gamma'=(m',n')$ in $M\rtimes N$; their horizontal compositions are
$$\begin{bmatrix} &\Gamma'\\ \Gamma& \end{bmatrix}=(m^{n}m',nn') \  \  \mathrm{and} \ \
\begin{bmatrix} \Gamma&\\ &\Gamma' \end{bmatrix}=\left(^{\partial_1m}{(^{n}m')}m,nn'\right).$$
We need to show that
$$
\lambda_{m^{n}m',nn'}=\begin{bmatrix} &\lambda_{m',n'}\\ \lambda_{m,n}& \end{bmatrix}\  \  \mathrm{and} \ \
\lambda_{^{\partial_1m}{(^{n}m')}m,nn'}=\begin{bmatrix} \lambda_{m,n}&\\ &\lambda_{m',n'}\end{bmatrix}.
$$
We show these equalities in section \ref{sect:App2} of Appendix. In $\mathbf{(Aut\overline{ \delta})}_{2}$, suppose that $(H,F)$ and $(K,G)$ are 1-homotopies as given in $\mathbf{Ch}^2_K$.
In $M\rtimes N$, the horizontal composition can be pictured as
$$
\begin{array}{c}
    \begin{tikzcd}[row sep=small,column sep=0.8cm]
                    &  \ar[dd,Rightarrow,"\scriptscriptstyle{(m,n)}"{description}] &    & \ar[dd,Rightarrow,"\scriptscriptstyle{(m',n')}"{description}]& \\
                 {*}\ar[rr,bend left=50,"\scriptscriptstyle n"] \ar[rr,bend right=50,"\scriptscriptstyle \partial_1mn"']& &{*}\ar[rr,bend left=50,"\scriptscriptstyle n'"] \ar[rr,bend right=50,"\scriptscriptstyle \partial_1m'n'"'] &  & {*} \\ \ & \ & \  & \
    \end{tikzcd}
\end{array}
{:=}
\begin{array}{c}
    \begin{tikzcd}[row sep=small,column sep=0.8cm]
                 & \ar[dd, Rightarrow, "\scriptscriptstyle{(m^{n}m',nn')}"{description}] & \\
                 {*}\ar[rr,bend left=50,"\scriptscriptstyle nn'"] \ar[rr,bend right=50,"\scriptscriptstyle \partial_1mn\partial_1m'n'"'] &  &{*}\\
                 & \
    \end{tikzcd}
\end{array}
$$
The image of this diagram under $\lambda$ in $\mathbf{Aut}\overline{\delta}$ can be represented  by the diagram;
$$
 \begin{tikzcd}
\mathrm{K}_3 \arrow[rrr,shift left=0.90ex,"{\scriptscriptstyle {{\lambda}^2_n}}"description,pos=.75]\arrow[rrr,shift left=-0.90ex,"{\scriptscriptstyle \lambda^2_{\partial_1mn}}"description,pos=.35]    \arrow[ddd,"{\overline{\tau}_3}"'] &&&\mathrm{K}_3\arrow[rrr,shift left=0.95ex,"{\scriptscriptstyle \lambda^2_{n'}}"description,pos=.70]\arrow[rrr,shift left=-0.95ex,"{\scriptscriptstyle{\lambda}^2_{\partial_1m'n'}}"description,pos=.35]    \arrow[ddd,"{\overline{\tau}_3}"]                &&& \mathrm{K}_3 \arrow[ddd,,"{\overline{\tau}_3}"]  \\ \\ \\
\mathrm{K}_2 \arrow[rrr,shift left=0.90ex,"{\scriptscriptstyle {{\lambda}^1_n}}"description,pos=.75]\arrow[rrr,shift left=-0.90ex,"{\scriptscriptstyle \lambda^1_{\partial_1mn}}"description,pos=.35]  \arrow[ddd,"{\overline{\tau}_2}"']\arrow[uuurrr,"{\scriptscriptstyle{H'_2}}"{sloped,above=-0.3ex,xshift=0.0em}]&&&\mathrm{K}_2 \arrow[rrr,shift left=0.95ex,"{\scriptscriptstyle \lambda^1_{n'}}"description,pos=.70]\arrow[rrr,shift left=-0.95ex,"{\scriptscriptstyle{\lambda}^1_{\partial_1m'n'}}"description,pos=.35]  \arrow[ddd,"{\overline{\tau}_2}"]\arrow[uuurrr,"{\scriptscriptstyle{K'_2}}"{sloped,above=-0.3ex,xshift=0.0em}]             &&&\mathrm{K}_2  \arrow[ddd,"{\overline{\tau}_2}"]\\ \\ \\
\mathrm{K}_1 \arrow[rrr,shift left=0.90ex,"{\scriptscriptstyle {{\lambda}^0_n}}"description,pos=.75]\arrow[rrr,shift left=-0.90ex,"{\scriptscriptstyle \lambda^0_{\partial_1mn}}"description,pos=.35] \arrow[uuurrr,"{\scriptscriptstyle {H'_1}}"{sloped,above=-0.3ex,xshift=0.0em}] &&&\mathrm{K}_1 \arrow[rrr,shift left=0.95ex,"{\scriptscriptstyle \lambda^{0}_{n'}}"description,pos=.70]\arrow[rrr,shift left=-0.95ex,"{\scriptscriptstyle{\lambda}^{0}_{\partial_{1}m'n'}}"description,pos=.35]\arrow[uuurrr,"{\scriptscriptstyle{K'_1}}"{sloped,above=-0.3ex,xshift=0.0em}]   &&&\mathrm{K}_1
\end{tikzcd}
$$
It must be that
$$\lambda_{ \begin{bsmallmatrix} \Gamma&\\ &\Gamma' \end{bsmallmatrix}}=\lambda_{((^{\partial_1mn}m')m,nn')}=((\lambda^1_{n'}H'_1+K'_1\lambda^0_{\partial_1mn}),
(\lambda^2_{n'}H'_2+K'_2\lambda^1_{\partial_1mn}),FG)$$
where $$F=(\lambda^0_n,\lambda^1_n,\lambda^2_n), \ \ G=(\lambda^0_{n'},\lambda^1_{n'},\lambda^2_{n'})$$
and
$$  F'=(\lambda^0_{\partial_1mn},\lambda^1_{\partial_1mn},\lambda^2_{\partial_1mn}), \ \ G'=(\lambda^0_{\partial_1m'n'},\lambda^1_{\partial_1m'n'},\lambda^2_{\partial_1m'n'})$$
and
$$H'_1=\lambda'_{m,n},\ \ H'_2=\lambda''_{m,n},\ \ K'_1=\lambda'_{m',n'},\ \ K'_2=\lambda''_{m',n'}.$$
In this diagram, for any $\mathbf{e}_{1,1,p}\in \mathrm{K}_1$, we obtain
\begin{align*}
(\lambda^1_{n'}H'_1+K'_1 \lambda^0_{\partial_1mn})(\mathbf{e}_{1,1,p})=&\lambda^1_{n'}H'_1(\mathbf{e}_{1,1,p})+K'_1\lambda^0_{\partial_1mn}(\mathbf{e}_{1,1,p})\\
=&\lambda^1_{n'}(\lambda'_{m,n}(\mathbf{e}_{1,1,p}))+\lambda^1_{m',n'}(\lambda^0_{\partial_1mn}(\mathbf{e}_{1,1,p}))\\
=&\lambda^1_{n'}(\overline{\mathbf{v}^{\scriptscriptstyle11}}_{1,^{p}m,pn})+\lambda^1_{m',n'}(\mathbf{e}_{1,1,p\partial_1mn})\\
=&\overline{\mathbf{v}^{\scriptscriptstyle11}}_{1,^{p}m,pnn'}+\overline{\mathbf{v}^{\scriptscriptstyle11}}_{1,^{p\partial_1mn}(m'),p\partial_1mnn'}\\
=&\overline{\mathbf{v}^{\scriptscriptstyle11}}_{1,^{p\partial_1mn}(m')^pm,pnn'} \ (\because \text{relation (\ref{3})})\\
=&\lambda'_{^{\partial_1m}(^{n}m')m,nn'}(\mathbf{e}_{1,1,p}).
\end{align*}
Thus, we have $$\lambda'_{ \begin{bsmallmatrix} \Gamma&\\ &\Gamma' \end{bsmallmatrix}}=(\lambda^1_{n'}H'_1+K'_1 \lambda^0_{\partial_1mn}).$$
In section  \ref{sect:App2} of Appendix, we show the equality $\lambda^2_{n'} H'_2+K'_2\lambda^1_{\partial_1mn}=\lambda''_{ \begin{bsmallmatrix} \Gamma&\\ &\Gamma' \end{bsmallmatrix}}.$

Similarly, it must be that
$$\lambda_{\begin{bsmallmatrix} &\Gamma'\\ \Gamma& \end{bsmallmatrix}}= \lambda_{m^nm',nn'}=((K'_1 \lambda^0_n+\lambda'_{\partial_1m'n'}H'_1),(K'_2\lambda^1_n+\lambda^2_{\partial_1m'n'}H'_2),FG).$$
For any $\mathbf{e}_{1,1,p}\in \mathrm{K}_1$, we obtain
\begin{align*}
(K'_1\lambda^0_n+\lambda'_{\partial_1m'n'}H'_1)(\mathbf{e}_{1,1,p})=&K'_1(\lambda^0_n(\mathbf{e}_{1,1,p}))+\lambda^1_{\partial_1m'n'}(H'_1(\mathbf{e}_{1,1,p})) \\
=&\lambda^1_{m',n'}(\mathbf{e}_{1,1,pn})+\lambda^1_{\partial_1m'n'}(\lambda^1_{m,n}(\mathbf{e}_{1,1,p})) \\
=&\overline{\mathbf{v}^{\scriptscriptstyle11}}_{1,^{pn}m',pnn'}+\lambda^1_{\partial_1m'n'}(\overline{\mathbf{v}^{\scriptscriptstyle11}}_{1,^{p}m,pn})\\
=&\overline{\mathbf{v}^{\scriptscriptstyle11}}_{1,^{pn}m',pnn'}+\overline{\mathbf{v}^{\scriptscriptstyle11}}_{1,^{p}m,pn\partial_1m'n'}\\
=&\overline{\mathbf{v}^{\scriptscriptstyle11}}_{1,^p{m}^{pn}m',pnn'} \ (\because \text{relation (\ref{3})}) \\
=&\overline{\mathbf{v}^{\scriptscriptstyle11}}_{1,^p{(m^{n}m')},pnn'}\\
=&\lambda'_{m^nm',nn'}(\mathbf{e}_{1,1,p}).
\end{align*}
Thus, we have $$\lambda'_{\begin{bsmallmatrix} &\Gamma'\\ \Gamma& \end{bsmallmatrix}}=(K'_1\lambda^0_n+\lambda'_{\partial_1m'n'}H'_1).$$
 We will show the equality $K'_2\lambda^1_n+\lambda^2_{\partial_1m'n'}H'_2=\lambda''_{m^nm',nn'}=\lambda''_{\begin{bsmallmatrix} &\Gamma'\\ \Gamma& \end{bsmallmatrix}}$  in section \ref{sect:App2} in Appendix.

Thus, we have
$$\lambda_{\begin{bsmallmatrix} &\Gamma'\\ \Gamma& \end{bsmallmatrix}}={\begin{bmatrix} &{\lambda_{\Gamma'}}\\ {\lambda_{\Gamma}}& \end{bmatrix}}
\ \text{ and }\lambda_{ \begin{bsmallmatrix} \Gamma&\\ &\Gamma' \end{bsmallmatrix}}={ \begin{bmatrix}{\lambda_{\Gamma}}&\\ &{\lambda_{\Gamma'}}\end{bmatrix}}.$$


\subsection{$\lambda$ over 3-cells}
In $\mathfrak{C}^2$, for any 3-cell $(l,m,n)$, we will show that $\lambda_{l,m,n}$ is a 2-homotopy between 1-homotopies: $\lambda_{m,n}$ and $\lambda_{\partial_2 lm,n}$. This 3-cell in $\mathbf{(Aut \overline{ \delta})}_{3}$ can be given by
$$\lambda_{l,m,n}:=(\alpha'_{l,m,n},(\lambda'_{m,n},\lambda''_{m,n}),(\lambda^0_{n},\lambda^1_{n},\lambda^2_{n})).$$
The 2-source of $\lambda_{l,m,n}$ is $\lambda_{m,n}:=((\lambda'_{m,n},\lambda''_{m,n}),\lambda_{n})$ is a 1-homotopy from $\lambda_{n}$ to $\lambda_{\partial_1mn}$ and 2-target of $\lambda_{l,m,n}$ is a 1-homotopy
$$\lambda_{\partial_2lm,n}:=((\lambda'_{\partial_2lm,n},\lambda''_{\partial_2lm,n}),(\lambda^0_{n},\lambda^1_{n},\lambda^2_{n}))$$
from  $\lambda_{n}$ to $\lambda_{\partial_1mn}$. Thus, $\lambda_{l,m,n}$ can be pictured as
$$
   \begin{tikzcd}[row sep=0.3cm,column sep=scriptsize]
                & \ar[dd, Rightarrow, "{\scriptscriptstyle \lambda_{(m,n)}}"{swap,name=f,description}]& &\ar[dd, Rightarrow, "{\scriptscriptstyle \lambda_{(\partial_2lm,n)}}"'{swap,name=g,description}]& \\
               {\overline{\delta}} \ar[rrrr,bend left=40,"{\scriptscriptstyle \lambda_{n}}"] \ar[rrrr,bend right=40,"{\scriptscriptstyle \lambda_{\partial_1mn}}"']& \tarrow["\lambda_{\scriptscriptstyle {(l,m,n)}}" ,from=f,to=g, shorten >= -1pt,shorten <= 1pt ]{rrr}&  & & {\overline{\delta}}. \\
                \ & \ & \ & \ & \
    \end{tikzcd}
$$
Using the chain complex of length-2
$
\begin{tikzcd}
             \overline{\delta}:=\mathrm{K}_3\ar[r,"\overline{\tau}_3"]&\mathrm{K}_2\ar[r,"\overline{\tau}_2"]& \mathrm{K}_1
\end{tikzcd}
$
constructed in section \ref{complex}, the last picture can be given more clearly from our black box as follows:
\begin{equation*}
\xymatrix{\mathrm{K}_3 \ar[rr]^-{\overline{\tau}_3}\ar@<0.5ex>[dd]^<<<<<<{\lambda^2_{\partial_1mn}}\ar@<-0.5ex>[dd]_-{{\lambda}^2_n} &&\mathrm{K}_2 \ar[rr]^-{\overline{\tau}_2}\ar@<0.5ex>[dd]^<<<{\lambda^1_{\partial_1mn}}\ar@<-0.5ex>[dd]_-{{\lambda}^1_n}\ar@<-0.6ex>[ddll]|<<<<<<<<<<<<<<<<{\scriptscriptstyle{H'_2}} \ar@<0.6ex>[ddll]|<<<<<<<<{\scriptscriptstyle{K'_2}}&&\mathrm{K}_1\ar@{-->}@/^{0.6pc}/[ddllll]_(.4){\ \ {\large\alpha'}} \ar@<0.5ex>[dd]^<<<<<<{\lambda^0_{\partial_1mn}}\ar@<-0.5ex>[dd]_<<<<<<<<<<<{{\lambda}^0_n}\ar@<-0.6ex>[ddll]|<<<<<<<<<<<<<<<<{\scriptscriptstyle{H'_1}} \ar@<0.6ex>[ddll]|<<<<<<<<{\scriptscriptstyle{K'_1}}\ar@<0.5ex>[dd]
\\
\\
\mathrm{K}_3\ar[rr]_{\overline{\tau}_3} && \mathrm{K}_2 \ar[rr]_{\overline{\tau}_2} \ar[rr] &&\mathrm{K}_1}
\end{equation*}
where $H'_1=\lambda'_{m,n},\ \ K'_1=\lambda'_{\partial_2lm,n}$,\ \ $H'_2=\lambda''_{m,n},\ \ K'_2=\lambda''_{\partial_2lm,n}$. Using the action of $(l,m,n)$ in $\mathfrak{C}^2$ on a 3-cell $\mathbf{e}_{1,1,p}$ in $\mathrm{K}_1$, the homotopy component map $\alpha'=\alpha'_{(l,m,n)}:\mathrm{K}_1\longrightarrow \mathrm{K}_3$ can be given by
$$\alpha'_{(l,m,n)}(\mathbf{e}_{1,1,p})=\overline{\mathbf{v}^{\scriptscriptstyle22}}_{(^{p}l,^{p}m,pn)}$$
on generators.  For any element $(\mathbf{e}_{1,1,p})\in \mathrm{K}_1$; we have
$$\overline{\tau}_3\alpha'_{l,m,n}(\mathbf{e}_{1,1,p})=\overline\tau_3\overline{\mathbf{v}^{\scriptscriptstyle22}}_{(^{p}l,^{p}m,pn)}
=\overline{\mathbf{v}^{\scriptscriptstyle11}}_{1,^{p}({\partial_2} lm),pn}-\overline{\mathbf{v}^{\scriptscriptstyle11}}_{1,^{p}m,pn}$$
and
\begin{align*}
(K'_1-H'_1)(\mathbf{e}_{1,1,p})=&K'_1(\mathbf{e}_{1,1,p})-H'_1(\mathbf{e}_{1,1,p})\\
=&\lambda'_{\partial_2lm,n}(\mathbf{e}_{1,1,p})-\lambda'_{m,n}(\mathbf{e}_{1,1,p})\\
=&\overline{\mathbf{v}^{\scriptscriptstyle11}}_{1,^{p}({\partial_2 lm}),pn}-\overline{\mathbf{v}^{\scriptscriptstyle11}}_{1,^{p}m,pn}.
\end{align*}
Therefore, we obtain that $K'_1=H'_1+\overline\tau_3 \alpha'$, that is, $\lambda'_{(\partial_2lm,n)}=\lambda'_{(m,n)}+\overline\tau_3 \alpha'$. This is the chain homotopy condition for  2-homotopy $\alpha'$. In section \ref{sect:App4} of Appendix, the second homotopy condition for $\alpha'$;
$$
K'_2=H'_2+\alpha'\overline\tau_2
$$
will be showed.

\subsubsection{$\lambda_{l,m,n}$ preserves the 2-vertical composition of 3-cells}
The 2-vertical composition in $\mathfrak{C}^2$ of 3-cells $(l,m,n)$ and $(l',\partial_2lm,n)$ is given by
$$(l',\partial_2lm,n)\#_3(l,m,n)=(l'l,m,n)\in C_3.$$
We will show that
$$\lambda_{(l',\partial_2lm,n)\#_3(l,m,n)}=\lambda_{l'l,m,n}=\lambda_{l',\partial_2lm,n}\#_3\lambda_{l,m,n}\in (\mathbf{Aut}\overline{\delta})_3.$$
We have
$$\lambda_{l'l,m,n}=(\alpha'_{l'l,m,n},(\lambda'_{m,n},\lambda''_{m,n}),(\lambda^0_{n},\lambda^1_{n},\lambda^2_{n}))$$
and
$$\lambda_{l',\partial_2lm,n}\#_3\lambda_{l,m,n}=(\alpha'_{l',\partial_2lm,n},(\lambda'_{\partial_2lm,n},\lambda''_{\partial_2lm,n}),(\lambda^0_{n},\lambda^1_{n},\lambda^2_{n}))
\#_3(\alpha'_{l,m,n},(\lambda'_{m,n},\lambda''_{m,n}),(\lambda^0_{n},\lambda^1_{n},\lambda^2_{n})).$$

We know from the vertical composition of 3-cell in $\mathbf{(Aut \overline{ \delta})}_{3}$;
$$\lambda_{l',\partial_2lm,n}\#_3\lambda_{l,m,n}=(\alpha'_{l',\partial_2lm,n}+\alpha'_{l,m,n},(\lambda'_{m,n},\lambda''_{m,n}),(\lambda^0_{n},\lambda^1_{n},\lambda^2_{n}))$$
In this equality, we must show that
$$\alpha'_{l'l,m,n}=\alpha'_{l',\partial_2lm,n}+\alpha'_{l,m,n}.$$
For any element $(\mathbf{e}_{1,1,p})\in K_1$ we have;
$$\alpha'_{l'l,m,n}(\mathbf{e}_{1,1,p})=\overline{\mathbf{v}^{\scriptscriptstyle22}}_{^{p}{(l'l)},^{p}m,pn}$$
and
$$\alpha'_{l',\partial_2lm,n}(\mathbf{e}_{1,1,p})+\alpha'_{l,m,n}(\mathbf{e}_{1,1,p})=\overline{\mathbf{v}^{\scriptscriptstyle22}}_{(^{p}l',^{p}{(\partial_2lm)},pn)}
+\overline{\mathbf{v}^{\scriptscriptstyle22}}_{(^{p}l,^{p}m,pn)}.$$
From relation (\ref{1}), we have
$$\overline{\mathbf{v}^{\scriptscriptstyle22}}_{l'l,m,n}=\overline{\mathbf{v}^{\scriptscriptstyle22}}_{l',\partial_2lm,n}+\overline{\mathbf{v}^{\scriptscriptstyle22}}_{l,m,n}$$
and thus we obtain that
\begin{align*}
(\alpha'_{l',\partial_2lm,n}+\alpha'_{l,m,n})(\mathbf{e}_{1,1,p})=&\overline{\mathbf{v}^{\scriptscriptstyle22}}_{^{p}l',^{p}{(\partial_2lm)},pn}+\overline{\mathbf{v}^{\scriptscriptstyle22}}_{(^{p}l,^{p}m,pn)}\\
=&\overline{\mathbf{v}^{\scriptscriptstyle22}}_{^{p}{(l'l)},^{p}m,pn}\\
=&\alpha'_{l'l,m,n}(\mathbf{e}_{1,1,p}).
\end{align*}
Therefore, we obtain the equality:
$$\lambda_{(l',\partial_2l,mn)\#_3(l,m,n)}=\lambda_{l'l,m,n}=\lambda_{l',\partial_2lm,n}\#_3\lambda_{l,m,n}.$$
\subsubsection{$\lambda_{l,m,n}$ preserves the 1-vertical composition of 3-cells }
Recall that the 1-vertical composition in $\mathfrak{C}^2$ of 3-cells $(l,m,n)$ and $(l',m',\partial_1mn)$ is given by
$$(l',m',\partial_1mn)\#_1(l,m,n)=(l'^{m'}l,m'm,n)$$
as pictured in the following diagram
$$
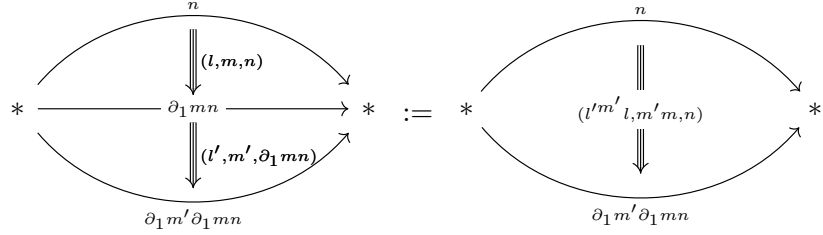
\begin{array}{c}
   \begin{tikzcd}[row sep=scriptsize,column sep=2.05cm]
    \vphantom{f}
    * \arrow[rr, "\scriptscriptstyle n",  ""{name=F, below}, bend left=50]
    \arrow[rr, "\scriptscriptstyle \partial_{1}m'\partial_{1}mn"',  ""{name=H, above}, bend right=50]
    \arrow[rr, "\scriptscriptstyle \partial_1mn" description, ""{name=GA,above}, ""{name=GB,below}]
     & &\vphantom{f}*
    \tarrow[from=F, to=GA, "{\scriptscriptstyle (l,m,n)}"]{d}
    \tarrow[from=GB, to=H, "{\scriptscriptstyle (l',m',\partial_1mn)}"]{d}
  \end{tikzcd}
\end{array}
{:=}
\begin{array}{c}
    \begin{tikzcd}[row sep=small,column sep=1.8cm]
                 &\tarrow["\scriptscriptstyle{(l'^{m'}l,m'm,n)}"description,shorten >= -8pt,shorten <= -8pt]{dd}& \\
                 {*} \ar[rr,bend left=50,"\scriptscriptstyle n"] \ar[rr,bend right=50,"\scriptscriptstyle\partial_{1}m'\partial_{1} mn"']& & {*}\\
                  & \
    \end{tikzcd}
\end{array}
$$
Under $\lambda$, these diagram can be pictured in $\mathbf{(Aut \overline{ \delta})}_{3}$ as
$$
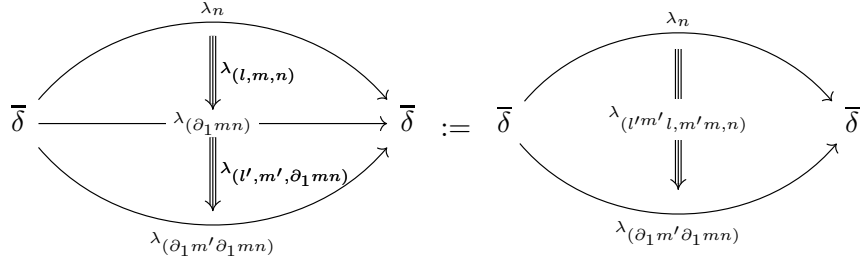
\begin{array}{c}
   \begin{tikzcd}[row sep=scriptsize,column sep=2.30cm]
    \vphantom{f}
    {\overline{\delta}}\arrow[rr, "\scriptscriptstyle \lambda_{n}",  ""{name=F, below}, bend left=50]
    \arrow[rr, "\scriptscriptstyle \lambda_{(\partial_{1}m'\partial_{1}mn)}"',  ""{name=H, above}, bend right=50]
    \arrow[rr, "\scriptscriptstyle \lambda_{(\partial_1mn)}" description, ""{name=GA,above}, ""{name=GB,below}]
     & &\vphantom{f} {\overline{\delta}}
    \tarrow[from=F, to=GA, "{\scriptscriptstyle \lambda_{(l,m,n)}}"]{d}
    \tarrow[from=GB, to=H, "{\scriptscriptstyle \lambda_{(l',m',\partial_1mn)}}"]{d}
  \end{tikzcd}
\end{array}
{:=}
\begin{array}{c}
    \begin{tikzcd}[row sep=small,column sep=1.8cm]
                 &\tarrow["\scriptscriptstyle{\lambda_{(l'^{m'}l,m'm,n)}}"description,shorten >= -8pt,shorten <= -8pt]{dd}& \\
                 {\overline{\delta}} \ar[rr,bend left=50,"\scriptscriptstyle \lambda_{n}"] \ar[rr,bend right=50,"\scriptscriptstyle\lambda_{(\partial_{1}m'\partial_{1} mn)}"']& &  {\overline{\delta}}\\
                  & \
    \end{tikzcd}
\end{array}
$$
We can write
$$\lambda_{l'^{m'}l,m'm,n}=(\alpha'_{l'^{m'}l,m'm,n},(\lambda'_{m'm,n},\lambda''_{m'm,n}),(\lambda^0_{n},\lambda^1_{n},\lambda^2_{n}))$$
and
$$\lambda_{l,m,n}=(\alpha'_{l,m,n},(\lambda'_{m,n},\lambda''_{m,n}),(\lambda^0_{n},\lambda^1_{n},\lambda^2_{n}))$$
and
$$\lambda_{l',m',\partial_1mn}=(\alpha'_{l',m',\partial_1mn},(\lambda'_{m',\partial_1mn},\lambda''_{m',\partial_1mn}),(\lambda^0_{\partial_1mn},\lambda^1_{\partial_1mn},\lambda^2_{\partial_1mn})).$$

Consider the diagram obtained from the black box:
 $$
 
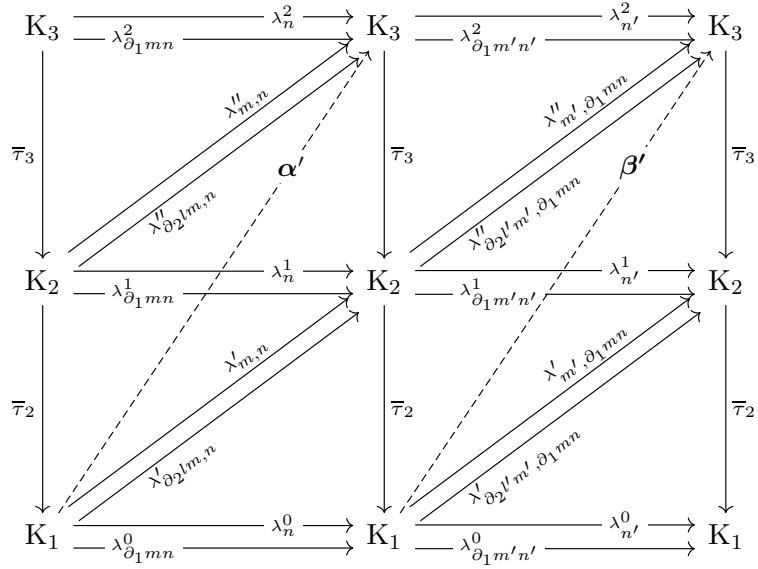
\begin{tikzcd}
\mathrm{K}_3 \arrow[rrrr,shift left=0.90ex,"{\scriptscriptstyle {{\lambda}^2_n}}"description,pos=.75]\arrow[rrrr,shift left=-0.90ex,"{\scriptscriptstyle \lambda^2_{\partial_1mn}}"description,pos=.25]    \arrow[dddd,"{\overline{\tau}_3}"'] &&&&\mathrm{K}_3\arrow[rrrr,shift left=0.95ex,"{\scriptscriptstyle \lambda^2_{n'}}"description,pos=.75]\arrow[rrrr,shift left=-0.95ex,"{\scriptscriptstyle{\lambda}^2_{\partial_1m'n'}}"description,pos=.30]    \arrow[dddd,"{\overline{\tau}_3}"]                &&&& \mathrm{K}_3 \arrow[dddd,,"{\overline{\tau}_3}"] \\ \\ \\ \\
\mathrm{K}_2 \arrow[rrrr,shift left=0.90ex,"{\scriptscriptstyle {{\lambda}^1_n}}"description,pos=.75]\arrow[rrrr,shift left=-0.90ex,"{\scriptscriptstyle \lambda^1_{\partial_1mn}}"description,pos=.25]  \arrow[dddd,"{\overline{\tau}_2}"']\arrow[uuuurrrr,shift left=-0.7ex,"\scriptscriptstyle {\lambda}''_{\partial_2lm,n}"{sloped,below=-0.3ex,xshift=-2.0em}]\arrow[uuuurrrr,shift left=0.7ex,"\scriptscriptstyle {\lambda}''_{m,n}"{sloped,above=-0.3ex,xshift=2.0em}]&&&&\mathrm{K}_2 \arrow[rrrr,shift left=0.95ex,"{\scriptscriptstyle \lambda^1_{n'}}"description,pos=.75]\arrow[rrrr,shift left=-0.95ex,"{\scriptscriptstyle{\lambda}^1_{\partial_1m'n'}}"description,pos=.30]  \arrow[dddd,"{\overline{\tau}_2}"]\arrow[uuuurrrr,shift left=-0.7ex,"{\scriptscriptstyle{\lambda}''_{\partial_2l'm',\partial_1mn}}"{sloped,below=-0.3ex,xshift=-2.0em}]\arrow[uuuurrrr,shift left=0.7ex,"{\scriptscriptstyle{\lambda}''_{m',\partial_1mn}}"{sloped,above=-0.3ex,xshift=2.0em}]              &&&&\mathrm{K}_2  \arrow[dddd,"{\overline{\tau}_2}"]\\ \\ \\ \\
\mathrm{K}_1 \arrow[rrrr,shift left=0.90ex,"{\scriptscriptstyle {{\lambda}^0_n}}"description,pos=.75]\arrow[rrrr,shift left=-0.90ex,"{\scriptscriptstyle \lambda^0_{\partial_1mn}}"description,pos=.25] \arrow[uuuurrrr,shift left=-0.7ex,"\scriptscriptstyle {\lambda}'_{\partial_2lm,n}"{sloped,below=-0.3ex,xshift=-2.0em}]\arrow[uuuurrrr,shift left=0.7ex,"\scriptscriptstyle {\lambda}'_{m,n}"{sloped,above=-0.3ex,xshift=2.0em}] \arrow[uuuuuuuurrrr,dashed,"\Large\boldsymbol\alpha'"description,pos=.75]         &&&&\mathrm{K}_1 \arrow[rrrr,shift left=0.95ex,"{\scriptscriptstyle \lambda^{0}_{n'}}"description,pos=.75]\arrow[rrrr,shift left=-0.95ex,"{\scriptscriptstyle{\lambda}^{0}_{\partial_{1}m'n'}}"description,pos=.30] \arrow[uuuurrrr,shift left=-0.7ex,"{\scriptscriptstyle{\lambda}'_{\partial_2l'm',\partial_1mn}}"{sloped,below=-0.3ex,xshift=-2.0em}]\arrow[uuuurrrr,shift left=0.7ex,"{\scriptscriptstyle{\lambda}'_{m',\partial_1mn}}"{sloped,above=-0.3ex,xshift=2.0em}]\arrow[uuuuuuuurrrr,dashed,"\Large\boldsymbol{\beta'}" description,pos=.75]          &&&&\mathrm{K}_1
\end{tikzcd}
$$
We can write
\begin{align*}
\lambda^2_{\partial_1m'n'}(\alpha'_{l,m,n})(\mathbf{e}_{1,1,p})+\beta'_{l',m',\partial_1mn}\lambda^0_{n}(\mathbf{e}_{1,1,p})=&\lambda^2_{\partial_1m'n'}(\overline{\mathbf{v}^{\scriptscriptstyle22}}_{^{p}l,^{p}m,pn})+\beta'_{l',m',\partial_1mn}(\mathbf{e}_{1,1,pn})\\
=&\overline{\mathbf{v}^{\scriptscriptstyle22}}_{^{p}l,^{p}m,pn\partial_1m'\partial_1mn}+\overline{\mathbf{v}^{\scriptscriptstyle22}}_{^{pn}l',^{pn}m',pn\partial_1mn}.
\end{align*}
Since $$\overline\tau_1\lambda'_{m,n}+\lambda^0_n=\lambda^0_{\partial_1mn},$$
$$\lambda''_{m,n} \overline\tau_2+\lambda^2_n=\lambda^2_{\partial_1mn}$$
and
$$\overline\tau_2\lambda''_{m,n}+\lambda'_{m,n} \overline\tau_2+\lambda'_{n}=\lambda'_{\partial_1mn},$$
we can say that
\begin{multline*}
((\lambda'_{m',\partial_1mn},\lambda''_{m',\partial_1mn}),(\lambda^0_{\partial_1mn},\lambda^1_{\partial_1mn},\lambda^2_{\partial_1mn}))\#_1((\lambda'_{m,n} ,\lambda''_{m,n}),(\lambda^0_{n},\lambda^1_{n},\lambda^2_{n}))\\
\begin{aligned}
=&((\lambda'_{m',\partial_1mn}+\lambda'_{m,n}\lambda''_{m',\partial_1mn}+\lambda''_{m,n}),(\lambda^0_{n},\lambda^1_{n},\lambda^2_{n}))\\
=&((\lambda'_{m'm,n},\lambda''_{m'm,n}),(\lambda^0_{n},\lambda^1_{n},\lambda^2_{n}))
\end{aligned}
\end{multline*}
Therefore, we obtain
$$\lambda_{l'^{m'}l,m'm,n}=(\alpha'_{l'^{m'}l,m'm,n},(\lambda'_{m'm,n},\lambda''_{m'm,n}),(\lambda^0_{n},\lambda^1_{n},\lambda^2_{n}))$$
and
$$\lambda'_{l',m',\partial_1mn}\#_1\lambda_{l,m,n}=(\alpha'_{l',m',\partial_1mn}+\alpha'_{l,m,n},(\lambda'_{m'm,n},
\lambda''_{m'm,n}),(\lambda^0_{n},\lambda^1_{n},\lambda^2_{n})).$$
Thus, to show this equality, we must prove that
$$\alpha'_{l'^{m'}l,m'm,n}=\alpha'_{l',m',\partial_1mn}+\alpha'_{l,m,n}.$$
For any element $(\mathbf{e}_{1,1,p})\in \mathrm{K}_1$, we obtain
$$\alpha'_{l'^{m'}l,m'm,n}(\mathbf{e}_{1,1,p})=\overline{\mathbf{v}^{\scriptscriptstyle22}}_{^{p}{(l')}^{p}{(^{m'}l)},^{p}{(m'm)},pn}$$
and
\begin{align*}
(\alpha'_{l',m',\partial_1mn}+\alpha'_{l,m,n})(\mathbf{e}_{1,1,p})=&\alpha'_{l',m',\partial_1mn}(\mathbf{e}_{1,1,p})+\alpha'_{l,m,n}(\mathbf{e}_{1,1,p})\\
=&\overline{\mathbf{v}^{\scriptscriptstyle22}}_{^{p}{(l')},^{p}{m'},p\partial_1mn}+\overline{\mathbf{v}^{\scriptscriptstyle22}}_{(^{p}l,^{p}m,pn)}.
\end{align*}
We know from relation (\ref{2});
$$\overline{\mathbf{v}^{\scriptscriptstyle22}}_{l'^{m'}l,m'm,n}=\overline{\mathbf{v}^{\scriptscriptstyle22}}_{l',m',\partial_1mn}+\overline{\mathbf{v}^{\scriptscriptstyle22}}_{l,m,n},$$
we obtain
$$\overline{\mathbf{v}^{\scriptscriptstyle22}}_{^{p}l,^{p}m,pn}+\overline{\mathbf{v}^{\scriptscriptstyle22}}_{^{p}{l'},^{p}{m'},p\partial_1mn}=\overline{\mathbf{v}^{\scriptscriptstyle22}}_{^{p}{l'}^{pm'}{(^{p}l)},^{p}{m'}^{p}{m},pn}=\overline{\mathbf{v}^{\scriptscriptstyle22}}_{^{p}{(l'^{m'}l)},^{p}{(m'm)},pn}$$
Thus, we can say that $\lambda_{l,m,n}$ preserves the 1-vertical composition of 3-cells from  $\mathfrak{C}^2$ to $\mathbf{Aut}\overline{\delta}$. Therefore, we have proven the following equality:
$$\lambda_{(l',m',\partial_1mn)\#_1(l,m,n)}=\lambda_{(l',m',\partial_1mn)}\#_1\lambda_{(l,m,n)}.$$

\subsubsection{$\lambda_{l,m,n}$ preserves the group operations}
For $J=(l,m,n)$ and $J^\prime=(l^\prime,m^\prime,n^\prime)$ in $L\rtimes M \rtimes N$, recall the semi-direct product of $J$ and $J'$ given by
$$J\cdot J'=\nabla=(l\{\partial_2(^{n}{l'}),m\}^{n}{(l')^{-1}},m^nm',nn')=(l^{m}{(^{n}{l'}),m^nm',nn'}).$$
This can be represented pictorially as
$$
\begin{array}{c}
\begin{tikzcd}[row sep=large]
 & {*} \ar[dr, Rightarrow,"{\scriptscriptstyle (\partial_2lm,n)}"{sloped,above=0.0ex,xshift=0.0em}]& &{*} \ar[dr, Rightarrow,"{\scriptscriptstyle (\partial_2l'm',n')}"{sloped,above=0.0ex,xshift=0.0em}]  \\
  {*}\ar[ur,"\scriptscriptstyle n"{sloped,above=0.0ex,xshift=-0.0em}] \ar[dr, Rightarrow,"{\scriptscriptstyle (m,n)}"{sloped,below=0.0ex,xshift=-0.0em}] & {J} & {*}\ar[ur,"\scriptscriptstyle n'"{sloped,above=0.0ex,xshift=-0.0em}] \ar[dr, Rightarrow,"{\scriptscriptstyle (m',n')}"{sloped,below=0.0ex,xshift=-0.0em}]& {J'} & {*} \\
  & {*}\ar[ur,"{\scriptscriptstyle \partial_1mn}"{sloped,below=0.0ex,xshift=-0.0em}]& & {*}\ar[ur,"{\scriptscriptstyle \partial_1m'n'}"{sloped,below=0.0ex,xshift=-0.0em}]
\end{tikzcd}
\end{array}
{:=}
\begin{array}{c}
\begin{tikzcd}[row sep=large]
 & {*} \ar[dr, Rightarrow,"{\scriptscriptstyle (\partial_2lm{^{n}({\partial_2l'm'}),nn'})}"] \\
  {*}\ar[ur,"\scriptscriptstyle nn'"{sloped,above=0.0ex,xshift=-0.0em}] \ar[dr, Rightarrow,"{\scriptscriptstyle (m^nm',nn')}"{sloped,below=0.0ex,xshift=-0.0em}] & {\nabla} & {*}\\
  & {*}\ar[ur,"{\scriptscriptstyle \partial_1mn\partial_1m'n'}"{sloped,below=0.0ex,xshift=-0.0em}]
\end{tikzcd}
\end{array}
$$
In $\mathbf{Ch}^2_K$, the product of 3-cells for $\mathbf{(Aut \overline{ \delta})}_{3}$ can be given by
\begin{multline*}
(\alpha',(H'_1,H'_2),(F_0,F_1,F_2))\cdot (\beta',(K'_1,K'_2),(G_0,G_1,G_2))\\
\begin{aligned}
=&(G_2\alpha'+\beta'F_0,(K'_1 F_0+G_1 H'_1,K'_2 F_1+G_2 H''_2),(F_0G_0,F_1G_1,F_2G_2)).
\end{aligned}
\end{multline*}
This is the group operation in $\mathbf{(Aut \overline{ \delta})}_{3}$.
We will show that
$$\lambda_{(l,m,n)\cdot (l',m',n')}=\lambda_{l(^{m_n}{l')},m^nm',nn'}=\lambda_{l,m,n}\cdot \lambda_{l',m',n'}.$$
In $\mathfrak{C}^2$, we can consider that $(l{^{m}(^{n}l')},m^nm',nn')$ is a 3-cell as
$$
    \xymatrix{J\cdot J'=\nabla=(l{^{m}(^{n}l')},m^nm',nn'):(1,m^nm',nn')\ar@3{->}[r]&(1,\partial_2lm{^{n}({\partial_2l'm'}),nn'})}
$$
For $J,J'$ in $\mathfrak{C}^2$ , we have;
$$\lambda_J=\lambda_{l,m,n}=(\alpha'_{l,m,n},(\lambda'_{m,n},\lambda''_{m,n}),(\lambda^0_{n},\lambda^1_{n},\lambda^2_{n}))$$
and
$$\lambda_{J'}=\lambda_{l',m',n'}=(\alpha'_{l',m',n'},(\lambda'_{m',n'},\lambda''_{m',n'}),(\lambda^0_{n'},\lambda^1_{n'},\lambda^2_{n'}))$$
and; we must show that
\begin{align*}
\lambda_{J\cdot  J'}=&\lambda_{(l{^{m}(^{n}l')},m^nm',nn')}\\
=&((\lambda^2_{\partial_1m'n'}\alpha'_{l,m,n}+\alpha'_{l',m',n'}\lambda^0_{n}),\\
&(\lambda'_{m',n'}\lambda^0_{n}+\lambda'_{\partial_1m'n'}\lambda'_{m,n},\lambda''_{m',n'}\lambda^1_{n}+\lambda^2_{\partial_1m'n'}\lambda''_{m,n}),\\
&(\lambda^0_{nn'},\lambda'_{nn'},\lambda^2_{nn'}))\\
=&(\alpha'_{l{^{m}(^{n}l')},m^nm',nn'},(\lambda'_{m^nm',nn'},\lambda''_{m^nm',nn'}),(\lambda^0_{nn'},\lambda^1_{nn'},\lambda^2_{nn'})).
\end{align*}
For any element $(\mathbf{e}_{1,1,p})\in \mathrm{K}_1$, we have
$$\alpha'_{l{^{m}(^{n}l')},m^nm',nn'}(\mathbf{e}_{1,1,p})=\overline{\mathbf{v}^{\scriptscriptstyle22}}_{(^{p}{(l(^{m_n}{l'})),^{p}{(m^nm')},pnn')}}$$
On the other hand;
\begin{align*}
(\lambda^2_{\partial_1m'n'}\alpha'_{l,m,n}+\alpha'_{l',m',n'}\lambda^0_{n})(\mathbf{e}_{1,1,p})=&\lambda^2_{\partial_1m'n'}\alpha'_{(l,m,n)}(\mathbf{e}_{1,1,p})+\alpha'_{(l',m',n')}\lambda^0_{n}(\mathbf{e}_{1,1,p})\\ =&\lambda^2_{\partial_1m'n'}(\overline{\mathbf{v}^{\scriptscriptstyle22}}_{^{p}l,^{p}m,pn})+\alpha'_{l',m',n'}(\mathbf{e}_{1,1,pn})\\
=&\overline{\mathbf{v}^{\scriptscriptstyle22}}_{^{p}l,^{p}m,pn\partial_1m'n'}+\overline{\mathbf{v}^{\scriptscriptstyle22}}_{^{p}{(^{n}l')},^{p}{(^{n}m')},pnn'}
\end{align*}
Using the relation (\ref{2}) in $J_2$ given by
$$\overline{\mathbf{v}^{\scriptscriptstyle22}}_{l'{^{m'}l},m'm,n}=\overline{\mathbf{v}^{\scriptscriptstyle22}}_{l,m,n}+\overline{\mathbf{v}^{\scriptscriptstyle22}}_{l',m',\partial_1mn}$$
and by taking in this equality; $m={^{p}{(^{n}m')},\ m'={^{p}m}}$, $n=pnn',\ \partial_1mn=pn\partial_1m'n',\ l={^{p}{(^{n}l')}}$ and $l'={^{p}l}$; we obtain
$$\overline{\mathbf{v}^{\scriptscriptstyle22}}_{^{p}l,^{p}m,pn\partial_1m'n'}+\overline{\mathbf{v}^{\scriptscriptstyle22}}_{^{p}{(^{n}l')},^{p}{(^{n}m')},pnn'}=\overline{\mathbf{v}^{\scriptscriptstyle22}}_{^{p}l{^{p}{(^{m_n}{l'})},^{p}m{^{p}{(^{n}m')}},pnn'}}$$
Thus, we have showed that $\lambda$ preserves the group operation from $C_3$ to $(\mathbf{Aut}\overline{\delta})_3$.

\section{Cayley's Theorem for 2-Crossed Modules and  Cat$^2$-Groups:}
The construction given above can be summarized in the following definition.

\begin{defn} Consider a 2-crossed module of groups:
$$\xymatrix{\mathfrak{X}:=(L\ar[r]^-{\partial_2}&M\ar[r]^-{\partial_1}&N,\{-,-\})}$$
and its associated Gray 3-group-groupoid with a single object or cat$^2$-group:
$$
\xymatrix{\mathfrak{C}^2:=L\rtimes M \rtimes N \ar@<1ex>[r]^-{s_3,t_3} \ar@<0ex>[r]&1\rtimes M \rtimes N \ar@<1ex>[r]^-{s_2,t_2} \ar@<0ex>[r]\ar@<1ex>[l]^-{e_3}& 1\rtimes 1 \rtimes N\ar@<1ex>[l]^-{e_2}\ar@<0.5ex>[r]\ar@<-0.5ex>[r]& \{*\}}.
$$
The right regular representation of $\mathfrak{C}^2$ is a lax 3-functor (contravariant on 1-cells)
$$\mathbf{\lambda}:\mathfrak{C}^2\longrightarrow\mathbf{Ch}^2_K$$
searching each $n\in N$ to the chain automorphism $\lambda_n=(\lambda^2_n,\lambda^1_n,\lambda^0_n)$ where
$$\lambda^0_n(\mathbf{e}_{1,1,n'})=\mathbf{e}_{1,1,nn'},\ \ \lambda^1_n(\overline{\mathbf{v}^{\scriptscriptstyle11}}_{1,m',n'})=\overline{\mathbf{v}^{\scriptscriptstyle11}}_{1,m',n'n}, \ \ \lambda^2_n(\overline{\mathbf{v}^{\scriptscriptstyle22}}_{l',m',n'})=\overline{\mathbf{v}^{\scriptscriptstyle22}}_{l',m',n'n}$$
and each $(1,m,n)\in 1\rtimes M\rtimes N$ to the 1-homotopy $\lambda_{m,n}:\lambda_n\Rightarrow \lambda_{\partial_1mn}$ with the chain homotopy components;
$$\lambda'_{m,n}(\mathbf{e}_{1,1,n'})=\overline{\mathbf{v}^{\scriptscriptstyle11}}_{1,^{n'}m,n'n},\ \ \lambda''_{m,n}(\overline{\mathbf{v}^{\scriptscriptstyle11}}_{1,m',n'})=\overline{\mathbf{v}^{\scriptscriptstyle22}}_{\{m',^{n'}m\},m'{^{n'}m},n'n}$$
where all chain automorphisms and homotopies reside in  $\mathbf{(Aut \overline{ \delta})}_{1}$ and  $\mathbf{(Aut \overline{ \delta})}_{2}$ for the linear transformations;
$\delta_2:=\overline{\tau_2}|_{\ke \overline {\sigma_2}}$ and
$\delta_3:=\overline{\tau_3}|_{\ke \overline {\sigma_3}}$ obtained from $\overline{K(\mathfrak{C}^2)}$ of $\mathfrak{C}^2$ and each $(l,m,n)\in L\rtimes M\rtimes N$ to the 2-homotopy $\lambda_{l,m,n}:\lambda_{m,n}\Rrightarrow \lambda_{\partial_2lm,n}$ with the chain homotopy component
$$\alpha'_{l,m,n}(\mathbf{e}_{1,1,n'})=\overline{\mathbf{v}^{\scriptscriptstyle22}}_{^{n'}l,^{n'}m,n'n}.$$
Since the construction may be applied to any Gray-3-group-groupoid (or any 2-crossed module), it gives us a Gray-3-group-groupoid with a single object or cat$^2$-group version of Cayley's theorem, in term of linear regular representation.
\end{defn}
Therefore, we can give the following result.
\begin{thm}$\boldsymbol{\mathrm{(Cayley)}}$
 For any 2-crossed module $\mathfrak{X}$ and its associated cat$^2$-group $\mathfrak{C}^2$, the right regular representation as given above, exists.
 \end{thm}
 It was shown that 2-crossed modules, Gray 3-groupoids with a single 0-cell and cat$^2$-groups are equivalent. We have a definition of regular representations for a cat$^2$-group $\mathfrak{C}^2$, and this may also be considered as a regular representation of corresponding 2-crossed module. Therefore, we might give a regular representation of the 2-crossed module $\mathfrak{X}$ to be a regular representation of $\mathfrak{C}^2(\mathfrak{X})$.

\section{Appendix}

\subsubsection{The proof of interchange law}\label{interchangelaw}

For any 3-cells;
$\alpha=(l_1,m_1,n_1)$, \ $\beta=(l'_1,\partial_2l_1m_1,n_1)$,\ $\gamma=(l_2,m_2,n_2)$ and $\delta=(l'_2,\partial_2l_2m_2,n_2)$ in $ L\rtimes M\rtimes N$;\\
We must show that
$$(\alpha \#_3 \beta)\cdot(\gamma \#_3 \delta)=(\alpha \cdot \gamma) \#_3(\beta \cdot \delta)$$ for the vertical composition $\#_3$ and semi-direct product of 3-cells in $\mathfrak{C}^2$. Therefore, we obtain;
\begin{align*}
(\alpha \#_3 \beta)\cdot(\gamma \#_3 \delta)=&((l_1,m_1,n_1)\#_3(l'_1,\partial_2l_1m_1,n_1))\cdot((l_2,m_2,n_2)\#_3(l'_2,\partial_2l_2m_2,n_2))\\
=&(l'_1l_1,m_1,n_1)\cdot(l'_2l_2,m_2,n_2)\\
=&(l'_1l_1(^{m_1}{(^{n_1}{(l'_2l_2)})),m_1^{n_1}m_2,n_1n_2})\\
=&(l'_1l_1(^{m_1}(^{n_1}{l'_2})){l^{-1}_1}l_1(^{m_1}(^{n_1}{l_2})),m_1^{n_1}m_2,n_1n_2)\\
=&(l'_1(^{\partial_2l_1m_1}(^{n_1}{l'_2}))l_1(^{m_1}(^{n_1}{l_2})),m_1^{n_1}m_2,n_1n_2)\\
=&(l_1(^{m_1}{(^{n_1}{l_2})),m_1^{n_1}m_2,n_1n_2}) \#_3(l'_1(^{\partial_2l_1m_1}{(^{n_1}{l'_2})),\partial_2l_1m_1(^{n_1}(\partial_2l_2m_2)),n_1n_2})\\
=&((l_1,m_1,n_1)\cdot(l_2,m_2,n_2))\#_3((l'_1,\partial_2l_1m_1,n_1)\cdot(l'_2,\partial_2l_2m_2,n_2))\\
=&(\alpha \cdot \gamma) \#_3(\beta \cdot \delta)
\end{align*}
\subsubsection{The proof of equality $H'_2\overline{\tau}_3=\lambda^2_{\partial_1mn}-\lambda^2_n.$} \label{sect:App5}

We need to add a new generator element for $J_2$ as
\begin{multline*}
{\mathbf{u}_{\scriptscriptstyle5}}=\mathbf{e}_{^{\partial_1m}{(^{n}{l'})},^{\partial_1m}{(^{n}{m'})},nn'}-\mathbf{e}_{1,^{\partial_1mn}{m'},nn'}
-\mathbf{e}_{^{m}{(^{n}{l'})},m^nm',nn'}+\mathbf{e}_{1,m^nm',nn'}\\
-\mathbf{e}_{l',m',n'\partial_1mn}+\mathbf{e}_{1,m',n'\partial_1mn}+\mathbf{e}_{l',m',n'n}-\mathbf{e}_{1,m',n'n}\in J_2
\end{multline*}
and
\begin{multline*}
{\mathbf{v}_{\scriptscriptstyle5}}=\mathbf{e}_{1,^{\partial_1mn}{(\partial_2l'm')},nn'}-\mathbf{e}_{1,^{\partial_1mn}{m'},nn'}-\mathbf{e}_{1,m{^{n}{(\partial_2l'm')}},nn'}+\mathbf{e}_{1,m^nm',nn'}\\ -\mathbf{e}_{1,\partial_2l'm',n'\partial_1mn}+\mathbf{e}_{1,m',n'\partial_1mn}+\mathbf{e}_{1,\partial_2l'm',n'n}-\mathbf{e}_{1,m',n'n}\in J_1.
\end{multline*}
Using the generator element ${\mathbf{u}_{\scriptscriptstyle5}}$ in $J_2$ we have the following relation for $\mathrm{K}_3$;\\
$\overline{\mathbf{v}^{\scriptscriptstyle22}}_{^{\partial_1m}{(^{n}{l'})},^{\partial_1m}{(^{n}{m'})},nn'}-\overline{\mathbf{v}^{\scriptscriptstyle22}}_{^{m}{(^{n}{l'})},m^nm',nn'}
-\overline{\mathbf{v}^{\scriptscriptstyle22}}_{l',m',n'\partial_1mn}+\overline{\mathbf{v}^{\scriptscriptstyle22}}_{l',m',n'n}=\overline{0}$ in $\mathrm{K}_3$.\\
Using the last equality we obtain that;
\begin{align*}
(H'_2\overline{\tau}_3-(\lambda^2_{\partial_1mn}-\lambda^2_n))(\overline{\mathbf{v}^{\scriptscriptstyle22}}_{l',m',n'})
=&\overline{\mathbf{v}^{\scriptscriptstyle22}}_{\{m,^{n}{(\partial_2l'm')\}},m{^{n}{(\partial_2l'm')}},nn'}-\overline{\mathbf{v}^{\scriptscriptstyle22}}_{\{m,^{n}{m'}\},m^nm',nn'}\\
&-\overline{\mathbf{v}^{\scriptscriptstyle22}}_{l',m',n'\partial_1mn}+\overline{\mathbf{v}^{\scriptscriptstyle22}}_{l',m',n'n}
\end{align*}
In relation $(B)$, if we take $l=1$; we have;
\begin{align*}
(H'_2\overline{\tau}_3-(\lambda^2_{\partial_1mn}-\lambda^2_n))(\overline{\mathbf{v}^{\scriptscriptstyle22}}_{l',m',n'})
=&\overline{\mathbf{v}^{\scriptscriptstyle22}}_{^{\partial_1m}{(^{n}{l'})},^{\partial_1mn}{m'},nn'}-\overline{\mathbf{v}^{\scriptscriptstyle22}}_{^{m}{(^{n}{l'})},m^nm',nn'}\\
&-\overline{\mathbf{v}^{\scriptscriptstyle22}}_{l',m',n'\partial_1mn}+\overline{\mathbf{v}^{\scriptscriptstyle22}}_{l',m',n'n}=\overline{0}
\end{align*}
Thus, we obtain;
$$H'_2\overline{\tau}_3=\lambda^2_{\partial_1mn}-\lambda^2_n.$$

\subsubsection{The proof of  equality $\lambda''_{(m',\partial_1mn)\#_2(m,n)}=K'_2+H'_2$} \label{sect:App1}

For any  $\overline{\mathbf{v}^{\scriptscriptstyle11}}_{1,d,p}\in \mathrm{K}_2$,
$$\lambda''_{(m',\partial_1mn)\#_2(m,n)}(\overline{\mathbf{v}^{\scriptscriptstyle11}}_{1,d,p})=\lambda''_{m'm,n}(\overline{\mathbf{v}^{\scriptscriptstyle11}}_{1,d,p})
=\overline{\mathbf{v}^{\scriptscriptstyle22}}_{\{d,^{p}(m'm)\},d^{p}(m'm),pn}$$
and
$$(K'_2+H'_2)(\overline{\mathbf{v}^{\scriptscriptstyle11}}_{1,d,p})=(\lambda''_{m',\partial_1mn}+\lambda''_{m,n})(\overline{\mathbf{v}^{\scriptscriptstyle11}}_{1,d,p})
=\overline{\mathbf{v}^{\scriptscriptstyle22}}_{\{d,^{p}m'\},d^{p}m',p\partial_1mn}+\overline{\mathbf{v}^{\scriptscriptstyle22}}_{\{d,^{p}m\},d^{p}m,pn}.$$
We need to add a new relations in $J_2$ and $J_1$ as follows:
\begin{multline*}
\mathbf{e}_{\{m,m'm''\},mm'm'',n''}-\mathbf{e}_{1,mm'm'',n''}-\mathbf{e}_{\{m,m'\},mm',\partial_1m''n''}\\
\begin{aligned}
&+\mathbf{e}_{1,mm',\partial_1m''n''}-\mathbf{e}_{\{m,m''\},mm'',n''}+\mathbf{e}_{1,mm',n''}\in J_2
\end{aligned}
\end{multline*}
and;
$$\mathbf{e}_{1,mm'm'',n''}-\mathbf{e}_{1,mm',\partial_1m''n''}-\mathbf{e}_{1,mm',n''}\in J_1.$$
Thus we obtain
$$\overline{\mathbf{v}^{\scriptscriptstyle22}}_{\{m,m'm''\},mm'm'',n''}=\overline{\mathbf{v}^{\scriptscriptstyle22}}_{\{m,m'\},mm',\partial_1m''n''}+\overline{\mathbf{v}^{\scriptscriptstyle22}}_{\{m,m''\},mm'',n''}.$$

If we take; $(1,m,n)=(1,d,p), \ (1,m',n')=(1,^{p}m',p\partial_1mn)$ and $(1,m'',n'')=(1,^pm,pn)$, then, $\partial_1m''n''=p\partial_1mn$.

We obtain
$$\overline{\mathbf{v}^{\scriptscriptstyle22}}_{\{d,^{p}m'^{p}m\},d^{p}m'^{p}m,pn}=\overline{\mathbf{v}^{\scriptscriptstyle22}}_{\{d,^{p}m'\},d^{p}m',p\partial_1mn}+\overline{\mathbf{v}^{\scriptscriptstyle22}}_{\{d,^{p}m\},d^{p}m,pn}$$
Therefore;
$$\lambda''_{(m',\partial_1mn)\#_2(m,n)}(\overline{\mathbf{v}^{\scriptscriptstyle11}}_{1,d,p})=(K'_2+H'_2)(\overline{\mathbf{v}^{\scriptscriptstyle11}}_{1,d,p}).$$
In general;
\begin{align*}
\lambda_{(m',\partial_1mn)\#_2(m,n)}=&\left((K'_1+H'_1,K'_2+H'_2),F\right)\\
=&\left((\lambda'_{m',\partial_1mn}+\lambda'_{m,n}, \ \lambda''_{m',\partial_1mn}+\lambda''_{m,n}),\lambda_n\right)\\
=&(K,G)\#_2(H,F).
\end{align*}

\subsubsection{The proof of equality: $K'_2\lambda^1_n+\lambda^2_{\partial_1m'n'}H'_2=\lambda''_{m^nm',nn'}$}\label{sect:App2}


We will show that
$$\lambda''_{\begin{bsmallmatrix} &\Gamma'\\ \Gamma& \end{bsmallmatrix}}={\begin{bmatrix} &{\lambda''_{\Gamma'}}\\ {\lambda''_{\Gamma}}& \end{bmatrix}}
\text{  and  }
\lambda''_{ \begin{bsmallmatrix} \Gamma&\\ &\Gamma' \end{bsmallmatrix}}={ \begin{bmatrix}{\lambda''_{\Gamma}}&\\ &{\lambda''_{\Gamma'}}\end{bmatrix}}.$$

For $\overline{\mathbf{v}^{\scriptscriptstyle11}}_{1,d,p}\in \mathrm{K}_2$, we obtain;

$$\lambda''_{\begin{bsmallmatrix} &\Gamma'\\ \Gamma& \end{bsmallmatrix}}(\overline{\mathbf{v}^{\scriptscriptstyle11}}_{1,d,p})=\lambda''_{m^nm',nn'}(\overline{\mathbf{v}^{\scriptscriptstyle11}}_{1,d,p})=\overline{\mathbf{v}^{\scriptscriptstyle22}}_{{\{d,^{p}(m^{n}m')\},d^{p}(m^{n}m'),pnn'}}$$
and
\begin{align*}
(K'_2 \lambda^1_{n}+\lambda^2_{\partial_1m'n'} H'_2)(\overline{\mathbf{v}^{\scriptscriptstyle11}}_{1,d,p})=&K'_2 \lambda^1_{n}(\overline{\mathbf{v}^{\scriptscriptstyle11}}_{1,d,p})+\lambda^2_{\partial_1m'n'} H'_2(\overline{\mathbf{v}^{\scriptscriptstyle11}}_{1,d,p}) \\
=&\lambda''_{m',n'}(\lambda^1_n(\overline{\mathbf{v}^{\scriptscriptstyle11}}_{1,d,p}))+\lambda^2_{\partial_1m'n'}(\lambda''_{m,n}(\overline{\mathbf{v}^{\scriptscriptstyle11}}_{1,d,p}))\\
=&\lambda''_{(m',n')}(\overline{\mathbf{v}^{\scriptscriptstyle11}}_{1,d,pn})+\lambda^2_{\partial_1m'n'}(\overline{\mathbf{v}^{\scriptscriptstyle22}}_{\{d,^{p}m\},d^{p}m,pn})\\
=&\overline{\mathbf{v}^{\scriptscriptstyle22}}_{\{d,^{pn}m'\},d^{pn}m',pnn'}+\overline{\mathbf{v}^{\scriptscriptstyle22}}_{\{d,^{p}m\},d^{p}m,pn\partial_1m'n'}.
\end{align*}
We need to add a new relation:
$$\overline{\mathbf{v}^{\scriptscriptstyle22}}_{{\{m,m'm''\},mm'm'',n''}}=\overline{\mathbf{v}^{\scriptscriptstyle22}}_{{\{m,m'\},mm',\partial_1m''n''}}+\overline{\mathbf{v}^{\scriptscriptstyle22}}_{{\{m,m''\},mm'',n''}}$$
in $J_2$. By taking in this relation $(1,m,n)=(1,d,p)$ and \\ $(1,^{p}m,pn)=(1,m',n')$,\ \ $(1,^{pn}m',pnn')=(1,m'',n'')$ \\ we have $$\partial_1m''n''=\partial_1(^{pn}m')pnn'=pn\partial_1m'n' \text{  and  } n''=pnn'.$$
Using this relation, we obtain
$$\overline{\mathbf{v}^{\scriptscriptstyle22}}_{{\{d,^{p}m\},d^{p}m,pn\partial_1m'n'}}+\overline{\mathbf{v}^{\scriptscriptstyle22}}_{{\{d,^p(^{n}m')\},d^p(^{n}m'),pnn'}}=\overline{\mathbf{v}^{\scriptscriptstyle22}}_{{\{d,^{p}m^p(^{n}m')\},d^{p}m^p(^{n}m'),pnn'}}$$
and thus,
$$\lambda''_{\begin{bsmallmatrix} &\Gamma'\\ \Gamma& \end{bsmallmatrix}}={\begin{bmatrix} &{\lambda''_{\Gamma'}}\\ {\lambda''_{\Gamma}}& \end{bmatrix}}=K'_2 \lambda^1_{n}+\lambda^2_{\partial_1m'n'} H'_2.$$
Similarly, we must show that
$$\lambda''_{ \begin{bsmallmatrix} \Gamma&\\ &\Gamma' \end{bsmallmatrix}}={ \begin{bmatrix}{\lambda''_{\Gamma}}&\\ &{\lambda''_{\Gamma'}}\end{bmatrix}}=\lambda^2_{n'} H'_2+K'_2 \lambda'_{\partial_1mn}.$$
For any $\overline{\mathbf{v}^{\scriptscriptstyle11}}_{1,d,p}\in \mathrm{K}_2$, we obtain
$$\lambda''_{\begin{bsmallmatrix} \Gamma&\\ &\Gamma' \end{bsmallmatrix}}(\overline{\mathbf{v}^{\scriptscriptstyle11}}_{1,d,p})=\lambda''_{(^{\partial_1m}(^{n}m')m,nn')}(\overline{\mathbf{v}^{\scriptscriptstyle11}}_{1,d,p})
=\overline{\mathbf{v}^{\scriptscriptstyle22}}_{(\{d,^{p}(^{\partial_1m}(^{n}m')m)\} ,d^{p}(^{\partial_1m}(^{n}m')m),pnn')}$$
and ;
\begin{align*}
(\lambda^2_{n'}H'_2+K'_2\lambda'_{\partial_1mn})(\overline{\mathbf{v}^{\scriptscriptstyle11}}_{({1,d,p})})=&\lambda^2_{n'}(\lambda''_{m,n}(\overline{\mathbf{v}^{\scriptscriptstyle11}}_{({1,d,p})})
+\lambda''_{m',n'}(\lambda^1_{\partial_1mn}(\overline{\mathbf{v}^{\scriptscriptstyle11}}_{({1,d,p})})\\ =&\lambda^2_{n'}(\overline{\mathbf{v}^{\scriptscriptstyle22}}_{({\{d,^{p}m\},d^{p}m,pn})})+\lambda''_{(m'n')}(\overline{\mathbf{v}^{\scriptscriptstyle11}}_{1,d,p\partial_1mn})\\
=&\overline{\mathbf{v}^{\scriptscriptstyle22}}_{({\{d,^{p}m\},d^{p}m,pnn'})}+\overline{\mathbf{v}^{\scriptscriptstyle22}}_{(\{d,^{p\partial_1mn}m'\},d^{p\partial_1m}(^{n}m'),p\partial_1mnn')}
\end{align*}
and by taking
\begin{align*}
(1,m,n)=&(1,d,p)\\
(1,m',n')=&(1,^{p\partial_1m}(^{n}m'),p\partial_1mn)\\
(1,m'',n'')=& (1,^{p}m,pnn')\\
\partial_1(m'')n''=&\partial_1(^{p}m)pnn'=p\partial_1mnn'
\end{align*}
in the last relation, we have

\begin{align*}
{\begin{bmatrix}{\lambda''_{\Gamma}}&\\ &{\lambda''_{\Gamma'}}\end{bmatrix}}=&(\lambda^2_{n'} H'_2+K'_2 \lambda'_{\partial_1mn})(\overline{\mathbf{v}^{\scriptscriptstyle11}}_{1,d,p})\\
=&\overline{\mathbf{v}^{\scriptscriptstyle22}}_{({\{d,^{p}m\},d^{p}m,pnn'})}+\overline{\mathbf{v}^{\scriptscriptstyle22}}_{(\{d,^{p\partial_1m}(^{n}m')\},d ^{p\partial_1m}(^{n}m'),p\partial_1mnn')}\\
=&\overline{\mathbf{v}^{\scriptscriptstyle22}}_{({\{d,^{p}(^{\partial_1m}{(^{n}m'))^{p}m\}},d ^{p}(^{\partial_1m}{(^{n}m'))^{p}m},pnn'})} \  ( \text{ Since the last relation} )\\
=&\lambda''_{^{\partial_1m}(^{n}m')m,nn'}(\overline{\mathbf{v}^{\scriptscriptstyle11}}_{(1,d,p)})\\
=&\lambda''_{ \begin{bsmallmatrix} \Gamma&\\ &\Gamma' \end{bsmallmatrix}}
\end{align*}
$\Box$
\subsubsection{The proof of equality $\alpha' \overline \tau_2=K'_2-H'_2.$} \label{sect:App4}

For any 3-cell $(l,m,n)$ in $C_3$ and $(l',m',n')$ in $C_3$; a generator element in $J_2$ can be taken as
\begin{multline*}
{\mathbf{u}_{\scriptscriptstyle3}}={\mathbf{e}_{^{n'}l,^{n'}m,n'n}}-{\mathbf{e}_{1,^{n'}m,n'n}}-{\mathbf{e}_{^{\partial_1m'n'}l,{^{\partial_1m'n'}m,\partial_1m'n'n}}}
+{\mathbf{e}_{1,{^{\partial_1m'n'}m,\partial_1m'n'n}}}\\+{\mathbf{e}_{l,m^{n}m',nn'}}-{\mathbf{e}_{1,m^{n}m',nn'}}-{\mathbf{e}_{^{\partial_1mn}{(^{m'}l)},{^{\partial_1mn}m',nn'}}}+{\mathbf{e}_{1,{^{\partial_1mn}m',nn'}}}\in {J_2}
\end{multline*}
and
\begin{multline*}
{\mathbf{u}_{\scriptscriptstyle4}}={\mathbf{e}_{l{^{m}{(^n{l'})},m^{n}m'},nn'}}-{\mathbf{e}_{1,m^{n}m',nn'}}+{\mathbf{e}_{\{\partial_2lm,^{n}{\partial_2l'm'}\},\partial_2lm^n(\partial_2l'm'),nn'}}
-{\mathbf{e}_{1,\partial_2lm^n(\partial_2l'm'),nn'}}\\-{\mathbf{e}_{\{m,^{n}m'\},m{^{n}m'},nn'}}+{\mathbf{e}_{1,m{^{n}m'},nn'}}-{\mathbf{e}_{^{\partial_1mn}{l'(^{m'}l)},{^{\partial_1mn}m'm,nn'}}}+{\mathbf{e}_{1,{^{\partial_1mn}m'm,nn'}}}\in {J_2}
\end{multline*}
We need to add new relations given below for $J_1$.
\begin{multline*}
\mathbf{v}_{\scriptscriptstyle3}=\mathbf{e}_{1,^{n'}{(\partial_2lm)},n'n}-\mathbf{e}_{1,^{n'}m,n'n}-\mathbf{e}_{(1,{^{\partial_1m'n'}{(\partial_2lm)}},\partial_1m'n'n)}
+\mathbf{e}_{1,^{\partial_1m'n'}m,\partial_1m'n'n}\\
+\mathbf{e}_{1,\partial_2lm{^{n}{m'}},nn'}-\mathbf{e}_{1,m^nm',nn'}-\mathbf{e}_{1,{^{\partial_1mn}{(\partial_2(^{m'}l)m')}},nn'}+\mathbf{e}_{1,^{\partial_1mn}m',nn'}\in J_1
\end{multline*}
\begin{multline*}
\mathbf{v}_{\scriptscriptstyle4}=\mathbf{e}_{1,\partial_2lm{^{n}{(\partial_2l'm')}},nn'}-\mathbf{e}_{1,m^nm',nn'}+\mathbf{e}_{1,\partial_2\{\partial_2lm,^{n}{\partial_2l'm'}\}\partial_2lm{^{n}{(\partial_2l'm')}},nn'}
-\mathbf{e}_{1,\partial_2lm{^{n}{(\partial_2l'm')}},nn'}\\
-\mathbf{e}_{1,\partial_2\{m,^{n}{m'}\}m^nm',nn'}+ \mathbf{e}_{1,m^nm',nn'}-\mathbf{e}_{1,^{\partial_1mn}{\partial_2}(l')\partial_2(^{m'}l)^{\partial_1mn}{m'}m,nn'}
+\mathbf{e}_{1,^{\partial_1mn}{m'}m,nn'}\in J_1.
\end{multline*}
Using the generator elements $\mathbf{u}_{\scriptscriptstyle3}$ of $J_2$, we obtain the following relation in $\mathrm{K}_3$
$$\overline{\mathbf{v}^{\scriptscriptstyle22}}_{^{n'}l,^{n'}m,n'n}-\overline{\mathbf{v}^{\scriptscriptstyle22}}_{^{\partial_1m'}{(^{n'}l)},^{\partial_1m'n'}m,\partial_1m'n'n}=\overline{\mathbf{v}^{\scriptscriptstyle22}}_{^{\partial_1mn}{(^{m'}l)},^{\partial_1mn}m',nn'}
-\overline{\mathbf{v}^{\scriptscriptstyle22}}_{l,m^nm',nn'} \cdots(A)$$
Similarly the generator element $\mathbf{u}_{\scriptscriptstyle4}$ of $J_2$,we obtain the following relation in $\mathrm{K}_3$
\begin{multline*}
$$\overline{\mathbf{v}^{\scriptscriptstyle22}}_{l{^{m}{(^{n}{l'})},m^nm',nn'}}+\overline{\mathbf{v}^{\scriptscriptstyle22}}_{\{\partial_2lm,^{n}{(\partial_2l'm')}\},\partial_2lm{^{n}{(\partial_2l'm')}},nn'}\\
=\overline{\mathbf{v}^{\scriptscriptstyle22}}_{\{m,^{n}{m'}\},m^nm',nn'}+\overline{\mathbf{v}^{\scriptscriptstyle22}}_{^{\partial_1mn}{(l'^{m'}l)},^{\partial_1mn}{m'}m,nn'} \cdots(B)$$
\end{multline*}
Using the relations; we have
\begin{align*}
\alpha'_{l,m,n}\overline{\tau}_2+H'_2-K'_2=&\alpha'_{l,m,n}\overline{\tau}_2(\overline{\mathbf{v}^{\scriptscriptstyle11}}_{1,m',n'})+H'_2(\overline{\mathbf{v}^{\scriptscriptstyle11}}_{1,m',n'})-K'_2(\overline{\mathbf{v}^{\scriptscriptstyle11}}_{1,m',n'})\\ =&\alpha'_{l,m,n}(\mathbf{e}_{1,1,n'}-\mathbf{e}_{1,1,\partial_1m'n'})+\overline{\mathbf{v}^{\scriptscriptstyle22}}_{\{m,^{n}{m'}\},m^nm',nn'}-\overline{\mathbf{v}^{\scriptscriptstyle22}}_{\{\partial_2lm,^{n}{m'}\},\partial_2lm^{n}{m'},nn'}\\
=&\overline{\mathbf{v}^{\scriptscriptstyle22}}_{^{n'}l,^{n'}m,n'n}-\overline{\mathbf{v}^{\scriptscriptstyle22}}_{^{\partial_1m'n'}l,^{\partial_1m'n'}m,\partial_1m'n'n}\\
&+\overline{\mathbf{v}^{\scriptscriptstyle22}}_{\{m,^{n}{m'}\},m^nm',nn'}-\overline{\mathbf{v}^{\scriptscriptstyle22}}_{\{\partial_2lm,^{n}{m'}\},\partial_2lm^{n}{m'},nn'}
\end{align*}
In relation $(B)$, if we take $l'=1$, we have;
$$\overline{\mathbf{v}^{\scriptscriptstyle22}}_{l,m^nm',nn'}+\overline{\mathbf{v}^{\scriptscriptstyle22}}_{\{\partial_2lm,^{n}{m'}\},\partial_2lm{^{n}{m'}},nn'}
=\overline{\mathbf{v}^{\scriptscriptstyle22}}_{\{m,^{n}{m'}\},m^nm',nn'}+\overline{\mathbf{v}^{\scriptscriptstyle22}}_{^{\partial_1mn}{(^{m'}l)},^{\partial_1mn}{m'}m,nn'}$$
and from this equality we have
$$\alpha'_{l,m,n}\overline{\tau}_2+H'_2-K'_2=\overline{\mathbf{v}^{\scriptscriptstyle22}}_{^{n'}l,^{n'}m,n'n}-\overline{\mathbf{v}^{\scriptscriptstyle22}}_{^{\partial_1m'n'}l,^{\partial_1m'n'}m,\partial_1m'n'n}
+\overline{\mathbf{v}^{\scriptscriptstyle22}}_{l,m^nm',nn'}-\overline{\mathbf{v}^{\scriptscriptstyle22}}_{^{\partial_1m}{(^{n}{(^{m'}l)})},^{\partial_1mn}{m'},nn'}$$
and this is of course the relation $A$ in $\mathrm{K}_3$ . Thus, we obtain
$$\alpha'_{l,m,n}\overline{\tau}_2+H'_2-K'_2=\overline{0}.$$

\begin{equation*}
\begin{array}{llllll}
\text{Murat SARIKAYA} \text{ and }  \text{Erdal ULUALAN} \\
 \text{Dumlup\i nar University} \\
\text{Science Faculty}  \\
 \text{Department of Mathematics} \\
 \text{K\"{u}tahya-TURKEY}
\\
\text{E-mails}: \text{e.ulualan@dpu.edu.tr}\\
\text{murat.sarikaya@dpu.edu.tr}
\end{array}
\end{equation*}

\begin{thebibliography}{99}

\bibitem{Jinan} Al-asady, J. : \textrm{Representations of Crossed Squares and Cat$^2$-Groups}, University of Leicester, Dep. of Math. (Ph. D. Thesis) (2018).

\bibitem{Arvasi} Arvasi, Z.: \textrm{Crossed Squares and 2-Crossed Modules of Commutative Algebras }. Theory and Applications of Categories, 3(7):160-181, (1997).

\bibitem{AU1} Arvasi, Z. and Ulualan, E.: \textrm{On algebraic models for homotopy 3-types}. Journal of
Homotopy and Related Structures, Vol.1, No \textbf{1}, pp.1-27,
(2006).

\bibitem{Baez1} Baez, J.C. and  Crans, A.S.: \textrm{Higher dimensional algebra VI:
Lie 2-algebras}. {Theory and Applications of Categories}
vol.12, No.15, pp. 492-538, (2004).

\bibitem{BL} Baez, J.C. and Lauda, A.D.: \textrm{Higher dimensional algebra V:2-groups}.
 {Theory and Applications of Categories} vol.12, pp. 423-492,
(2004).

\bibitem{Barker} Barker, M.F. : \textrm{Representations of Crossed Modules and Cat$^1$-Groups}, University of Wales, Bangor (Ph. D. Thesis) (2003).

\bibitem{Be}  Berger, C.: \textrm{Double loop spaces, braided monoidal categories and algebraic 3-type of space}. Higher homotopy structures in topology and mathematical Physics (Poughkeepsie, NY, 1996), 49-66, Contemp. Math. 227, Amer. Math. Soc., Providence, RI, (1999).

\bibitem{Brown}  Brown, R.: \textrm{Possible connections between
whiskered categories and groupoids, Leibniz algebras, automorphism
structures and local-to-global questions}. {Journal of
Homotopy and Related Structures} Vol.5, No.1, pp. 305-318, (2010).

\bibitem{BG}  Brown, R. and Gilbert, N.D.: \textrm{ Algebraic
models of 3-types and automorphism structures for crossed
modules}. {Proc. London Math. Soc.}, (3) \textbf{59}, 51-73,
(1989).

\bibitem{Brown-loday}  Brown, R. and Loday, J.L.: \textrm{van
Kampen theorems for diagrams of spaces}. {Topology},
\textbf{26}, 311-335, (1987).

\bibitem{bs} Brown, R. and Spencer, C.B.: \textrm{$\mathcal{%
G}$-Groupoids, crossed modules and the fundamental groupoid of a
topological group}. {Proc. Konink. Neder. Akad. van
Wetenschappen Amsterdam},  Vol:79 (\textbf{4}) (1976).

\bibitem{Burrow} Burrow, M.: \textrm{Representation Theory of Finite Groups}.
{Academic Paperbacks},  (1965).

\bibitem{cc}  Carrasco, P. and  Cegarra, A.M.: \textrm{
Group-theoretic algebraic models for homotopy types}.
{Journal Pure Appl. Algebra}, \textbf{75}, 195-235, (1991).

\bibitem{Con}  Conduch{\'{e}}, D.: \textrm{Modules crois{\'{e}}s g{%
\'{e}}n{\'{e}}ralis{\'{e}}s de longueur 2}. {Jour. Pure and
Applied Algebra},  \textbf{34}, 155-178, (1984).

\bibitem{Con1} Conduch\'{e}, D.: \textrm{Simplicial crossed modules
and mapping cones}, \emph{Georgian Mathematical Journal},
\textbf{10}, No. 4, 623-636, (2003).

\bibitem{Cr} Crans, S.E.: \textrm{ A tensor product for Gray categories}, Theory Appl. Categ., 5, 12-69, (1999).

\bibitem{Curtis} Curtis, C.W. and Reiner, I.: \textrm{Representation Theory of Finite Groups and Associative Algebras},  Interscience, (1962).

\bibitem{Dornhoff} Dornhoff, L.: \textrm{Group Representation Theory},  Marcel Dekker, New York, (1971).


\bibitem{Elgueta}  Elgueta,J.: \textrm{Representation Theory of 2-groups on Kapranov and Voevodsky's 2-vector spaces}, Advances in Mathematics, 213, 53-92,2007.

\bibitem{Ellis} Ellis, G.J.: Crossed modules and their higher dimensional analogues,Ph. D. Thesis, UCNW, Bangor,  (1984).

\bibitem{Gray}Gray, J.W.: \textrm{Formal category theory: Adjointness for 2-categories}, Lecture Notes in Math. 391, Sipringer Berlin, (1974).

\bibitem{HKK} Hardie, K.A.; Kamps, K.H. and Kieboom, R.W.: \textrm{A Homotopy 2-groupoid of a Haussdorff Space. Papers in honour of Bernhard Banaschewski (Cape Town,1996). Appl. Categ. Structures 8 (2000), no. 1-2, 209-234. Representation Theory of 2-groups on Kapranov and Voevodsky's 2-vector spaces}, Advances in Mathematics, 213, 53-92,2007.

\bibitem{JT}Joyal, A. and Tierney, M.: \textrm{ Algebraic homotopy types, Handwritten lecture notes}, (1984).

\bibitem{KP}Kamps, K.H. and  Porter, T.: \textrm{2-groupoid enrichments in homotopy theory and algebra, K-Theory, 25, 4, 373-409}, (2002).

\bibitem{ppp} Kamps, K.H., Pumplün, D. and Tholen, W.(eds): \textrm{ Category Theory} (Gummersbach 1981), no.962 in Lecture Notes in Mathematics, Berlin, 1982, Springer-Verlag.

\bibitem{Kapranov} Kapranov, M. and  Voevodsky, V.: \textrm{2-categories and Zamelodchikov tetrahedra equations}, Proc. Symp. Pure math. \textbf{56}:177-260, (1994).

\bibitem{Kelly} Kelly, G.M. and Street, R.: \textrm{Review of the elements of 2-categories} in Sydney Category Seminar (1972/3), no.420 in Lecture Notes in Mathematics, Springer-Verlag, (1974).

\bibitem{Ledermann} Ledermann, W : \textrm{Intorduction to Group Characters}, Cambridge University Press,(1977).

\bibitem{Loday} Loday, J.-L.: \textrm{Spaces with Finitely Many Non-trivial Homotopy Groups}. {Journal of Pure and Applied Algebra}, \textbf{24} (2):179-202, (1982).



\bibitem{Martins} Martins, J.F. and Picken, R.: \textrm{The Fundamental Gray 3-groupoid of a smooth manifold and local 3-dimensional holonomy based on a 2-crossed module}. {Differential Geometry and its Applications}, \textbf{29} (2):179-206, (2011).



\bibitem{mutpor2}   Mutlu, A. and Porter, T.: \textrm{Applications of Peiffer pairings in the
Moore complexes of a simplicial group}. {Theory and
Applications of Categories}, Vol.4, No.\textbf{7},  pp.148-173,
(1998).


\bibitem{mutpor3}  Mutlu, A. and Porter, T.: \textrm{Freeness conditions for  crossed squares and squared complexes }. {K-Theory
}, 20,345–368 (2000).

\bibitem{Passman}  Passman, D.S.: \textrm{Infinite Group Rings }. {Marcel Dekker}, New York (1971).

\bibitem{Porter} Porter, T.: \textrm{Internal Categories and Crossed Modules}, in Kamps et al., \cite{ppp}, pp. 249-255.

\bibitem{Serre}  Serre, J.-P.: \textrm{Représentations Linéaries des Groupes Finis }. {Hermann}, Paris, (1971).

\bibitem{Nizar}  Shammu, N.M.: \textrm{Algebraic and Categorical Structure of Categories of Crossed Modules of Algebras }. {Ph. D. Thesis}, UCNW, Bangor (1992).


\bibitem{WL} Walery, D.G., Loday, J.L.: \textrm{Obstructions a l'excision en K-theorie algebrique}, Sipringer Lecture Notes in Math., 854, 179-216 (1981).

\bibitem{Wang} Wang, W.: \textrm{On 3-gauge transformations, 3-curvatures, and Gray-categories} J. Math. Phys. 55, 4, 043506 (2014).

\bibitem{W}  Whitehead, J.H.C.: \textrm{Combinatorial homotopy II}.
 {Bull. Amer. Math. Soc.}, \textbf{55}, 453-496, (1949).
\end{thebibliography}
\end{document}